\documentclass{article}
\usepackage[a4paper,top=4cm,bottom=4cm,left=2.5cm,right=2.5cm,bindingoffset=5mm]{geometry}
\usepackage{xcolor}

\usepackage[hyphens]{url}\usepackage[breaklinks=true,colorlinks,linkcolor={blue},citecolor={red},urlcolor={purple},]{hyperref}

\usepackage[font=small,format=default,indention=2em,labelsep = period]{caption} 
\usepackage[subrefformat=parens,labelformat=parens]{subfig}
\usepackage{enumitem}
\usepackage{tikz}\tikzset{>=latex}
\usepackage{overpic}

\usepackage{amssymb,bm,amsmath,amsthm,mathrsfs}
\usepackage{graphicx,algorithm,algorithmic,mathtools}	

\usepackage{enumitem}
\usepackage{todonotes}

\numberwithin{equation}{section}
\theoremstyle{definition}

\newtheorem{remark}{Remark}

\newcommand{\R}{\mathbb R}

\newcommand{\dt}{\Delta t}
\newcommand{\disp}{\displaystyle}
\newcommand\norm[1]{\left\lVert#1\right\rVert}

\def\rm{r^{\mathrm{max}}}

\def\cy{c^{y}}
\def\cx{c^{x}}
\def\ux{u^{x}}
\def\uy{u^{y}}
\def\rhot{\widetilde\rho}
\def\mut{\widetilde\mu}
\def\uxt{\widetilde u^{x}}
\def\uyt{\widetilde u^{y}}
\def\qx{q^{x}}
\def\qy{q^{y}}
\def\qxr{q^{x}_{\rho}}
\def\qxm{q^{x}_{\mu}}
\def\qyr{q^{y}_{\rho}}
\def\qym{q^{y}_{\mu}}
\def\uxr{u^{x}_{\rho}}
\def\uxm{u^{x}_{\mu}}
\def\uyr{u^{y}_{\rho}}
\def\uym{u^{y}_{\mu}}
\def\qxrt{\widetilde q^{x}_{\rho}}
\def\qxmt{\widetilde q^{x}_{\mu}}
\def\qyrt{\widetilde q^{y}_{\rho}}
\def\qymt{\widetilde q^{y}_{\mu}}
\def\lx{L^{x}}
\def\ly{L^{y}}
\def\nx{N_{x}}
\def\ny{N_{y}}
\def\rhob{\overline{\rho}}
\def\mub{\overline{\mu}}
\newcommand\rhotr[1]{\rho^{\mathrm{true #1}}}
\newcommand\mutr[1]{\mu^{\mathrm{true #1}}}
\def\esse{\mathscr{S}}

\def\kuno{\kappa_{1}}
\def\kdue{\kappa_{2}}
\def\vp{v_{+}}
\def\vm{v_{-}}
\def\vzero{v_{0}}
\def\fv{\hat{f}}
\def\gv{\hat{g}}
\def\lx{L^{x}}
\def\ly{L^{y}}
\def\nx{N_{x}}
\def\ny{N_{y}}
\def\deltax{\Delta x}
\def\deltay{\Delta y}
\def\deltat{\Delta t}
\def\omij{\Omega_{ij}}
\def\uij{U_{ij}}

\def\minute{\mathrm{min}}
\def\meter{\mathrm{m}}
\def\second{\mathrm{s}}
\def\vehkm{\mathrm{veh/km}}

\title{\LARGE\textbf{A two-dimensional multi-class traffic flow model}}
\author{\normalsize{Caterina Balzotti}\thanks{SBAI Department, Sapienza University of Rome (\href{mailto:caterina.balzotti@sbai.uniroma1.it}{caterina.balzotti@sbai.uniroma1.it})}
\and \normalsize{Simone Göttlich}\thanks{Department of Mathematics, University of Mannheim (\href{mailto:goettlich@uni-mannheim.de}{goettlich@uni-mannheim.de}).
}}
\date{\today}

\begin{document}
\maketitle

\begin{abstract}
\noindent
The aim of this work is to introduce a two-dimensional macroscopic traffic model for multiple populations of vehicles. Starting from the paper \cite{herty2018NHM}, where a two-dimensional model for a single class of vehicles is proposed, we extend the dynamics to a multi-class model leading to a coupled system of conservation laws in two space dimensions. Besides the study of the Riemann problems we also present a Lax-Friedrichs type discretization scheme recovering the theoretical results by means of numerical tests. We calibrate the multi-class model with real data and compare the fitted model to the real trajectories. Finally, we test the ability of the model to simulate the overtaking of vehicles. 
\end{abstract}

\begin{description}
\item[\textbf{Keywords.}] Macroscopic traffic flow, two-dimensional model, multi-class model, Riemann problems, data-fitting.
\item[\textbf{Mathematics Subject Classification.}]  90B20; 35L65; 35Q91.
\end{description}

\section{Introduction}
In this paper, we are concerned with the study of a two-dimensional multi-class traffic model. This work is placed in the constantly evolving framework of mathematical models for traffic flow. The goal of traffic models is to provide tools capable of helping traffic management, in order to optimize transport and obtain economic and environmental benefits, such as the reduction of vehicles queues and pollution.

Traffic models are divided into three main categories, which depend on the scale of observation: microscopic, macroscopic and kinetic models. Microscopic models follow the dynamics of each vehicle and are described by ordinary differential equations (ODEs), see e.~g.~\cite{aw2002SIAP, gazis1961OR, helbing1998springer,nelson1995TTSP, Newell1961,Pipes1953}. Macroscopic models, based on fluid dynamics, consider aggregated quantities such as the density of vehicles and are governed by partial differential equations (PDEs), see e.~g.~\cite{AwRascle2000,colombo2002MCM, DAGANZO2006TRB, goatin2006MCM, LighthillWhitham1955, Richards1956, Zhang2002}. Kinetic models \cite{IllnerKlarMaterne2003,klar2000SIAP,nelson1995TTSP,phillips1979TPT, prigogine1960OR,puppo2017KRM} are between the previous two classes since they can be derived by microscopic models while macroscopic models can be derived by kinetic descriptions.
We refer to \cite{albi2019M3AS,GaravelloHanPiccoli2016} and references therein for a more complete review on traffic models. 

In recent years, the ever-increasing amount of real data, due to new technologies, has widely influenced the research on mathematical models for traffic flow. The common goal of researchers is to exploit the available real data to build \emph{ad hoc} traffic models, capable of simulating increasingly realistic scenarios. We refer to \cite{colombo2016M3AM,CristianiDeFabritiisPiccoli2010,FanHertySeibold2013,PiccoliKeFrieszYao2012} for some inspiring examples of data-fitted traffic models. The common feature of these models is the application of vehicles trajectory data collected in datasets such as \cite{germanDataset,TrafficNGSIM}. Datasets of this type generally contain data on multi-lane highways and are able to distinguish the type of vehicle.

The focus in this work is on macroscopic traffic models. In particular, we propose a multi-class generalization in two space dimensions of the well-known first order Lighthill-Whitham-Richards (LWR) model \cite{LighthillWhitham1955,Richards1956}. First order models such as the LWR are described by a single conservation law $\rho_{t}+f(\rho)_{x} = 0,$ 
where $\rho$ is the density of vehicles and $f(\rho)$ is the flux function. The shortcomings of first order models are well-known in literature, for instance the infeasible solutions with unbounded acceleration \cite{lebacque2003TRB} or the inability to reproduce complex traffic phenomena like stop-and-go waves \cite{laval2010PTRSA,schonhof2007TS}. However, the extension of first order models to multi-class \cite{benzoni2003EJAM,herty2006CMS,wong2002TRS}, multi-lane \cite{holden2019SIAM} or multi-dimensional \cite{herty2018NHM} models has proven to be suitable to improve the deficiencies of the LWR model and to be able to describe also complex traffic phenomena.

As we have already mentioned, traffic datasets contain information related to multi-lane highways with different types of vehicles. Most traffic models refer to dynamics of single-lane traffic and therefore do not consider the movements related to lane changes. Our aim is to exploit now all the available data, including the line-changing behavior and the different vehicles classes. To this end, we propose an extension of the work by Herty, Fazekas and Visconti for a single-class traffic model~\cite{herty2018NHM} to a multi-class traffic model in two space dimensions.
The most common approaches which include lane-changing are the two-dimensional models and the multi-lane models. The first approach is an emerging topic, and we refer to \cite{chetverushkin2006,herty2018NHM,herty2018SIAM} for some examples. The second approach has been used for instance in  \cite{klar1999SIAM1,klar1999SIAM2}, where the authors propose a microscopic, a kinetic and a fluid dynamic model with lane changing. Here, we stick to the two-dimensional approach and incorporate two types of vehicles interacting through the flux functions.
%
%
The proposed model is then defined by the coupling of LWR-type models for two classes of vehicles in the $x$ and $y$ direction. The interaction between the two classes of vehicles is modeled by means of the flux functions which depend on the sum of vehicle densities as in \cite{daganzo1997TRB,herty2006CMS}. 
With suitable assumptions on the flux functions, we study the two-dimensional Riemann problems and validate the model comparing the theoretical results with the solutions given by a numerical approximation of Lax-Friedrichs type. Then, we calibrate the flux and velocity functions with the German dataset \cite{germanDataset} and  compare the results of our model with real trajectories data. We also test the ability of the model of capturing vehicles overtaking.

The paper is organized as follows. In Section \ref{sec:modello} we introduce the traffic model and study the Riemann problems. In Section \ref{sec:numerica}, we describe the numerical scheme and validate the model via numerical tests. In Section \ref{sec:dati}, we calibrate the model with a German dataset and compare the results with the real trajectories of vehicles. In Section \ref{sec:sorpasso}, we propose a modified version of the model calibrated with real data, and finally we investigate on the ability of the model to simulate vehicles overtaking compared to a multi-lane model.

\section{Two-dimensional multi-class model}\label{sec:modello}

In this section, we introduce the traffic model used throughout the paper.
Let us consider two classes of vehicles, whose densities are denoted by $\rho$ and $\mu$, respectively. Our aim is to describe the dynamics of the two classes by means of a two-dimensional multi-class model. To this end, following \cite{herty2018NHM}, we introduce a LWR-type model in two dimensions for two classes of vehicles, i.e., 
\begin{equation}
\begin{cases}
\rho_{t}+(\qxr(\rho,\mu))_{x}+(\qyr(\rho,\mu))_{y} = 0\\
\mu_{t}+(\qxm(\rho,\mu))_{x}+(\qym(\rho,\mu))_{y} = 0,\\
\end{cases}
\label{eq:lwr2D}
\end{equation}
where $\qx_{\rho,\mu}$ are the fluxes of $\rho$ and $\mu$ along the $x$-direction, and $\qy_{\rho,\mu}$ are the fluxes of $\rho$ and $\mu$ along the $y$-direction. Similarly to \cite{herty2006CMS}, we define the flux functions as
\begin{equation}
\begin{split}
\qxr(\rho,\mu)&=\rho \cx\left(1-\left(\frac{\rho+\mu}{\rm}\right)\right) \qquad 
\qyr(\rho,\mu)=\rho \cy\left(1-\left(\frac{\rho+\mu}{\rm}\right)\right)\\
\qxm(\rho,\mu)&=\mu \cx\left(1-\left(\frac{\rho+\mu}{\rm}\right)\right) \qquad 
\qyr(\rho,\mu)=\mu \cy\left(1-\left(\frac{\rho+\mu}{\rm}\right)\right),
\label{eq:flussi}
\end{split} 
\end{equation}
where $\cx$ and $\cy$ are parameters to be calibrated and $\rm$ is the maximum density of vehicles. 
The velocity functions in $x$ and $y$ directions coincide for $\rho$ and $\mu$, and are defined by 
%
\[	\ux = \cx\left(1-\left(\frac{\rho+\mu}{\rm}\right)\right), \quad \uy = \cy\left(1-\left(\frac{\rho+\mu}{\rm}\right)\right).
\]
Hence, $\cx$ and $\cy$ represent the maximum velocity in $x$ and $y$ direction.
Note that we assume that the two classes of vehicles have the same velocity $\cx$ and $\cy$, and they have the same maximum density $\rm$.

First of all, we present the properties of model \eqref{eq:lwr2D}. To simplify the notation, we normalize $\rho$ and $\mu$ in order to fix $\rm=1$. We introduce the following vectors 
\allowdisplaybreaks
\[U = \begin{pmatrix}\rho\\\mu\end{pmatrix},\quad
f(U) = \begin{pmatrix}
\rho\cx\left(1-(\rho+\mu)\right)\\[0.1em]
\mu\cx\left(1-(\rho+\mu)\right)
\end{pmatrix}, \quad
g(U) = \begin{pmatrix}
\rho\cy\left(1-(\rho+\mu)\right)\\[0.1em]
\mu\cy\left(1-(\rho+\mu)\right)
\end{pmatrix}
\]
and matrices
\[A(U) = Df(U), \quad
B(U) = Dg(U).\]
%
Therefore, we can rewrite system \eqref{eq:lwr2D} as 
\begin{equation}
U_{t}+AU_{x}+BU_{y} = 0.
\label{eq:2Dvett}
\end{equation}
System \eqref{eq:2Dvett} is hyperbolic if any linear combination of $A$ and $B$ is diagonalizable. Thus, for $(\kuno,\kdue) \in \R^{2}$, we define $C=\kuno A + \kdue B$. The eigenvalues of $C$ are

\begin{align*}
\lambda_{1} =(\kuno\cx+\kdue\cy)(1-(\rho+\mu)), \qquad
\lambda_{2} = (\kuno\cx+\kdue\cy)(1-2(\rho+\mu))
\end{align*}  
which are real for any couple $(\kuno,\kdue)$, and they coincide if and only if $(\rho,\mu)= (0,0)$ or $\kuno= -\cy\kuno/\cx$. The associated eigenvectors are
%
\[\gamma_{1} = \begin{pmatrix}
-1\\1
\end{pmatrix}, \quad
\gamma_{2} = \begin{pmatrix}
\rho/\mu\\1
\end{pmatrix}.
\]
The first eigenvalue is linearly degenerate, i.e., $\nabla\lambda_{1}\cdot \gamma_{1} = 0$, while the second one is genuinely nonlinear, i.e., $\nabla\lambda_{2}\cdot \gamma_{2} \neq 0$. The Riemann invariants are 
%
\[	z_{1} = \rho+\mu, \quad
	z_{2} = \log\left(\rho/\mu\right).
\]

\subsection{Two-dimensional Riemann problems}\label{sec:riemann}

Next, we describe the two-dimensional Riemann problem \cite{dafermos2016springer,zhang1989TAMS} associated with \eqref{eq:lwr2D}. To simplify the computations, we introduce the variable $r=\rho+\mu$, so that problem \eqref{eq:lwr2D} can be rewritten as
\begin{equation}\label{eq:r}
\begin{cases}
r_{t}+\left(r\cx\left(1-r\right)\right)_{x}+\left(r\cy\left(1-r\right)\right)_{y} = 0\\
r(0,x,y) = r_{0}(x,y).
\end{cases}
\end{equation} 
The Riemann problem in two dimensions is given by \eqref{eq:r} with constant initial data in the four quadrants, i.e. 
\begin{equation}\label{eq:RPinit}
r_{0}(x,y) = \begin{cases}
v_{1} &\quad 0<x<\infty,\quad 0<y<\infty\\
v_{2} &\quad -\infty<x<0,\quad 0<y<\infty\\
v_{3} &\quad -\infty<x<0,\quad -\infty<y<0\\
v_{4} &\quad 0<x<\infty,\quad -\infty<y<0.\\
\end{cases}
\end{equation}
For convenience, we define $\fv(r)=r\cx\left(1-r\right)$ and $\gv(r) = r\cy\left(1-r\right)$.
\begin{remark}\label{rem:negative}
The treatment of the two-dimensional Riemann problem proposed in \cite{dafermos2016springer} assumes convex flux functions $f$ and $g$. In order to recover this hypothesis in our case, it is sufficient to choose the parameters $\cx$ and $\cy$ negative. However, the concave case can be recovered from the following analysis through proper sign changes.
\end{remark}
We look for self-similar solutions $r(t,x,y) = v(x/t,y/t)$ and therefore introduce $\xi = \tfrac{x}{t}$ and $\eta=\tfrac{y}{t}.$
We can rewrite the first equation of \eqref{eq:r} as
\begin{equation}
(\cx(1-2v)-\xi)v_{\xi}+\left(\cy(1-2v)-\eta\right)v_{\eta}=0
\label{eq:selfsimilar}
\end{equation}
which leads us to
%
\[(\cx(1-2v)-\xi)d\eta+\left(\cy(1-2v)-\eta \right)d\xi=0,
\]
where $\cx(1-2v) = \fv'(v)$ and $\cy(1-2v)=\gv'(v)$.

The set of singular points parametrized by $v$ is the straight line
%
\[	\esse = \left\{(\xi,\eta) \,|\, \xi =\cx(1-2v), \,\eta = \cy(1-2v)\right\}.
\]
Defining
\begin{equation}
\begin{split}
\gamma(\vm,\vp) &= \frac{\fv(\vp)-\fv(\vm)}{\vp-\vm}=\cx(1-\vp-\vm) \\ 
\nu(\vm,\vp) &= \frac{\gv(\vp)-\gv(\vm)}{\vp-\vm}=\cy(1-\vp-\vm),
\end{split}
\label{eq:gamma}
\end{equation}
the Rankine-Hugoniot jump condition is
\begin{equation}
\frac{d\eta}{d\xi} =-\frac{\nu(\vm,\vp)-\eta}{\gamma(\vm,\vp)-\xi}= -\frac{\cy(1-\vp-\vm)-\eta}{\cx(1-\vp-\vm)-\xi}.
\label{eq:RH}
\end{equation}
Assuming that the normal vector $(d\eta,d\xi)$ is directed towards the positive side of the shock curves, the Oleinik's entropy condition is
\begin{equation}
\begin{split}
&\left(\gamma(\vm,\vzero)-\gamma(\vm,\vp)\right)d\eta+\left(\nu(\vm,\vzero)-\nu(\vm,\vp)\right)d\xi  \\&=\,\cx(1-\vp-\vzero)d\eta+\cy(1-\vp-\vzero)d\xi\\
&\geq 0
\end{split}
\label{eq:oleinik}
\end{equation}
for $\vzero$ between $\vm$ and $\vp$.

\begin{figure}[h!]
\centering
\subfloat[][Case of no shocks.]
{\label{fig:caso1}
\begin{tikzpicture}[scale=0.8]
\draw [very thick](1.2,1.2) -- (4.8,3.8);
\node at (3.4,3.1) {$\mathscr{S}$};
\draw (0,1.2) -- (1.2,1.2);
\draw (1.2,1.2) -- (1.2,0);
\draw (2.4,0) -- (2.4,5);
\draw (0,2.07) -- (2.4,2.07);
\draw (1.8,0) -- (1.8,1.65);
\draw (0,1.65) -- (1.8,1.65);
\draw(3,0) -- (3,5);
\draw(3,2.5) -- (6,2.5);
\draw(4,3.22) -- (6,3.22);
\draw(4,3.2) -- (4,5);
\draw(4.8,3.8) -- (6,3.8);
\draw(4.8,3.8) -- (4.8,5);
\draw[->] (0,4) -- node[above = 0.5] {\small$\eta$}(0,5);
\draw[->] (0,4) -- node[right = 0.5] {\small$\xi$}(1,4);
\node at (0.8,0.8) {\small$v_{3}$};
\node at (2,2.4) {\small$v_{2}$};
\node at (5.2,4.2) {\small$v_{1}$};
\node at (3.4,2) {\small$v_{4}$};
\end{tikzpicture}
} \qquad\qquad
\subfloat[][Case of no rarefaction waves.]
{\label{fig:caso2}
\begin{tikzpicture}[scale=0.8]
\draw[->] (0,4) -- node[above = 0.5] {\small$\eta$}(0,5);
\draw[->] (0,4) -- node[right = 0.5] {\small$\xi$}(1,4);
\draw[fill=black] (3,4) circle (0.05) node[above right] {$A$};
\draw (3,4) -- (3,5);
\draw (3,4) -- (2,4);
\draw (3,4) -- (3.4,2.6);
\draw[fill=black] (3.4,2.6) circle (0.05) node[below left] {$O$};
\draw (3.4,2.6) -- (5,1);
\draw[fill=black] (5,1) circle (0.05) node[above right] {$B$};
\draw (5,1) -- (5,0);
\draw (5,1) -- (6,1);
\node at (2,1) {\small$v_{3}$};
\node at (2.5,4.5) {\small$v_{2}$};
\node at (5,4) {\small$v_{1}$};
\node at (5.5,0.5) {\small$v_{4}$};
\end{tikzpicture}
}
\caption{Representation of no shocks \protect\subref{fig:caso1} and no rarefaction waves \protect\subref{fig:caso2}.}
\end{figure}
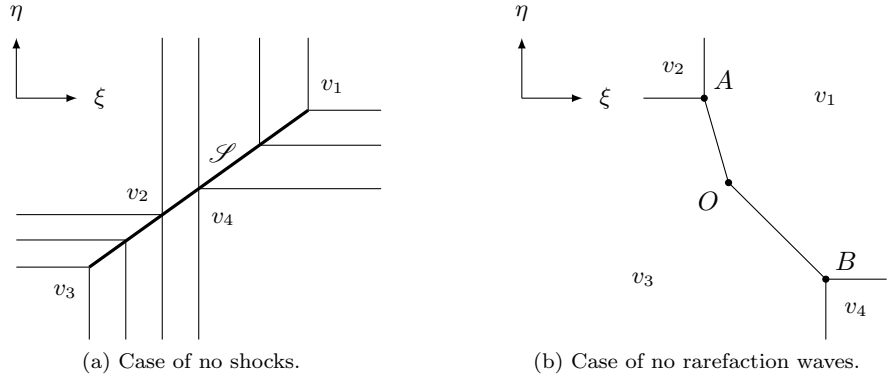

The initial data in \eqref{eq:RPinit} for problem \eqref{eq:r} in the variables $(\xi,\eta)$ is given by  
\begin{equation}
\lim_{\substack{\xi/\eta=const,\\\xi^{2}+\eta^{2}\to\infty}}v(\xi,\eta) = \begin{cases}
v_{1} &\quad \xi>0,\quad\eta>0\\
v_{2} &\quad \xi<0,\quad\eta>0\\
v_{3} &\quad \xi<0,\quad\eta<0\\
v_{4} &\quad \xi>0,\quad\eta<0.
\label{eq:datoIniz}
\end{cases}
\end{equation}
The solution of problem \eqref{eq:selfsimilar} with initial data \eqref{eq:datoIniz} is composed of elementary waves. There are five possible cases:
(1) no shocks, (2) no rarefaction waves, (3) exactly one shock, (4) exactly one rarefaction wave, (5) two rarefaction waves and two shocks.
In this work we skip the full details of the possible cases and refer to \cite{zhang1989TAMS} for a detailed discussion. Let us highlight the five cases now:\\[1ex]
\noindent
{\bf (1) No shocks:} This case is verified when $v_{3}<v_{2}< v_{4}<v_{1}$. Each couple $(v_{2},v_{1}), \,(v_{3},v_{4})\, (v_{1},v_{4})$ and $(v_{2},v_{3})$ is connected by rarefaction waves and the straight line $\esse$ defines the points of connection between them. The solution is represented in Figure \subref*{fig:caso1}.\\[1ex]
\noindent
{\bf (2) No rarefaction waves:} This case is verified when $v_{3}>v_{4}> v_{2}>v_{1}$. The couples $(v_{2},v_{1})$ and $(v_{2},v_{3})$ are connected by two shocks which collide in $A = (\gamma(v_{1},v_{2}),\nu(v_{2},v_{3}))$ while the couples $(v_{3},v_{4})$ and $(v_{4},v_{1})$ are connected by two shocks colliding in $B = (\gamma(v_{3},v_{4}),\nu(v_{1},v_{4}))$. Then, we have two shocks which connect $v_{1}$ and $v_{3}$. They start from the point $O = (\gamma(v_{1},v_{3}),\nu(v_{1},v_{3}))$ and terminate either in $A$ or in $B$. The solution is represented in Figure \subref*{fig:caso2}. 

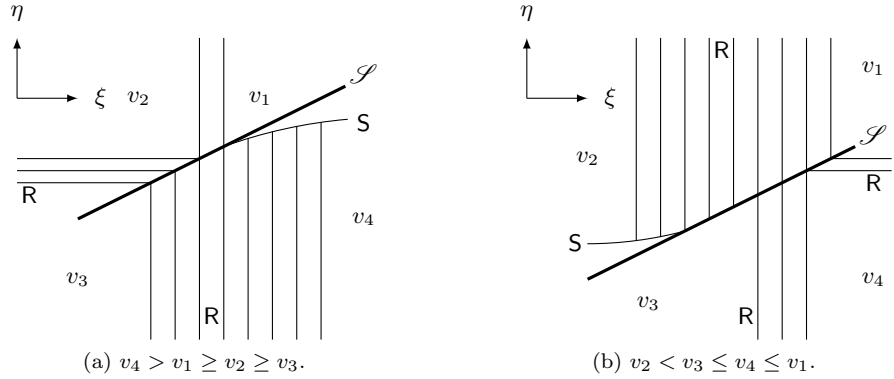
\begin{figure}[h!]
\centering
\subfloat[][$v_{4}>v_{1}\geq v_{2}\geq v_{3}$.]
{\label{fig:casoc1}
\begin{tikzpicture}[scale=0.8]
\draw[->] (0,4) -- node[above = 0.5] {\small$\eta$}(0,5);
\draw[->] (0,4) -- node[right = 0.5] {\small$\xi$}(1,4);
\draw[very thick] (1,2) -- (5.4,4.2);
\draw (3,0) -- (3,5);
\draw (3.4,0) -- (3.4,5);
\draw(0,3) -- (3,3);
\draw (2.6,0) -- (2.6,2.8);
\draw(0,2.8) -- (2.6,2.8);
\draw (2.2,0) -- (2.2,2.6);
\draw(0,2.6) -- (2.2,2.6);
\draw[black] (3.4,3.2) arc (110:95:8);
\draw(3.8,0) -- (3.8,3.34);
\draw(4.2,0) -- (4.2,3.45);
\draw(4.6,0) -- (4.6,3.54);
\draw(5,0) -- (5,3.6);
\node at (0.2, 2.4) {\small$\mathsf{R}$};
\node at (3.2, 0.4) {\small$\mathsf{R}$};
\node at (5.7, 3.6) {\small$\mathsf{S}$};
\node at (5.7, 4.4) {$\mathscr{S}$};
\node at (1,1) {\small$v_{3}$};
\node at (2,4) {\small$v_{2}$};
\node at (4,4) {\small$v_{1}$};
\node at (5.7,2) {\small$v_{4}$};
\end{tikzpicture}
} \qquad\qquad
\subfloat[][$v_{2}<v_{3}\leq v_{4}\leq v_{1}$.]
{\label{fig:casoc2}
\begin{tikzpicture}[scale=0.8]
\draw[->] (0,4) -- node[above = 0.5] {\small$\eta$}(0,5);
\draw[->] (0,4) -- node[right = 0.5] {\small$\xi$}(1,4);
\draw[very thick] (1,1) -- (5.4,3.2);
\draw (2.6,1.8) -- (2.6,5);
\draw[black] (2.6,1.8) arc (285:270:6.2);
\draw (2.2,1.7) -- (2.2,5);
\draw (1.8,1.64) -- (1.8,5);
\draw (3,2) -- (3,5);
\draw (3.4,2.2) -- (3.4,5);
\draw (3.8,0) -- (3.8,5);
\draw (4.2,0) -- (4.2,5);
\draw (4.6,0) -- (4.6,5);
\draw (4.6,2.8) -- (6,2.8);
\draw (5,3) -- (6,3);
\draw (5,3) -- (5,5);
\node at (3.6, 0.4) {\small$\mathsf{R}$};
\node at (5.7, 2.6) {\small$\mathsf{R}$};
\node at (3.2, 4.8) {\small$\mathsf{R}$};
\node at (0.8, 1.6) {\small$\mathsf{S}$};
\node at (5.7, 3.4) {$\mathscr{S}$};
\node at (2,0.6) {\small$v_{3}$};
\node at (1,3) {\small$v_{2}$};
\node at (5.7,4.5) {\small$v_{1}$};
\node at (5.7,1) {\small$v_{4}$};
\end{tikzpicture}
}
\caption{Representation of exactly one shock, where $\mathsf{R}$ denotes rarefaction waves and $\mathsf{S}$ shocks.}
\end{figure}

\noindent
{\bf (3) Exactly one shock:} This case is verified when
\[v_{4}>v_{1}\geq v_{2}\geq v_{3} \quad \text{or} \quad
v_{2}<v_{3}\leq v_{4}\leq v_{1}.
\]
The first sub-case is represented in Figure \subref*{fig:casoc1}. Using the Rankine-Hugoniot condition \eqref{eq:RH} it can be shown that the shock curve is concave, monotonically increasing in $(v_{1},v_{4})$, bounded by the base curve $\esse$, tangentially intersects $\esse$ and satisfies the entropy condition \eqref{eq:oleinik}. A similar analysis holds for the second sub-case which is represented in Figure \subref*{fig:casoc2}.

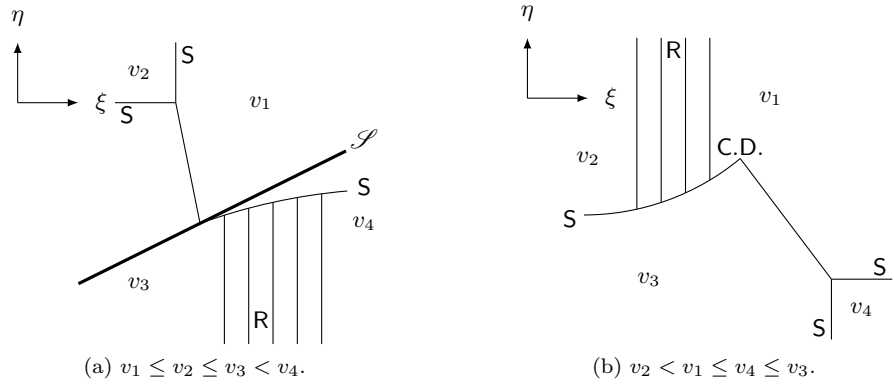
\begin{figure}[h!]
\centering
\subfloat[][$v_{1}\leq v_{2}\leq v_{3}< v_{4}$.]
{\label{fig:caso4a}
\begin{tikzpicture}[scale=0.8]
\draw[->] (0,4) -- node[above = 0.5] {\small$\eta$}(0,5);
\draw[->] (0,4) -- node[right = 0.5] {\small$\xi$}(1,4);
\draw[very thick] (1,1) -- (5.4,3.2);
\draw(2.6,4) -- (3,2);
\draw(2.6,4) -- (1.6,4);
\draw(2.6,4) -- (2.6,5);
\draw[black] (3,2) arc (110:95:9.5);
\draw(3.4,0) -- (3.4,2.14);
\draw(3.8,0) -- (3.8,2.25);
\draw(4.2,0) -- (4.2,2.35);
\draw(4.6,0) -- (4.6,2.43);
\draw(5,0) -- (5,2.49);
\node at (4, 0.4) {\small$\mathsf{R}$};
\node at (1.8, 3.8) {\small$\mathsf{S}$};
\node at (2.8, 4.8) {\small$\mathsf{S}$};
\node at (5.7, 2.6) {\small$\mathsf{S}$};
\node at (5.7, 3.4) {$\mathscr{S}$};
\node at (2,1) {\small$v_{3}$};
\node at (2,4.5) {\small$v_{2}$};
\node at (4,4) {\small$v_{1}$};
\node at (5.7,2) {\small$v_{4}$};
\end{tikzpicture}
} \qquad\qquad
\subfloat[][$v_{2}<v_{1}\leq v_{4}\leq v_{3}$.]
{\label{fig:caso4b}
\begin{tikzpicture}[scale=0.8]
\draw[->] (0,4) -- node[above = 0.5] {\small$\eta$}(0,5);
\draw[->] (0,4) -- node[right = 0.5] {\small$\xi$}(1,4);
\draw (1.8,2.16) -- (1.8,5);
\draw (2.2,2.28) -- (2.2,5);
\draw (2.6,2.43) -- (2.6,5);
\draw (3,2.64) -- (3,5);
\draw[black] (3.5,3) arc (310:270:4);
\draw (3.5,3) -- (5,1);
\draw (5,1) -- (6,1);
\draw (5,1) -- (5,0);
\node at (3.5, 3.2) {\small$\mathsf{C.D.}$};
\node at (2.4, 4.8) {\small$\mathsf{R}$};
\node at (0.7, 2) {\small$\mathsf{S}$};
\node at (4.8,0.2) {\small$\mathsf{S}$};
\node at (5.8, 1.2) {\small$\mathsf{S}$};
\node at (2,1) {\small$v_{3}$};
\node at (1,3) {\small$v_{2}$};
\node at (4,4) {\small$v_{1}$};
\node at (5.5,0.5) {\small$v_{4}$};
\end{tikzpicture}
}
\caption{Representation of exactly one rarefaction wave, where $\mathsf{R}$ denotes rarefaction waves, $\mathsf{S}$ shocks and $\mathsf{C.D.}$ contact discontinuities.}
\label{fig:caso4}
\end{figure}

\noindent
{\bf (4) Exactly one rarefaction wave:}
This case is verified when
\[v_{1}\leq v_{2}\leq v_{3}< v_{4} \quad \text{or} \quad
v_{2}<v_{1}\leq v_{4}\leq v_{3}.
\]
The first possibility of initial data gives results similar to the previous case of exactly one shock wave. For the second initial datum there exist several sub-cases but we omit the details. In Figure \ref{fig:caso4} we show two examples for the two initial data configurations.\\[1ex]

\begin{figure}[h!]
\centering
\subfloat[][$v_{4}> v_{1}\geq v_{3}> v_{2}$.]
{\label{fig:caso5a}
\begin{tikzpicture}[scale=0.8]
\draw[->] (0,4) -- node[above = 0.5] {\small$\eta$}(0,5);
\draw[->] (0,4) -- node[right = 0.5] {\small$\xi$}(1,4);
\draw[very thick] (2,2) --  node[below]{$\mathscr{S}$}(4,3);
\draw (1.6,1.7) -- (1.6,5);
\draw[black] (4,3) arc (110:95:7.8);
\draw (2,2) -- (2,5);
\draw (2.4,2.2) -- (2.4,5);
\draw (2.8,2.4) -- (2.8,5);
\draw (3.2,2.6) -- (3.2,5);
\draw (3.6,0) -- (3.6,5);
\draw (4,0) -- (4,5);
\draw (4.4,0) -- (4.4,3.14);
\draw (4.8,0) -- (4.8,3.24);
\draw (5.2,0) -- (5.2,3.34);
\draw[black] (2,2) arc (310:285:4);
\node at (3, 4.8) {\small$\mathsf{R}$};
\node at (4.2, 0.2) {\small$\mathsf{R}$};
\node at (0.3, 1.1) {\small$\mathsf{S}$};
\node at (5.8, 3.6) {\small$\mathsf{S}$};
\node at (2,1) {\small$v_{3}$};
\node at (1,3) {\small$v_{2}$};
\node at (5,4.5) {\small$v_{1}$};
\node at (5.8,1.5) {\small$v_{4}$};
\end{tikzpicture}
} \qquad\qquad
\subfloat[][$v_{4}> v_{3}> v_{1}> v_{2}$.]
{\label{fig:caso5b}
\begin{tikzpicture}[scale=0.8]
\draw[->] (0,4) -- node[above = 0.5] {\small$\eta$}(0,5);
\draw[->] (0,4) -- node[right = 0.5] {\small$\xi$}(1,4);
\draw[very thick] (1.5,1.5) --  (5,4);
\draw (1.8,1.87) -- (1.8,5);
\draw (2.2,2) -- (2.2,5);
\draw (2.6,2.3) -- (2.6,5);
\draw (3,2.57) -- (3,5);
\draw (3.4,2.85) -- (3.4,5);
\draw[black] (2.2,2) arc (290:270:5);
\draw[black] (3,2.58) arc (180:200:7.5);
\draw (3,2.57) -- (6,2.57);
\draw (3.4,2.85) -- (6,2.85);
\draw (3.01,2.25) -- (6,2.25);
\draw (3.02,1.95) -- (6,1.95);
\node at (2.4, 4.8) {\small$\mathsf{R}$};
\node at (5.8, 2.4) {\small$\mathsf{R}$};
\node at (0.5, 1.5) {\small$\mathsf{S}$};
\node at (3.6, 0.2) {\small$\mathsf{S}$};
\node at (5.2, 4.2) {$\mathscr{S}$};
\node at (2,1) {\small$v_{3}$};
\node at (1,3) {\small$v_{2}$};
\node at (4,4.5) {\small$v_{1}$};
\node at (5,1) {\small$v_{4}$};
\end{tikzpicture}
}
\caption{Representation of two shocks and two rarefaction waves, where $\mathsf{R}$ denotes rarefaction waves and $\mathsf{S}$ shocks.}
\label{fig:caso5}
\end{figure}
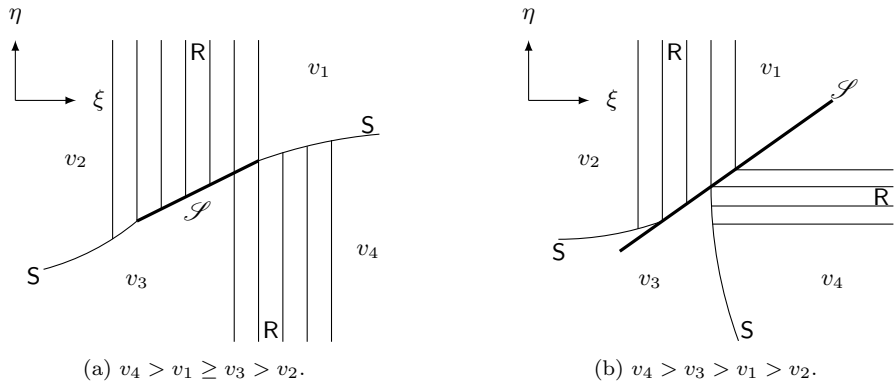

\noindent
{\bf (5) Two shocks and two rarefaction waves:}
This case is verified when
\[v_{4}> v_{1}\geq v_{3}> v_{2} \quad \text{or} \quad
v_{4}> v_{3}> v_{1}> v_{2}.
\]
The main difference between the two options of initial data is that in the first case the shock curves are not neighbors while in the second case they are neighbors. 
There are again several sub-cases, we only show an example of the two possible initial data sets in Figure \ref{fig:caso5}.

In the next section, we introduce a suitable discretization for the multi-class model. For validation purposes of the proposed scheme,
we aim to recover the theoretical results from above. 

\section{Numerical discretization}\label{sec:numerica}

The numerical analysis of system \eqref{eq:2Dvett} can be done using the dimensional splitting approach \cite{leveque2002CTAM} which means that the two-dimensional problem is split into two one-dimensional problems. Thus, equation \eqref{eq:2Dvett} is split into
\[U_{t}+AU_{x} = 0, \qquad
U_{t}+BU_{y}  = 0.
\]
We uniformly divide the two-dimensional space $[0,\lx]\times[0,\ly]$ into a grid $\Omega=[0,\nx]\times[0,\ny]$ with $x$-step $\deltax$ and $y$-step $\deltay$. We refer to the cell of the grid as $\omij$. Defining $\lambda_{1,2}$ and $\gamma_{1,2}$ the eigenvalues of $A$ and $B$ respectively, the time step $\deltat$ is determined by
\begin{equation}
\frac{\deltat}{\deltax}\leq \frac{1}{2}\Big(\max_{i,j=1,2}\{|\lambda_{i}|,|\gamma_{j}|\}\Big)^{-1}.
\label{eq:CFL}
\end{equation}
Then, the time interval $[0,T]$ is divided into time steps of length $\deltat$. 

Starting from a given initial datum $U^{0}_{ij}$, the numerical scheme is defined by the Strang splitting as
\[\begin{split}
\uij^{*} &= \uij^{n}-\frac{\deltat}{2\deltax}(F^{n}_{i+1/2,j}-F^{n}_{i-1/2,j})\smallskip\\
\uij^{**} &= \uij^{*}-\frac{\deltat}{\deltay}(G^{*}_{i,j+1/2}-G^{*}_{i,j-1/2})\smallskip\\
\uij^{n+1} &= \uij^{n}-\frac{\deltat}{2\deltax}(F^{**}_{i+1/2,j}-F^{**}_{i-1/2,j}).
\end{split}\]
%
We use the Local Lax-Friedrichs flux (also known as Rusanov flux) \cite[Chapter 3]{bertoluzza2009ACM} for $F$ and $G$, i.e.,

\begin{equation*}
F_{i+1/2,j} = \frac{1}{2}(f(U_{i+1,j})+f(U_{i,j})-\alpha_{i+1/2,j}(U_{i+1,j}-U_{i,j})),
\end{equation*}
where $\disp\alpha_{i+1/2,j}$ is the maximum modulus of the eigenvalues of the
Jacobian matrix in the interval $(U_{i,j},U_{i+1,j})$.

\subsection{Validation}
We now test the discretization method for the two-dimensional multi-class model \eqref{eq:2Dvett}
while comparing the numerical results to the theoretical solutions of the Riemann problems introduced in Section \ref{sec:riemann}. 

Our test setting is given by $\Omega=[-5,5]\times[-5,5]$ with $\deltax=\deltay=0.02$. We fix the parameters of \eqref{eq:flussi} to $\cx=\cy=-1$. As already observed in Remark \ref{rem:negative}, we fix negative parameters to recover convex flux functions.
The time interval $[0,T]=[0,1]$ is divided into time steps of length $\deltat$ satisfying condition \eqref{eq:CFL}. The initial datum for the two classes $\rho$ and $\mu$ is taken as in \eqref{eq:datoIniz}
\[\rho_{0}(x,y) = \begin{cases}
\rho_{1} &\quad (x,y)\in(0,5)\times(0,5)\\
\rho_{2} &\quad (x,y)\in(-5,0)\times(0,5)\\
\rho_{3} &\quad (x,y)\in(-5,0)\times(-5,0)\\
\rho_{4} &\quad (x,y)\in(0,5)\times(-5,0).
\end{cases}
\]
and $\mu_{0}(x,y)=\rho_{0}(x,y)/2$. For simplicity of notation we take the values $\rho_{i}\in\{1,2,3,4\}$ and then normalize $\rho$ and $\mu$ dividing by $\rm=\rho^{\mathrm{max}}+\mu^{\mathrm{max}}=6$.

Since we aim to recover the results of the analysis done in Section \ref{sec:riemann}, where the plots are defined for the plane $(\xi,\eta)$ with $\xi=x/t$ and $\eta=y/t$, we note that for $t=1$ the variables $\xi$ and $\eta$ coincide with $x$ and $y$. Therefore, we plot the contours of the numerical solution $U^{n}_{ij}$ at time $t^{n}=1$ in order to identify the plane $(x,y)$ with the plane $(\xi,\eta)$ and 
for a better comparison.

\begin{figure}[h!]
\centering
\subfloat[][No shocks.]
{\label{fig:caso1NUM}\includegraphics[width=.31\columnwidth]{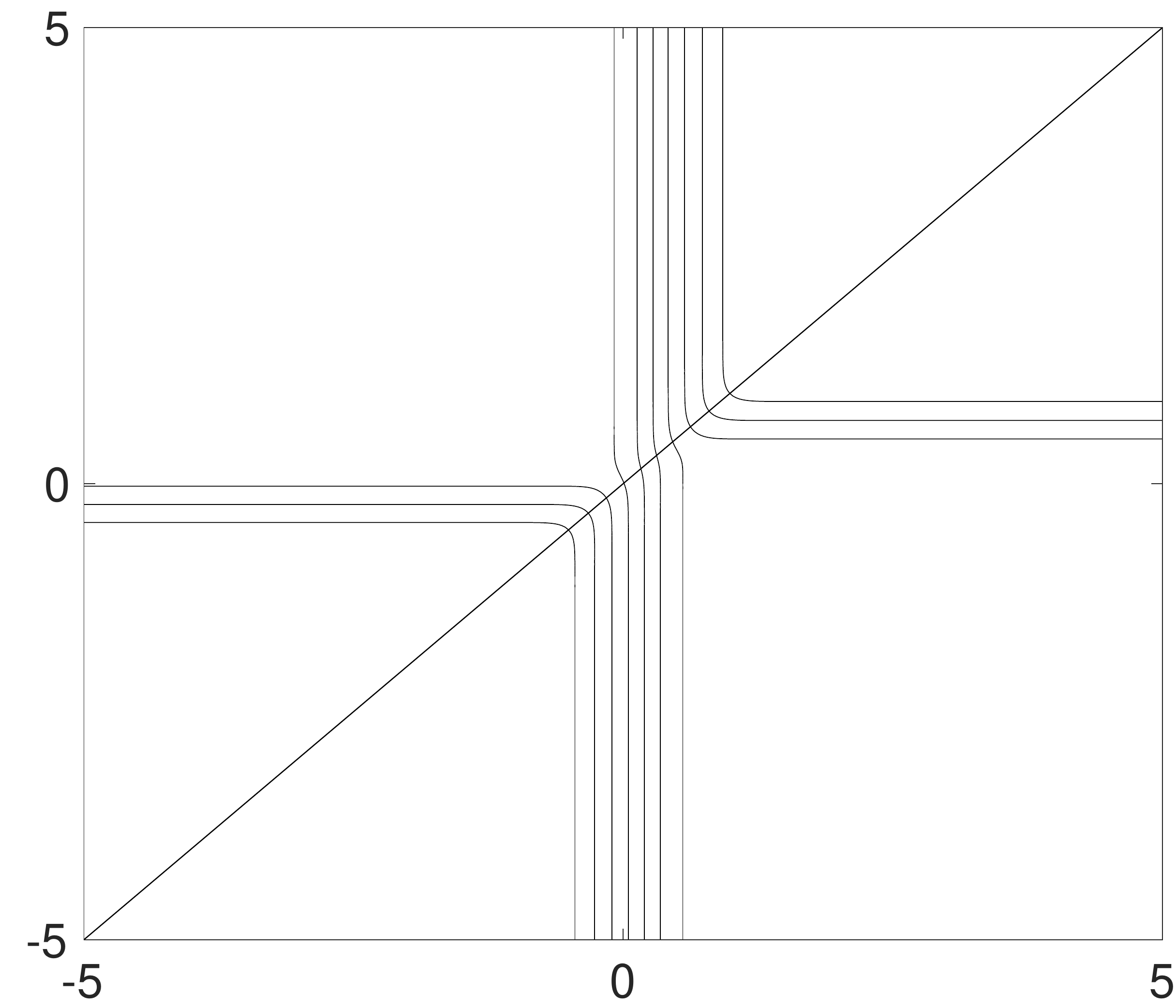}} \,
\subfloat[][No rarefaction waves.]
{\label{fig:caso2NUM}\includegraphics[width=.31\columnwidth]{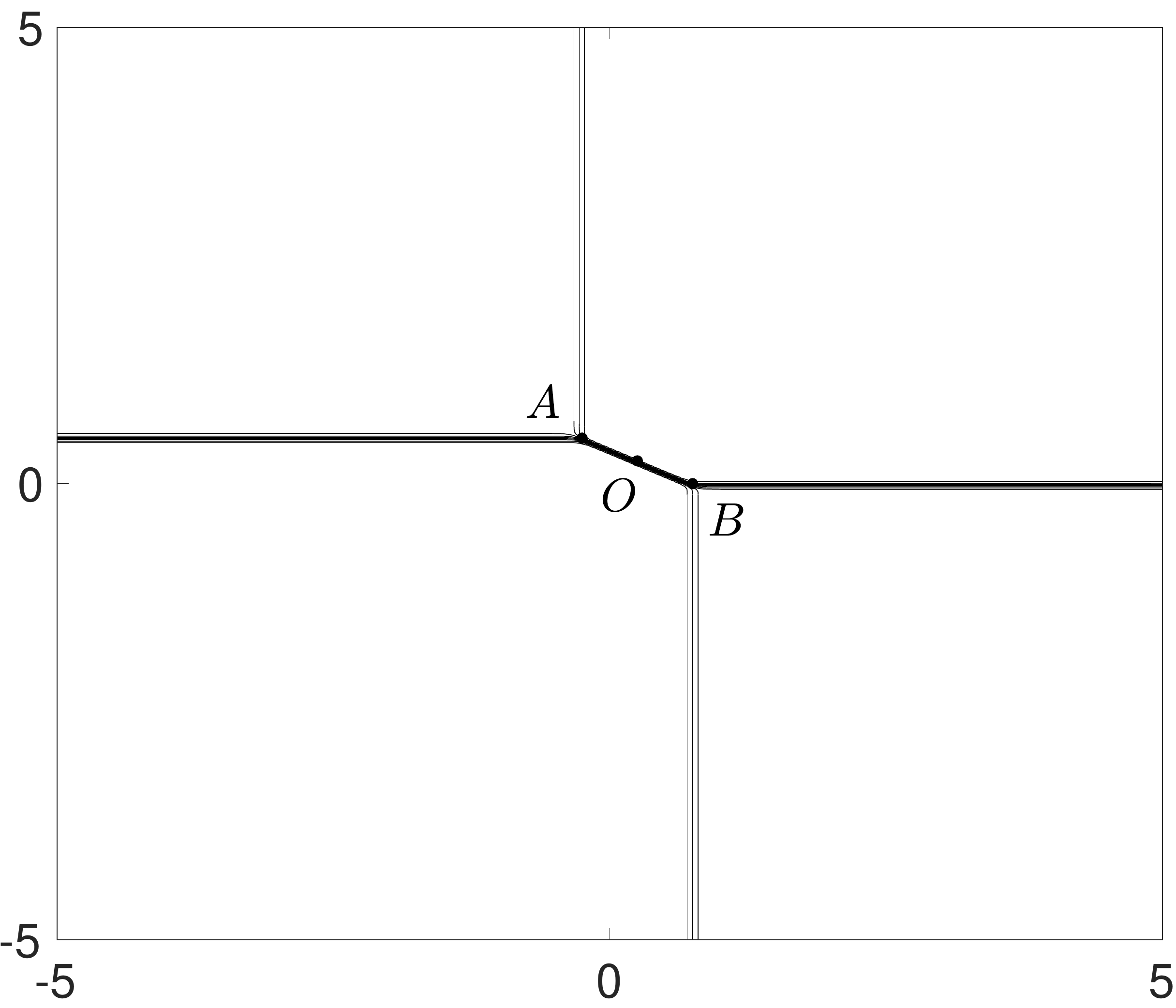}}
\,
\subfloat[][Exactly one shock.]
{\label{fig:caso3NUM}\includegraphics[width=.31\columnwidth]{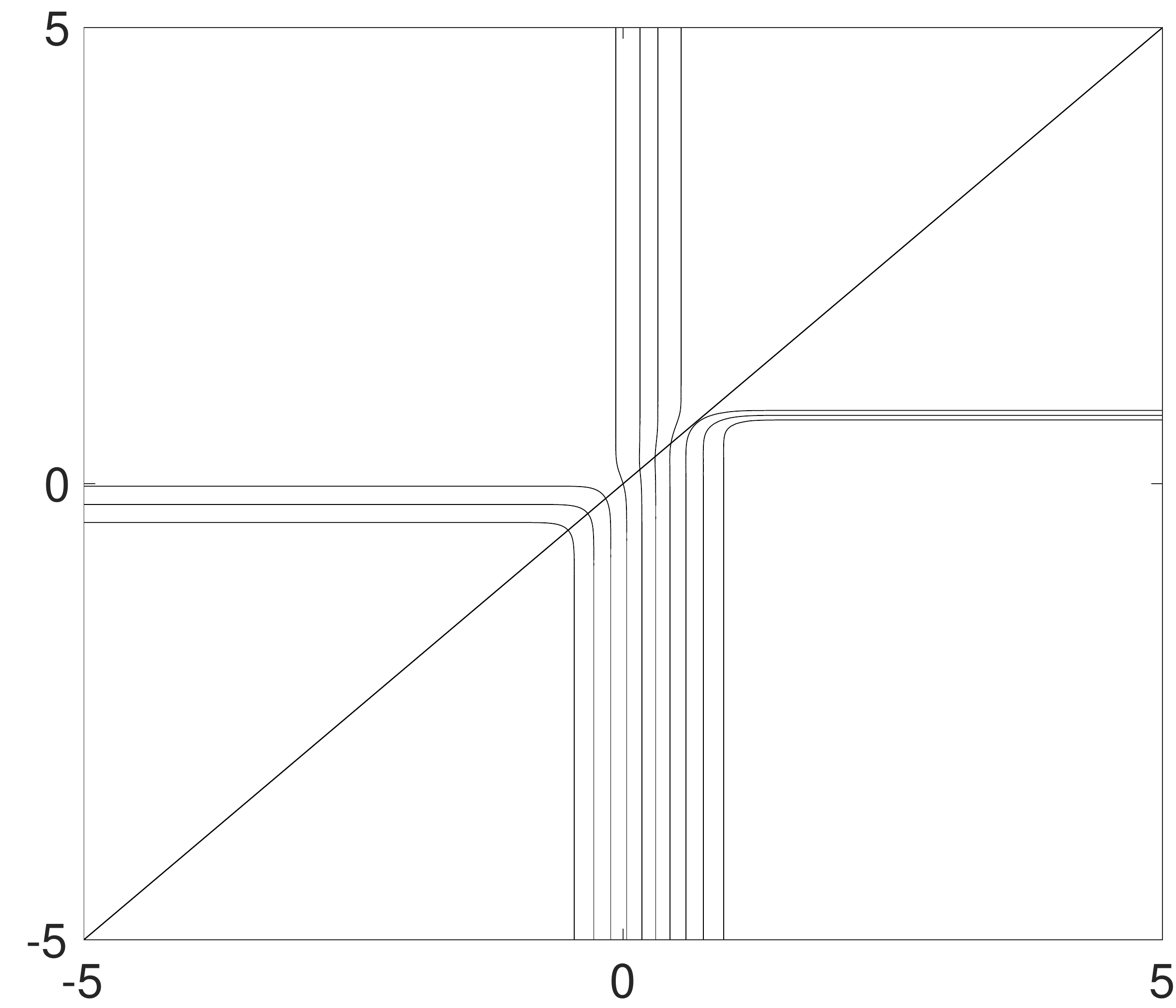}}
\\
\subfloat[][Exactly one rarefaction wave.]
{\label{fig:caso4NUM}\includegraphics[width=.31\columnwidth]{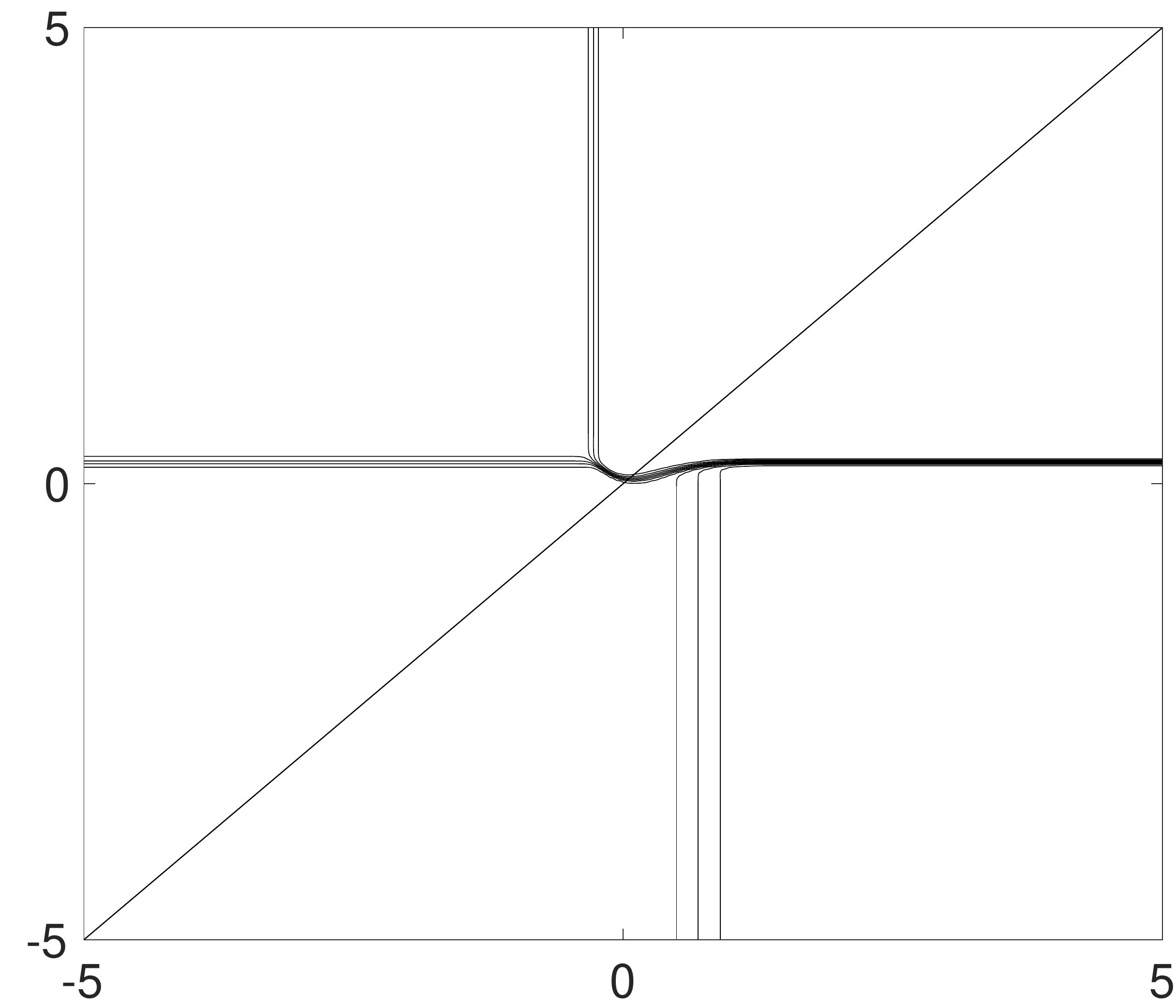}}
\,
\subfloat[][Two shocks and two rarefaction waves.]
{\label{fig:caso5NUM}\includegraphics[width=.31\columnwidth]{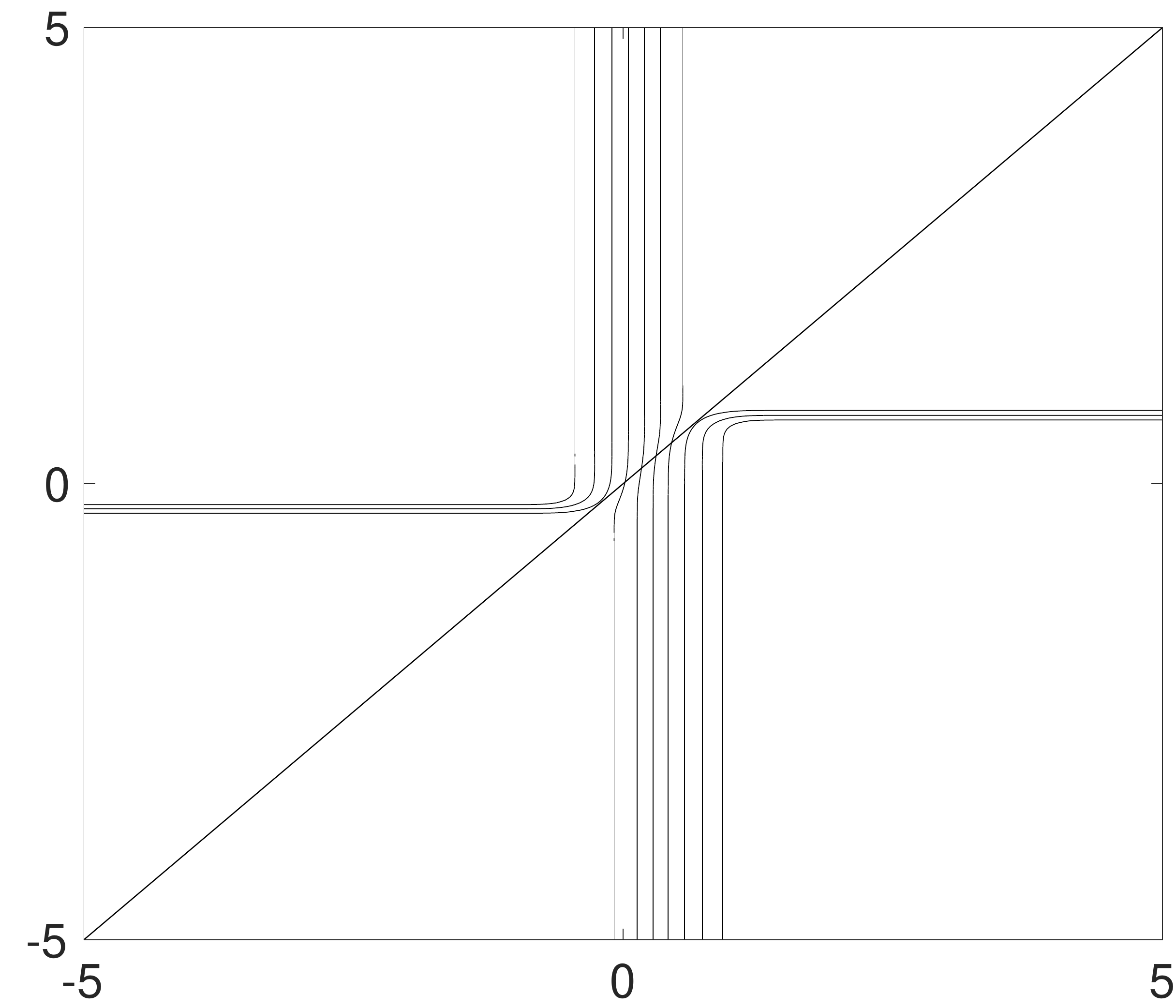}}
\caption{Numerical solutions of Riemann problems depending on the initial datum.}
\end{figure}

As we have explained in Section \ref{sec:riemann}, there are only five possible configurations of the solution, which are determined by the initial values $\rho_{i}$, $i=1,\dots,4$.\\[1ex]
\noindent
{\bf (1) No shocks:} We fix $\rho_{1}=4,\, \rho_{2}=2,\, \rho_{3}=1$ and $\rho_{4}=3$. As shown in Figure \subref*{fig:caso1NUM}, we have only rarefaction waves connected by the straight line $\esse=\{(x,y) \,|\, y=x\}$. The results in Figure \subref*{fig:caso1NUM} coincide with the theoretical solution shown in Figure \subref*{fig:caso1}. \\[1ex]
{\bf (2) No rarefaction waves:} We fix $\rho_{1}=1,\, \rho_{2}=2, \,\rho_{3}=4$ and $\rho_{4}=3$. In Figure \subref*{fig:caso2NUM}, the points of connection between the shocks are $A = (\gamma(\rho_{1}+\mu_{1},\rho_{2}+\mu_{2}),\nu(\rho_{2}+\mu_{2},\rho_{3}+\mu_{3}))=(0.25,0.5), B = (\gamma(\rho_{3}+\mu_{3},\rho_{4}+\mu_{4}),\nu(\rho_{4}+\mu_{4},\rho_{1}+\mu_{1}))=(0.75,0)$ and $O = (\gamma(\rho_{1}+\mu_{1},\rho_{3}+\mu_{3}),\nu(\rho_{1}+\mu_{1},\rho_{3}+\mu_{3}))=(0.25,0.25)$ with $\gamma$ and $\nu$ defined in \eqref{eq:gamma}. The results in Figure \subref*{fig:caso2NUM} coincide with the theoretical solution shown in Figure \subref*{fig:caso2}. \\[1ex]
{\bf (3) Exactly one shock:} We fix $\rho_{1}=3,\, \rho_{2}=2,\, \rho_{3}=1$ and $\rho_{4}=4$. As shown in Figure \subref*{fig:caso3NUM}, we consider the first sub-case described in Section \ref{sec:riemann}, and the shock wave is below the straight line $\esse=\{(x,y) \,|\, y=x\}$. The results in Figure \subref*{fig:caso3NUM} coincide with the theoretical solution shown in Figure \subref*{fig:casoc1}. \\[1ex]
{\bf (4) Exactly one rarefaction wave:} We fix $\rho_{1}=1,\, \rho_{2}=2, \,\rho_{3}=3$ and $\rho_{4}=4$. The results in Figure \subref*{fig:caso4NUM} coincide with the theoretical solution shown in Figure \subref*{fig:casoc1}, and similarly to the previous case we have that the only rarefaction wave is below the straight line $\esse=\{(x,y) \,|\, y=x\}$. \\[1ex]
{\bf (5) Two shocks and two rarefaction waves:} We fix $\rho_{1}=3,\, \rho_{2}=1, \,\rho_{3}=2$ and $\rho_{4}=4$. As shown in Figure \subref*{fig:caso5NUM}, the shock waves are not neighbors, but they are separated by the rarefaction waves and the straight line $\esse=\{(x,y) \,|\, y=x\}$. The results in Figure \subref*{fig:caso5NUM} coincide with the theoretical solution shown in Figure \subref*{fig:caso5a}. 

\section{Data-driven multi-class model in 2D}\label{sec:dati}

\begin{figure}[h!]
\centering
\includegraphics[scale=0.9]{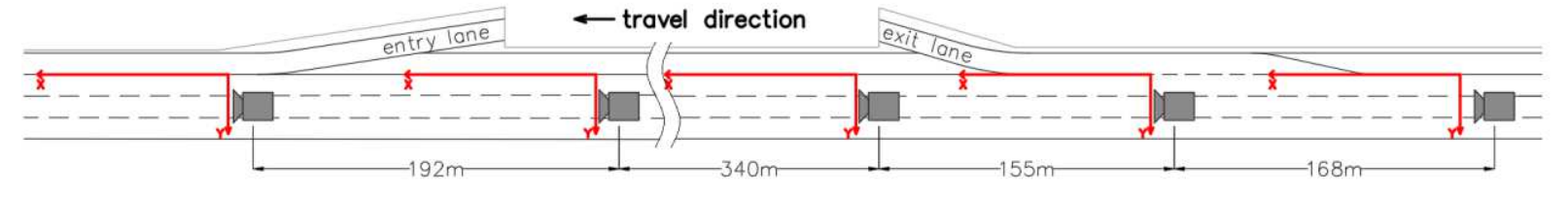}
\caption{German motorway A3 structure, cf.~\cite{germanDataset}.}
\label{fig:strada}
\end{figure}

In this section, we calibrate the two-dimensional multi-class model with a dataset of real trajectories data.
We employ the public German dataset \cite{germanDataset} which contains vehicle trajectories data recorded on the German motorway A3, nearby Frankfurt am Main. The analyzed area is a three lanes highway of about 900 meters in length and 12 meters in width, depicted in Figure \ref{fig:strada}. A system of five video cameras recorded the vehicles passing through the study area, collecting trajectory data for 20 minutes with a sampling period of about 0.2 seconds. We refer to \cite{germanDataset} for a detailed description of the dataset and of the data collection method. 
We observe that the dataset distinguishes several types of vehicles, and particularly in this work we focus on the dynamics of cars and trucks. 

\subsection{Fundamental diagrams}\label{sec:FD}
In order to calibrate the two-dimensional model with the German dataset, we need first to derive macroscopic quantities from the microscopic information provided by the dataset. Following \cite{herty2018NHM}, we describe how we derive the density of cars and trucks, $\rho$ and $\mu$, the speed in the two directions, $\ux$ and $\uy$, and the flux in the two directions, $q^{x}_{\rho,\mu}$ and $q^{y}_{\rho,\mu}$. Note that, as we have already observed in Section \ref{sec:modello}, the velocity functions coincide for $\rho$ and $\mu$ since $\cx$ and $\cy$ do not distinguish the class of vehicles. 

We consider the data from the second camera from the right of Figure \ref{fig:strada}, thus we work with 20 minutes of real data. We introduce the time interval ($t_{0},t_{M})$, with $t_{0}=0$ and $t_{M}=20\,\minute$, and uniformly divide it with a time step $dt$. Note that $dt$ is used to derive the macroscopic quantities from the microscopic ones, and is independent of the time step $\dt$ of the numerical scheme. 
We call $N_{\rho,\mu}(t_{k})$ the total number of cars and trucks at time $t_{k}$ and $\lx$ the length of the road along the main direction of travel. Then, we define
\begin{equation}\label{eq:density}
\rhot(t_{k}) = \frac{N_{\rho}(t_{k})}{\lx}, \qquad \mut(t_{k}) = \frac{N_{\mu}(t_{k})}{\lx}.
\end{equation} 

The German dataset only provides the position of vehicles with respect to the two directions, thus we need to derive the speed of vehicles from their positions. We assume that each vehicle travels at constant speed which corresponds to the slope of a linear approximation in the least square sense of the vehicle positions. We denote by $v^{x,y}_{i}$ the resulting microscopic speed of car $i$ and by $w^{x,y}_{i}$ the analogous speed of truck $i$. Since we assume that the two classes have the same speed function, we define the average speed as a function of the two classes
\begin{equation}\label{eq:u1}
	\begin{split}
		\uxt(t_{k})&=\frac{1}{N_{\rho}(t_{k})}\sum_{i=1}^{N_{\rho}(t_{k})}v^{x}_{i} + 	\frac{1}{N_{\mu}(t_{k})}\sum_{i=1}^{N_{\mu}(t_{k})}w^{x}_{i}\\
		\uyt(t_{k})&=\frac{1}{N_{\rho}(t_{k})}\sum_{i=1}^{N_{\rho}(t_{k})}v^{y}_{i}+\frac{1}{N_{\mu}(t_{k})}\sum_{i=1}^{N_{\mu}(t_{k})}w^{y}_{i}.
	\end{split}
\end{equation}
Finally, we combine equations \eqref{eq:density} and \eqref{eq:u1} to define the flux functions as
\begin{equation}\label{eq:q1}
	\begin{split}
		\qxrt(t_{k}) &= \rhot(t_{k})\uxt(t_{k}), \qquad\, 
		\qyrt(t_{k}) = \rhot(t_{k})\uyt(t_{k})\\  
		\qxmt(t_{k}) &= \mut(t_{k})\uxt(t_{k}), \qquad
		\qymt(t_{k}) = \mut(t_{k})\uyt(t_{k}).
	\end{split}
\end{equation}

\begin{figure}[H]
\centering
\subfloat[][Ground-truth speed $\uxt$.]{
\includegraphics[width=0.31\columnwidth]{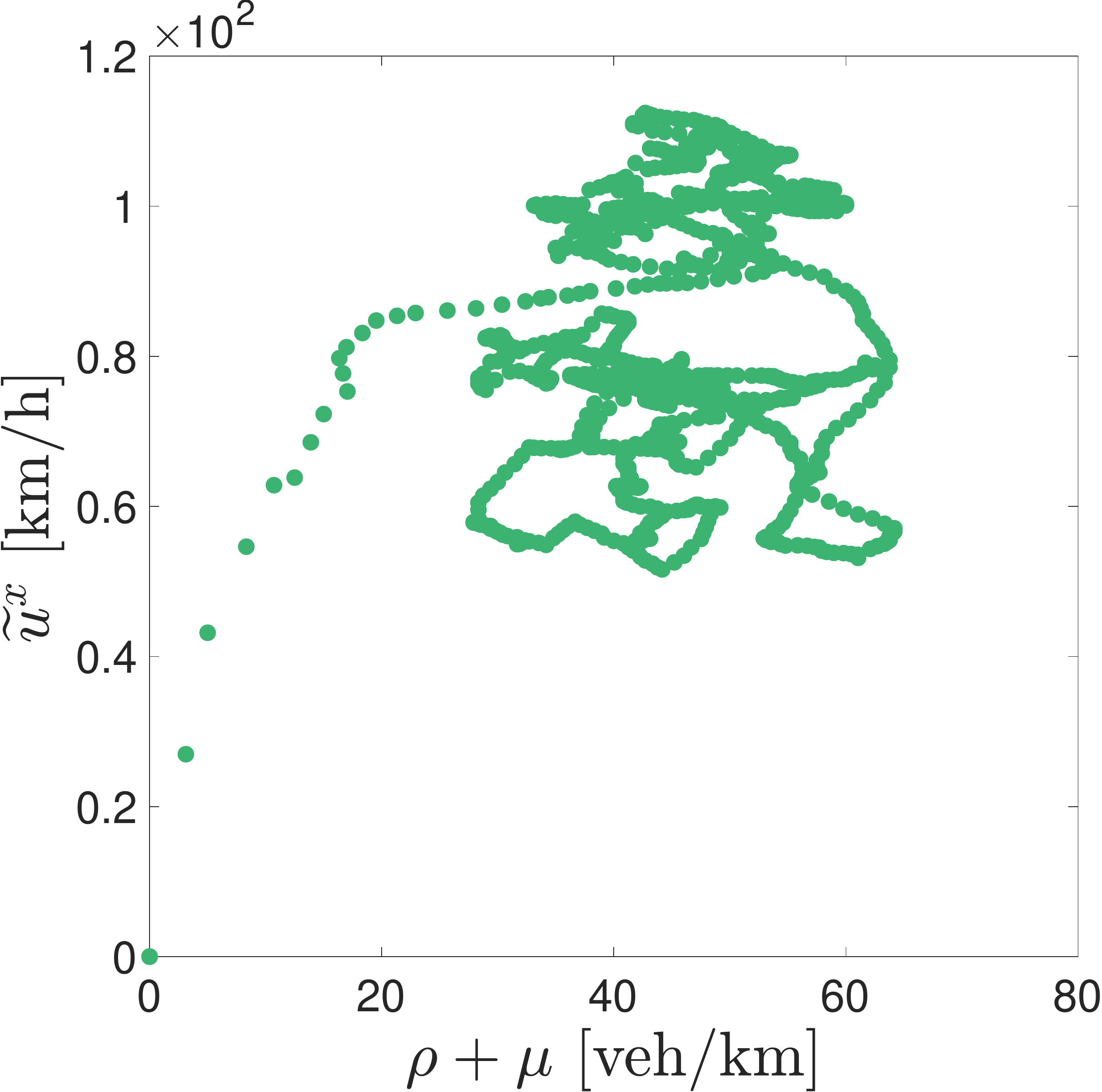}
}\,
\subfloat[][Ground-truth flux $\qxrt$.]{
\includegraphics[width=0.31\columnwidth]{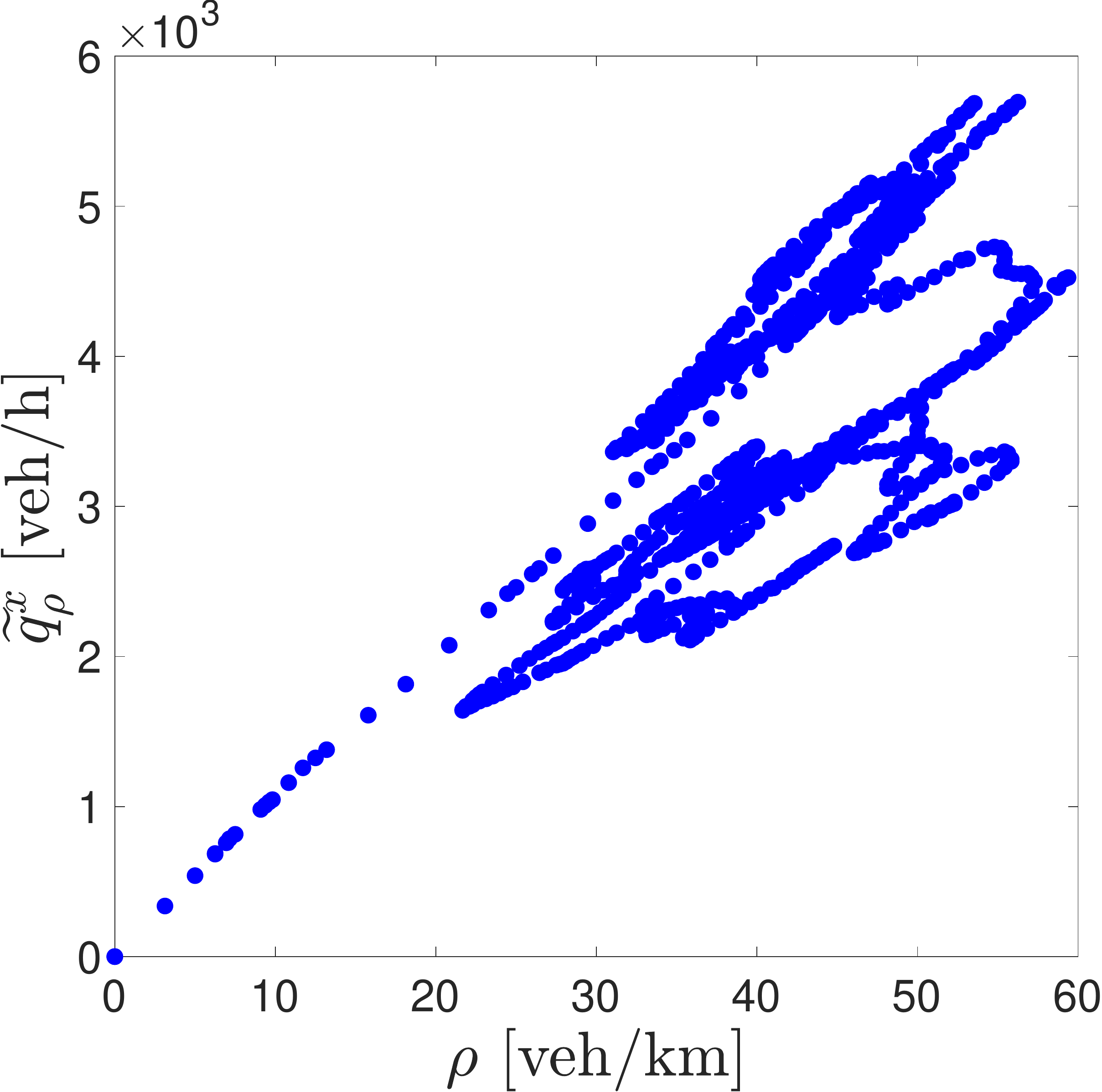}
}\,
\subfloat[][Ground-truth flux $\qxmt$.]{
\includegraphics[width=0.31\columnwidth]{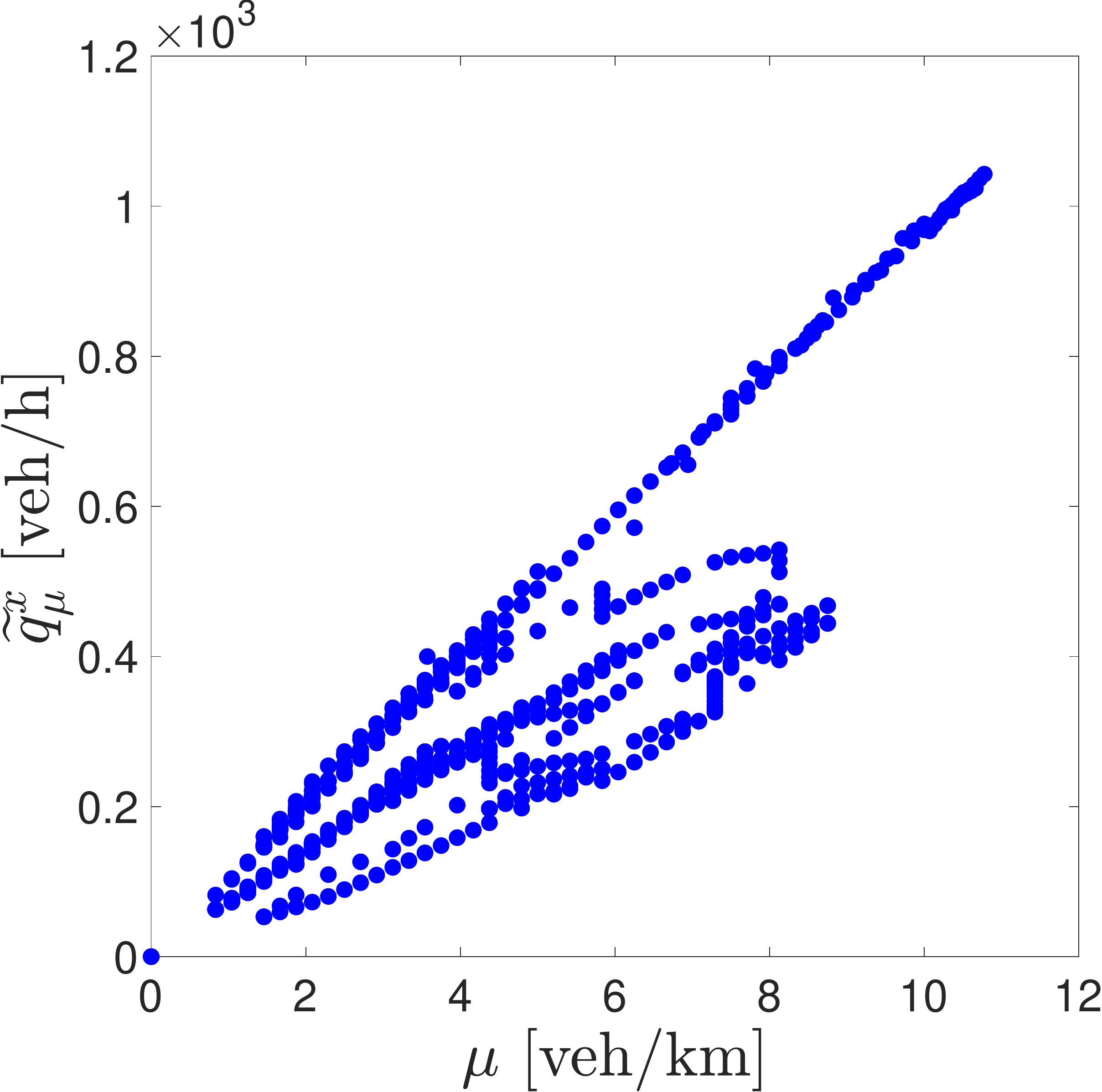}
}\\
\subfloat[][Ground-truth speed $\uyt$.]{
\includegraphics[width=0.31\columnwidth]{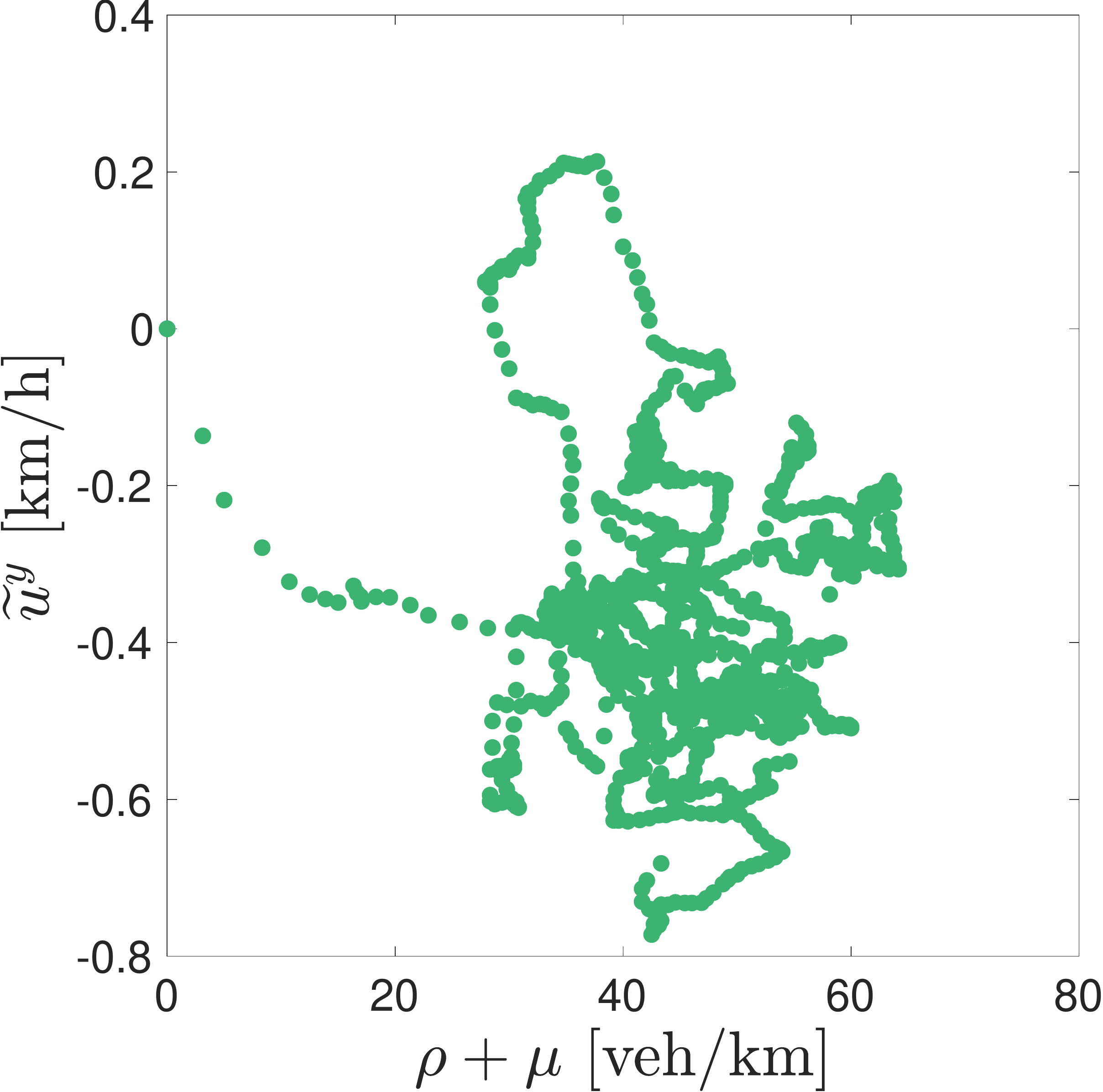}
}\,
\subfloat[][Ground-truth flux $\qyrt$.]{
\includegraphics[width=0.31\columnwidth]{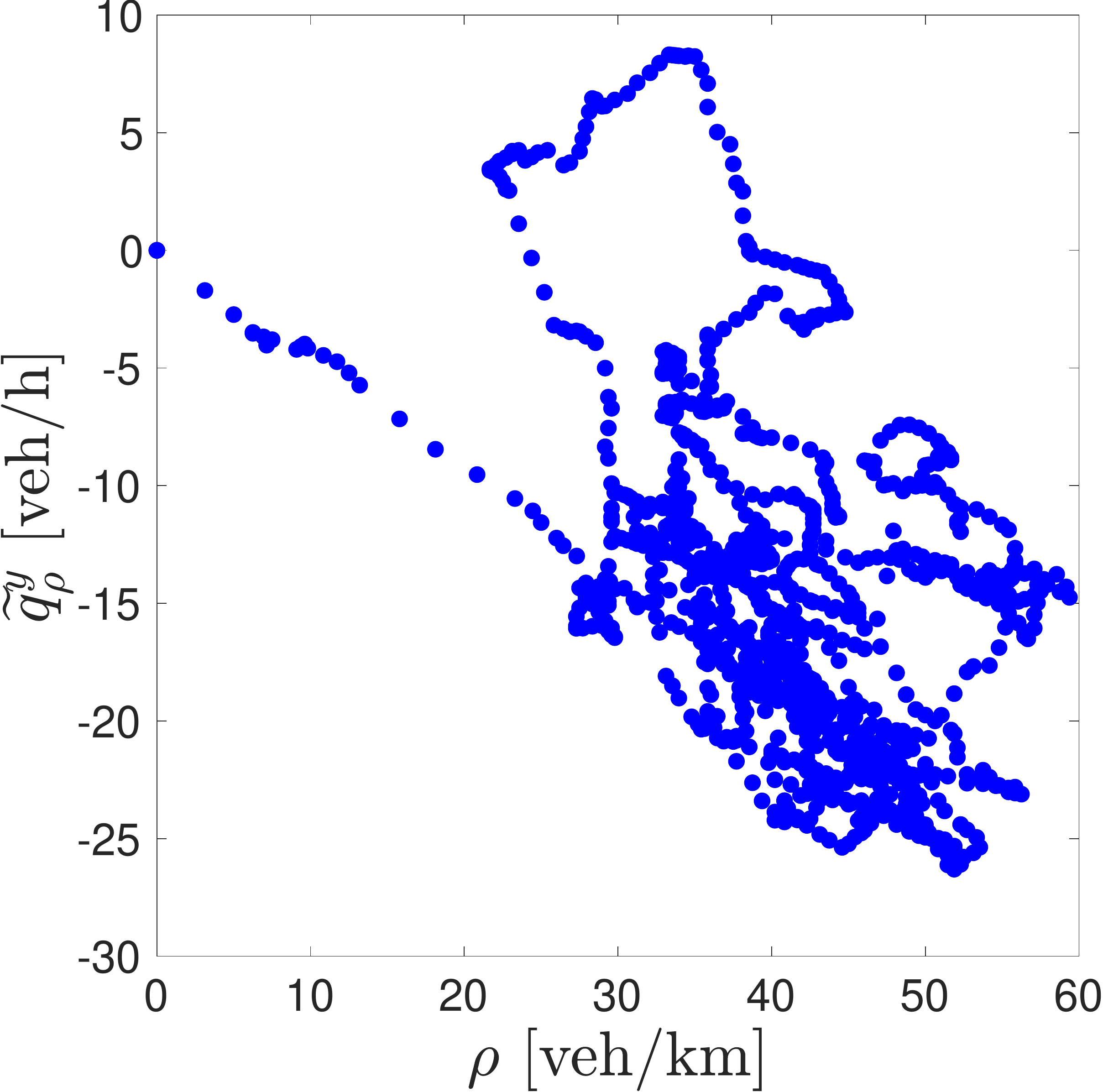}
}\,
\subfloat[][Ground-truth flux $\qymt$.]{
\includegraphics[width=0.31\columnwidth]{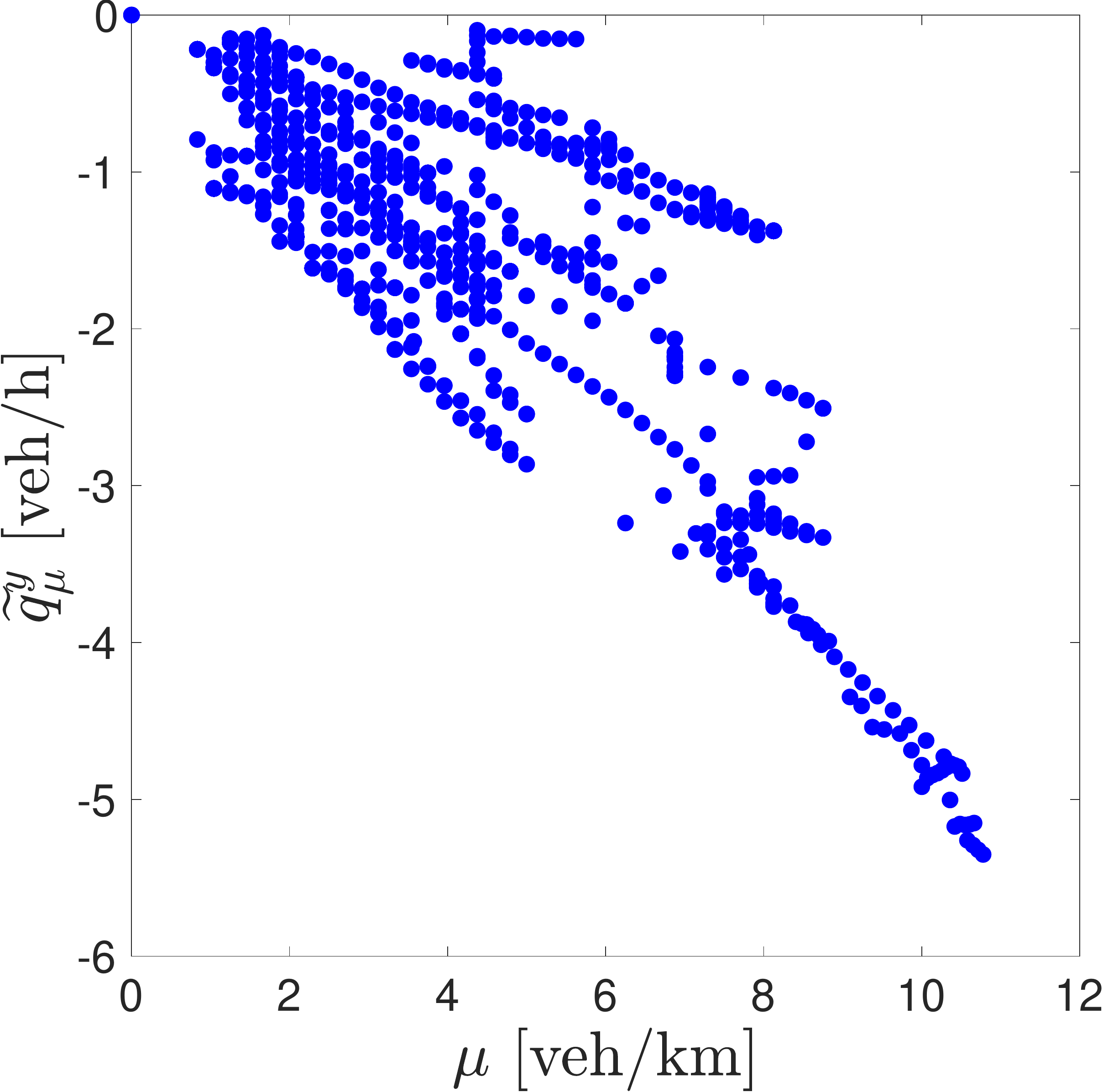}}
\caption{Speed-density and flux-density diagrams for the two classes related to the $x$-direction in the first row, and to the $y$-direction on the second row.}
\label{fig:fluxspeedTrue}
\end{figure}

Once we have the density, speed and flux data as functions of time, we aggregate them  with respect to a certain time period $\widetilde T=\kappa dt$. In particular, we fix $dt=1\,\second$ and $\kappa=60$. 
In Figure \ref{fig:fluxspeedTrue}, we show the speed-density and flux-density diagrams for the two classes of vehicles in the $x$ and $y$ directions. 

The graphs show that the main direction of the flow is along the $x$-axis, according to the structure of the analyzed road, while the movements along the $y$-axis represent the lane changes. Note that both the flux and the velocity along the $y$-direction show negative values, due to the lane change that can occur in both directions.
The maximum density value reached in Figure \ref{fig:fluxspeedTrue} is $60\,\vehkm$ for cars and $12\,\vehkm$ for trucks, and both values are much smaller than the maximum density of the road given by $\rm=400\,\vehkm$. More specifically, 
as we observe from Figure \ref{fig:fluxspeedTrue}, there are more cars than trucks along the road, thus we calibrate the maximum density fixing the length of vehicles as if there are only cars on the road. Hence, we assume that the length of vehicles plus the safety distance is $7.5\,\meter$, thus $\rm$ is defined as
\begin{equation}\label{eq:rmax}
\rm = \frac{\text{\# lanes}}{\text{length of vehicles + safety distance}} = \frac{3}{7.5 \text{ m}} = 400\, \frac{\text{veh}}{\text{km}}.
\end{equation}
However, it should be noted that the dataset only contains data in free-flow regimes without capturing congested traffic phase.

\begin{figure}[h!]
\centering
\subfloat[][Ground-truth speed $\uxt$ and family of speed functions $\ux$ as $\mu$ changes.]{\label{fig:speedx}
\includegraphics[width=0.31\columnwidth]{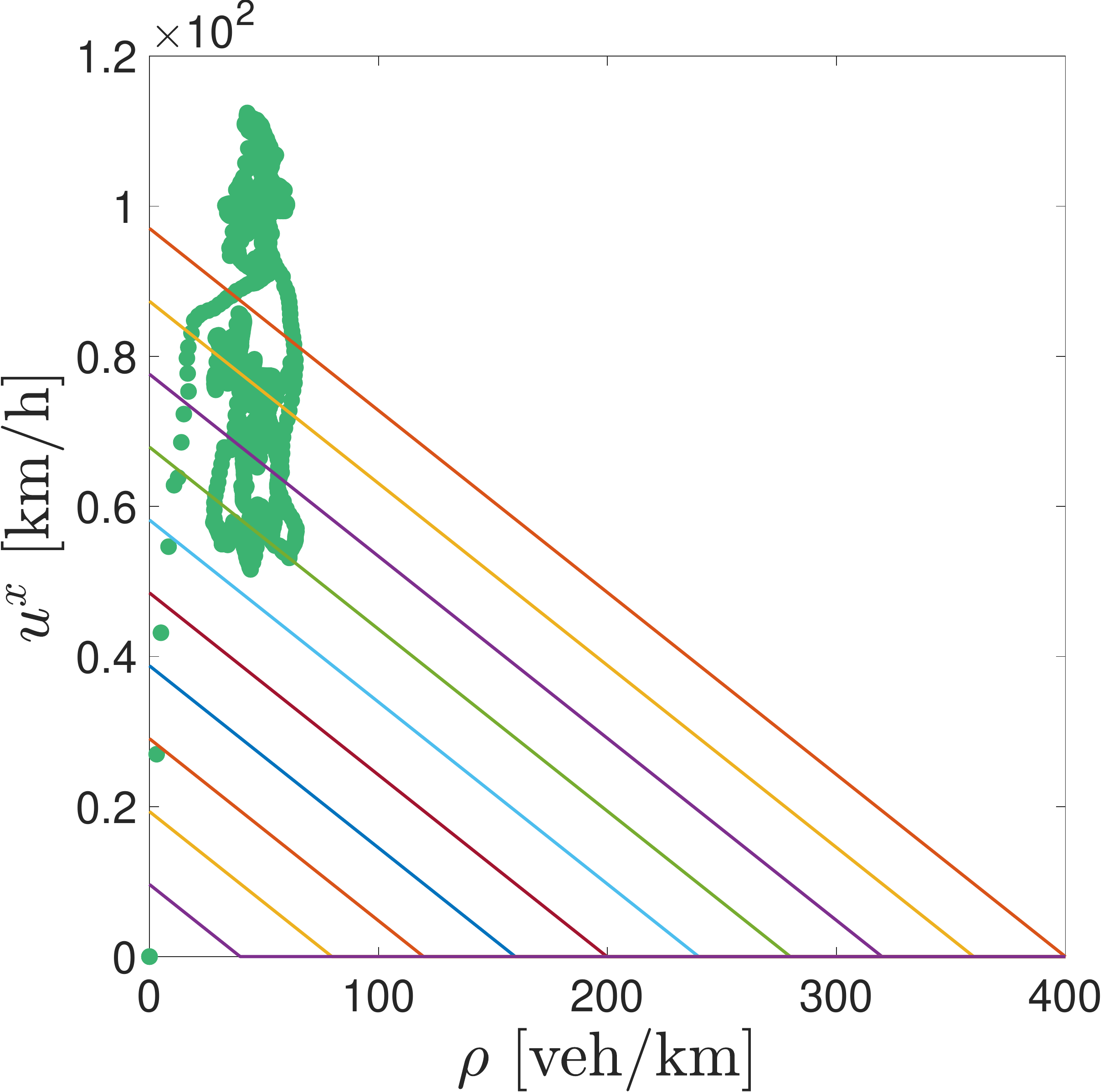}
}\,
\subfloat[][Ground-truth flux $\qxrt$ and family of flux functions $\qxr$ as $\mu$ changes.]{\label{fig:flussoxr}
\includegraphics[width=0.31\columnwidth]{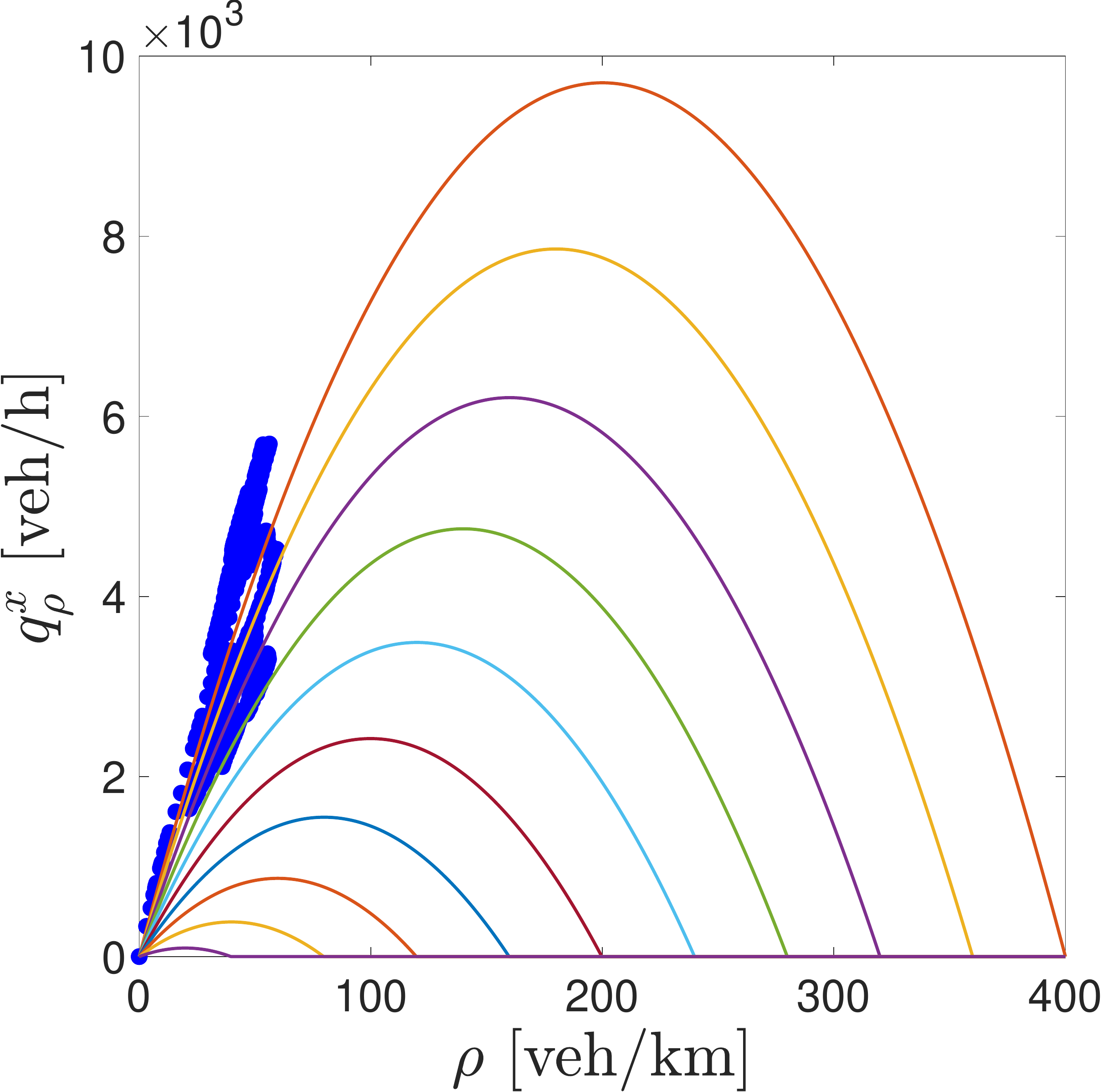}
}\,
\subfloat[][Ground-truth flux $\qxmt$ and family of flux functions $\qxm$ as $\rho$ changes.]{\label{fig:flussoxm}
\includegraphics[width=0.31\columnwidth]{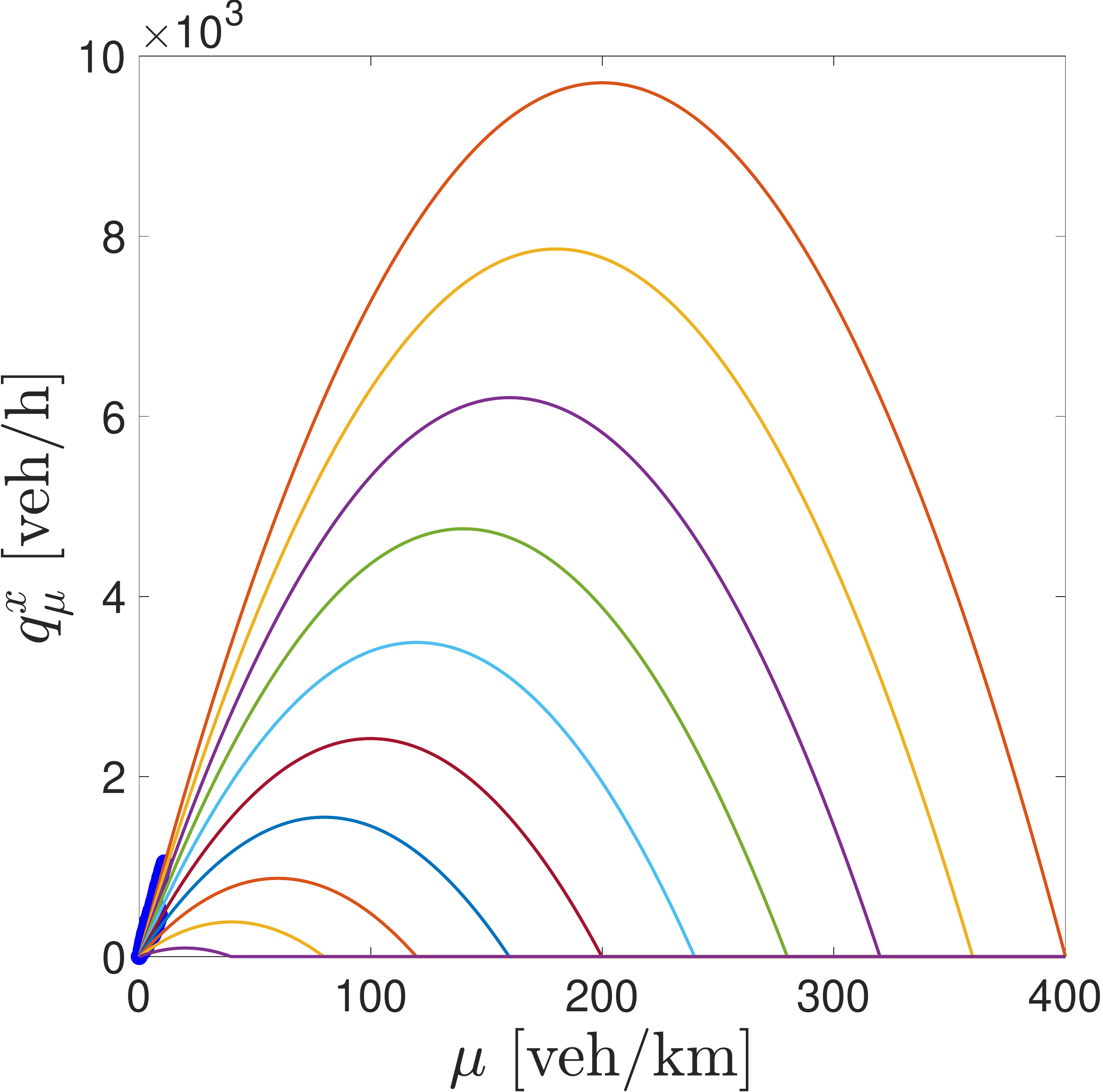}
}\\
\subfloat[][Ground-truth speed $\uyt$ and family of speed functions $\uy$ as $\mu$ changes.]{\label{fig:speedy}
\includegraphics[width=0.31\columnwidth]{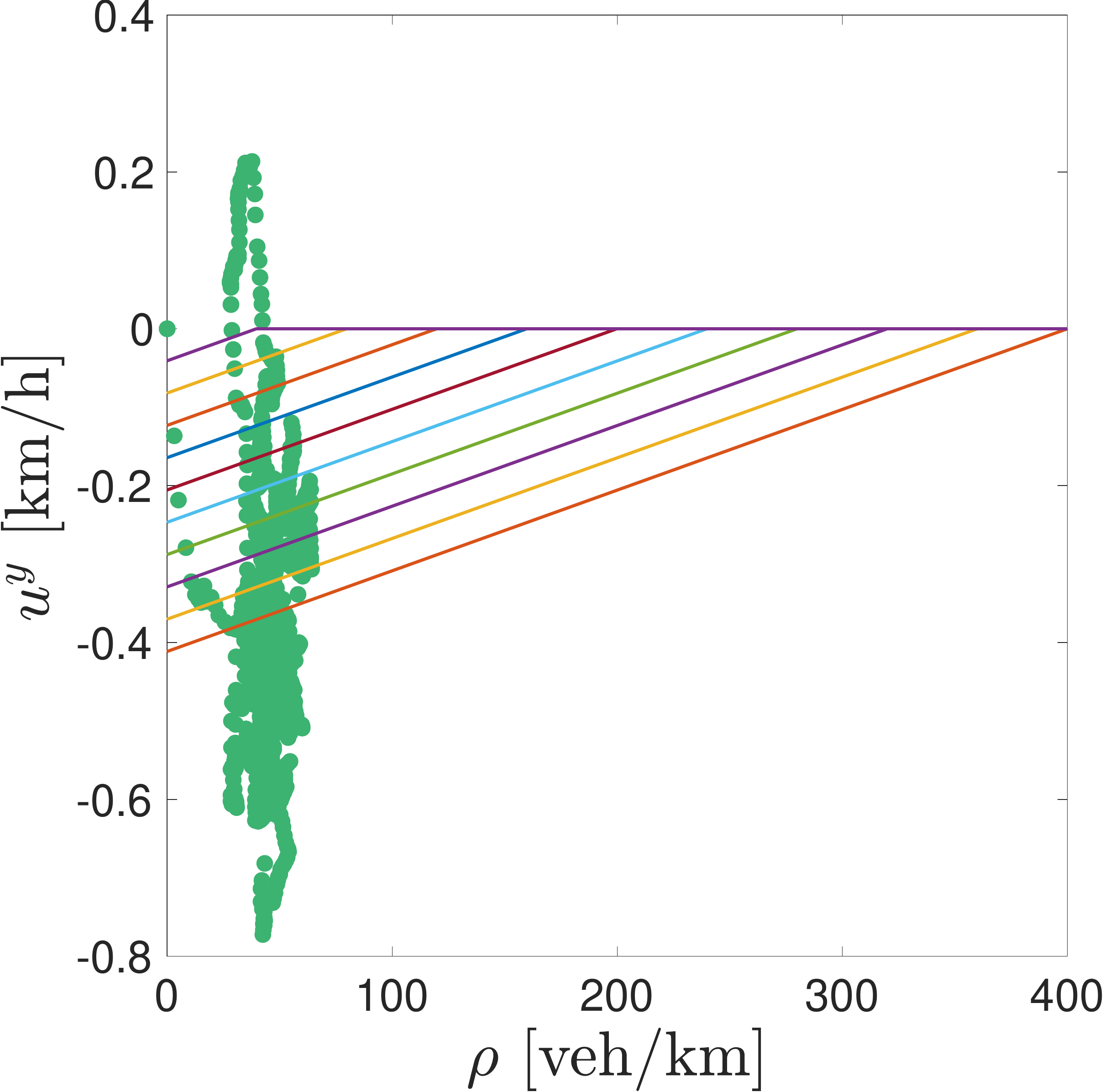}
}\,
\subfloat[][Ground-truth flux $\qyrt$ and family of flux functions $\qyr$ as $\mu$ changes.]{\label{fig:flussoyr}
\includegraphics[width=0.31\columnwidth]{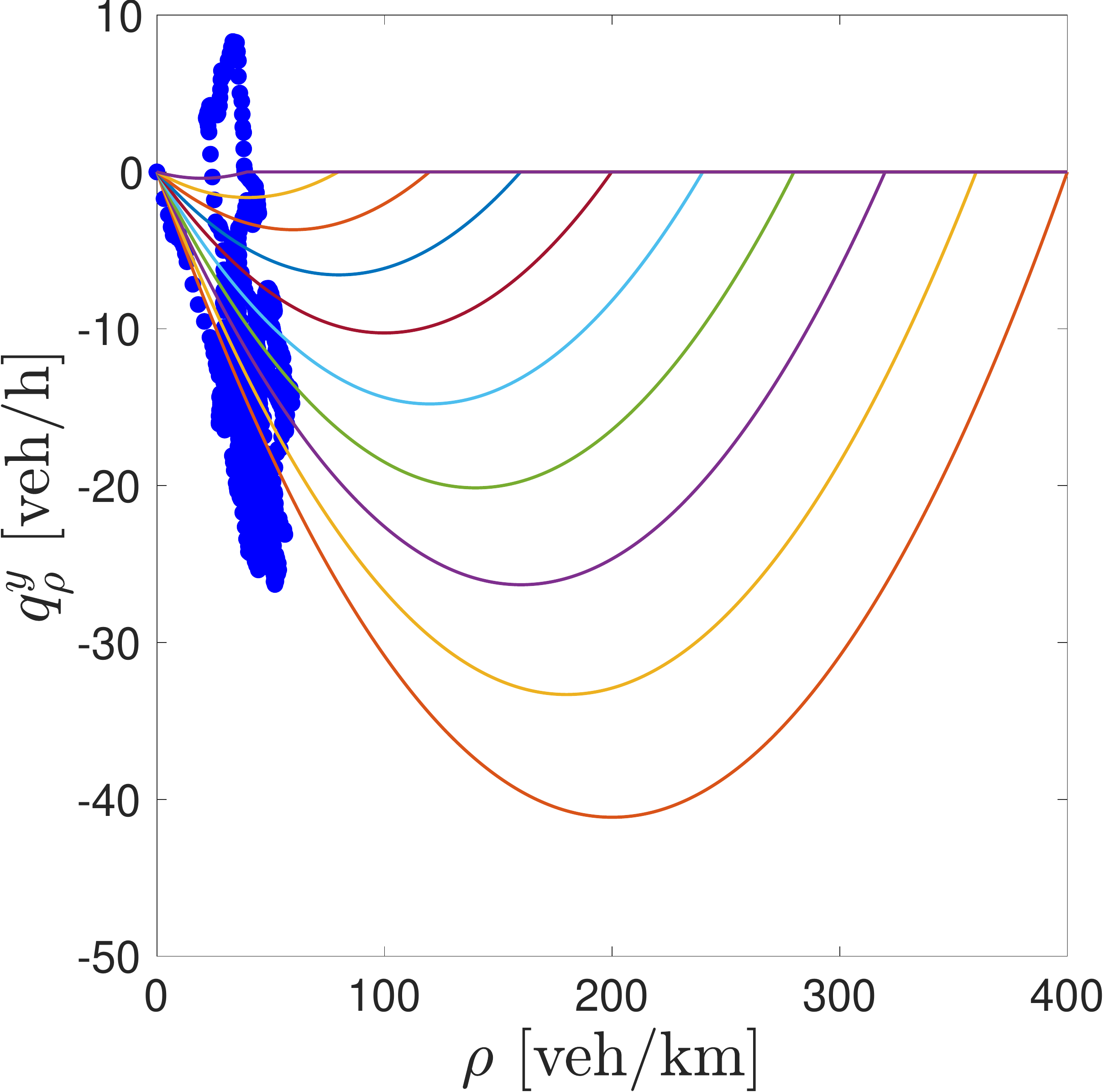}
}\,
\subfloat[][Ground-truth flux $\qymt$ and family of flux functions $\qym$ as $\rho$ changes.]{\label{fig:flussoym}
\includegraphics[width=0.31\columnwidth]{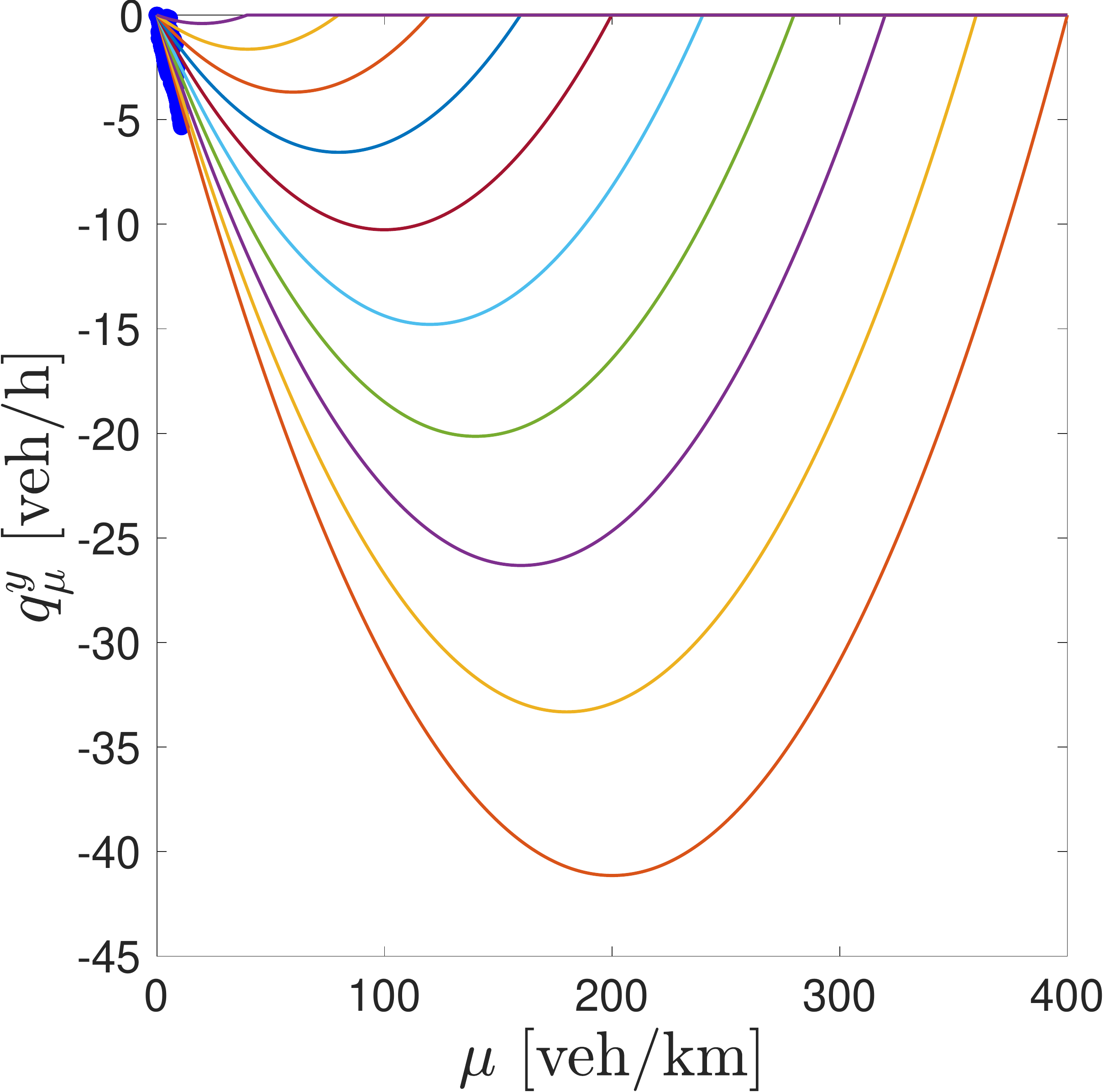}}
\caption{Speed-density and flux-density diagrams for the two classes defined from real data (green and blue circles) and family of speed and velocity functions related to the $x$-direction in the first row, and to the $y$-direction in the second row.}
\label{fig:fluxspeed2C}
\end{figure}

Now we need to compute the parameters $\cx$ and $\cy$. 
The parameters $\cx$ and $\cy$ are chosen in order to minimize the $L^{2}$-norm between the flux functions defined in \eqref{eq:flussi} and the fluxes derived from data in \eqref{eq:q1}, i.e., we consider
\[\min_{\cx} \left(\norm{\qxrt-\qxr(\rho,\mu)}^{2}_{2}+\norm{\qxmt-\qxm(\rho,\mu)}^{2}_{2}\right),\quad
\min_{\cy}\left(\norm{\qyrt-\qyr(\rho,\mu)}^{2}_{2}+\norm{\qymt-\qym(\rho,\mu)}^{2}_{2}\right).
\]
The computation is performed using the MATLAB $\mathtt{fminbnd}$ tool, which is a specific solver for minimization problems. We obtain $\cx=97.04$ and $\cy=-0.41$. Note that $\cy$ is negative, since the lane changes occur mainly towards the rightmost lane. The resulting speed and flux functions are shown in Figure \ref{fig:fluxspeed2C}. 
Since the flux functions depend on both $\rho$ and $\mu$, we have a family of flux and velocity functions. 
In particular, in Figure \subref*{fig:flussoxr} we show the family of flux functions $\qxr$ as $\mu$ changes. 
This means that at fixed value of $\mu$ we can move only along one of the flux curves.
For instance, if $\mu=0$, i.e. there are no trucks, then the fundamental diagram $\qxr$ corresponds to the maximum flux curve in Figure \subref*{fig:flussoxr}, if $\mu=\rm$ then no car can enter into the road, and thus $\qxr\equiv0$. A similar discussion holds for the other plots of Figure \ref{fig:fluxspeed2C}.

We observe that the advantage of the multi-class model is that we can cover quite well the clouds of real data by means of the family of flux and velocity functions. However, we note that, since the German dataset contains data which refer only to the not congested phase of traffic, we do not have enough data to better calibrate congested traffic situations. In particular, the choice of $\cx$ and $\cy$ equal for both of the classes seems to overestimate the flux for the class $\mu$, in both the directions.

\subsection{Reconstruction of density from data}\label{sec:KDE}
In this section we describe how to treat the microscopic data to define the initial density for the numerical scheme and the reference solution
for the comparison of the results. The German dataset gives information about the position of vehicles every 0.2 seconds, thus we work with pointwise data. In order to define a density function $\rho(t,x,y)$ on a domain $\mathcal{D},$ we use a kernel density estimation, the Parzan-Rosenblatt window method \cite{parzen1962AMS,rosenblatt1956AMS}. The idea of this method is to consider the data points as a density distribution and then recover the global density by summing these distributions. 

Let $N(t)$ be the number of cars at time $t$ and $(x_{i}(t),y_{i}(t))$ their positions, we define
\[\rhot(x,y) = \sum_{i=1}^{N(t)}\delta(x-x_{i}(t))\delta(y-y_{i}(t)).
\]
In order to recover the smooth function $\rho$, we introduce a two-dimensional Gaussian kernel
\[K(x,y) = \frac{1}{2\pi h_{x}h_{y}}\exp{\left(-\frac{x^{2}}{2h_{x}^{2}}-\frac{y^{2}}{2h_{y}^{2}}\right)},
\]
and then define
\begin{equation}
\rho(t,x,y) = \int_{\mathcal{D}}K(x-\xi,y-\eta)\tilde\rho(\xi,\eta)d\xi d\eta = \sum_{i=1}^{N(t)}K(x-x_{i}(t),y-y_{i}(t)).
\label{eq:KDE}
\end{equation}
We follow a similar procedure to estimate the density of trucks.
The parameters $h_{x}$ and $h_{y}$ are bandwidths chosen in order to obtain an almost constant density profile for equidistant vehicles \cite{fan2013AX}. These parameters depend on the dimensions of the road, i.e., on a road of dimensions $\lx\times\ly$ we fix $h_{x} = \lx/20$ and $h_{y} = \ly/20$, with $\lx$ being the length of the road along the $x$-axis and $\ly$ the length along the $y$-axis.

For each video camera of the German dataset we work with records data for about $80\,\meter$ in length and $12\,\meter$ in width. The average speed of vehicles is such that they exit from the recording area after a few seconds. In order to test longer simulations and compare them with real data, we assume that the trajectory of each vehicle can be approximated by a linear movement. Indeed, let us consider a vehicle $i$ which crosses the road between a time interval $[t_{0},t_{1}]$. We compute the coefficients $a_{i}^{x,y}$ and $b_{i}^{x,y}$ such that we can approximate the $x$ and $y$ position as $x(t) = a_{i}^{x}+b_{i}^{x}t$ and $y(t) = a_{i}^{y}+b_{i}^{y}t$ minimizing the $L^{2}$-norm of the difference with the real positions. In this way, we are able to compute the ``real'' position of vehicles even when they exit the supervised area. The computed positions also allow for a comparison to the numerical results.

\subsection{Numerical test}\label{sec:test1}
Now, we compare the numerical simulations of model \eqref{eq:lwr2D} with the real data computed from equation \eqref{eq:KDE}. The simulation refers to the data recorded by the second video camera of the German dataset. 

\begin{figure}[H]
\centering
\subfloat[][Cars density at time $t=0$.]
{\includegraphics[width=.3\columnwidth]{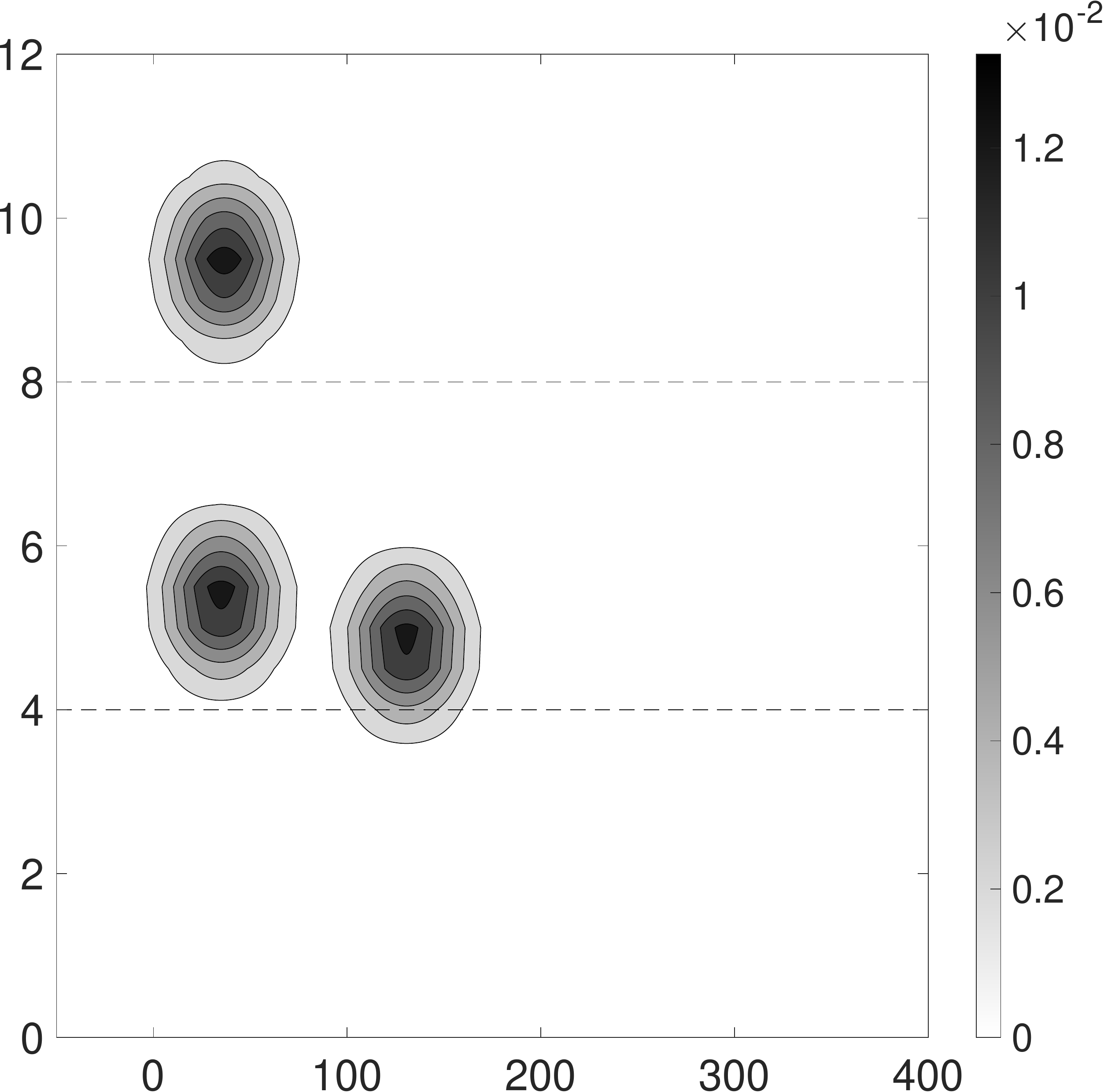}} \qquad
\subfloat[][Simulated cars density at time $t=T$.]
{\includegraphics[width=.3\columnwidth]{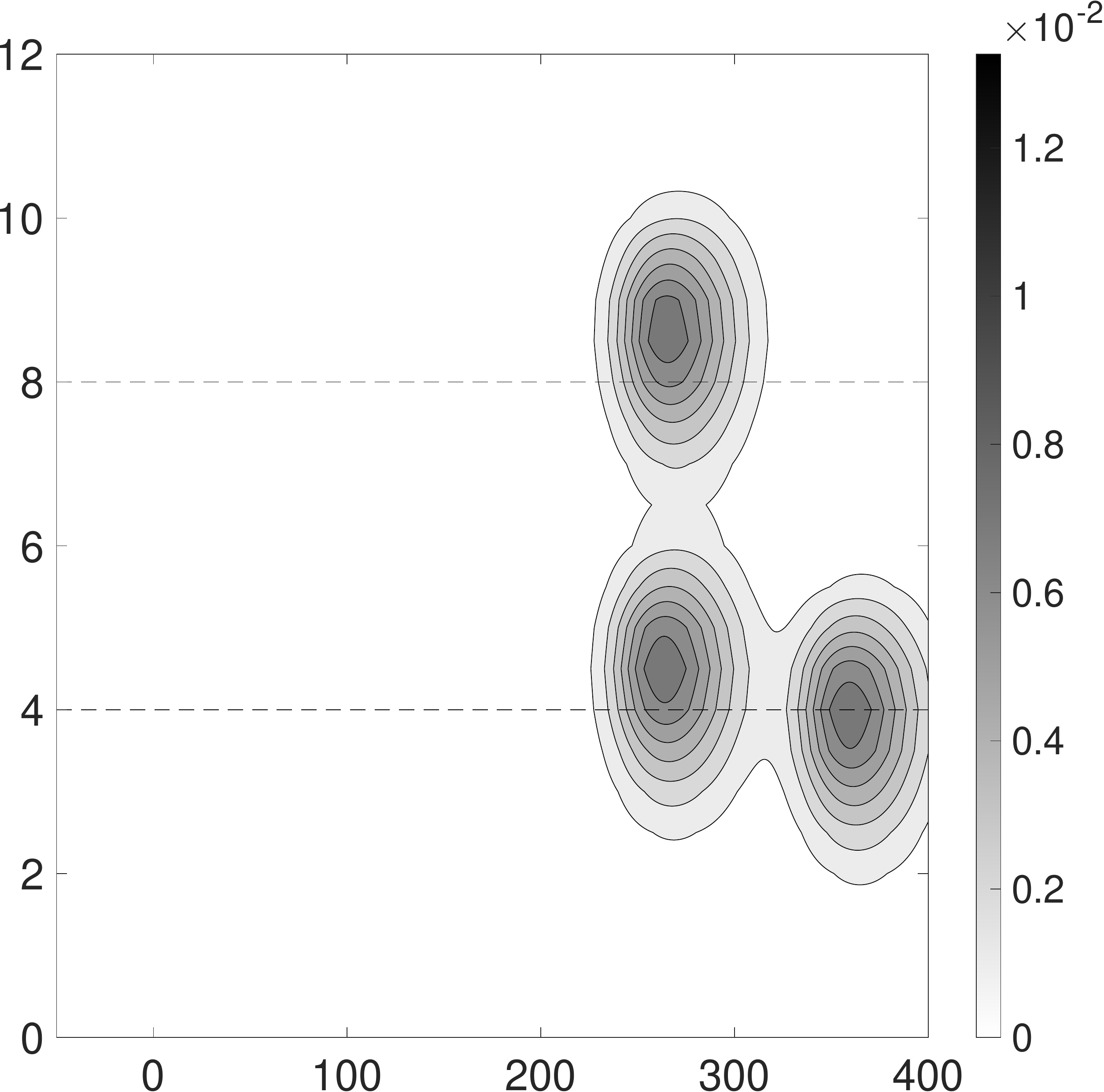}} \qquad
\subfloat[][Real cars density at time $t=T$.]
{\includegraphics[width=.3\columnwidth]{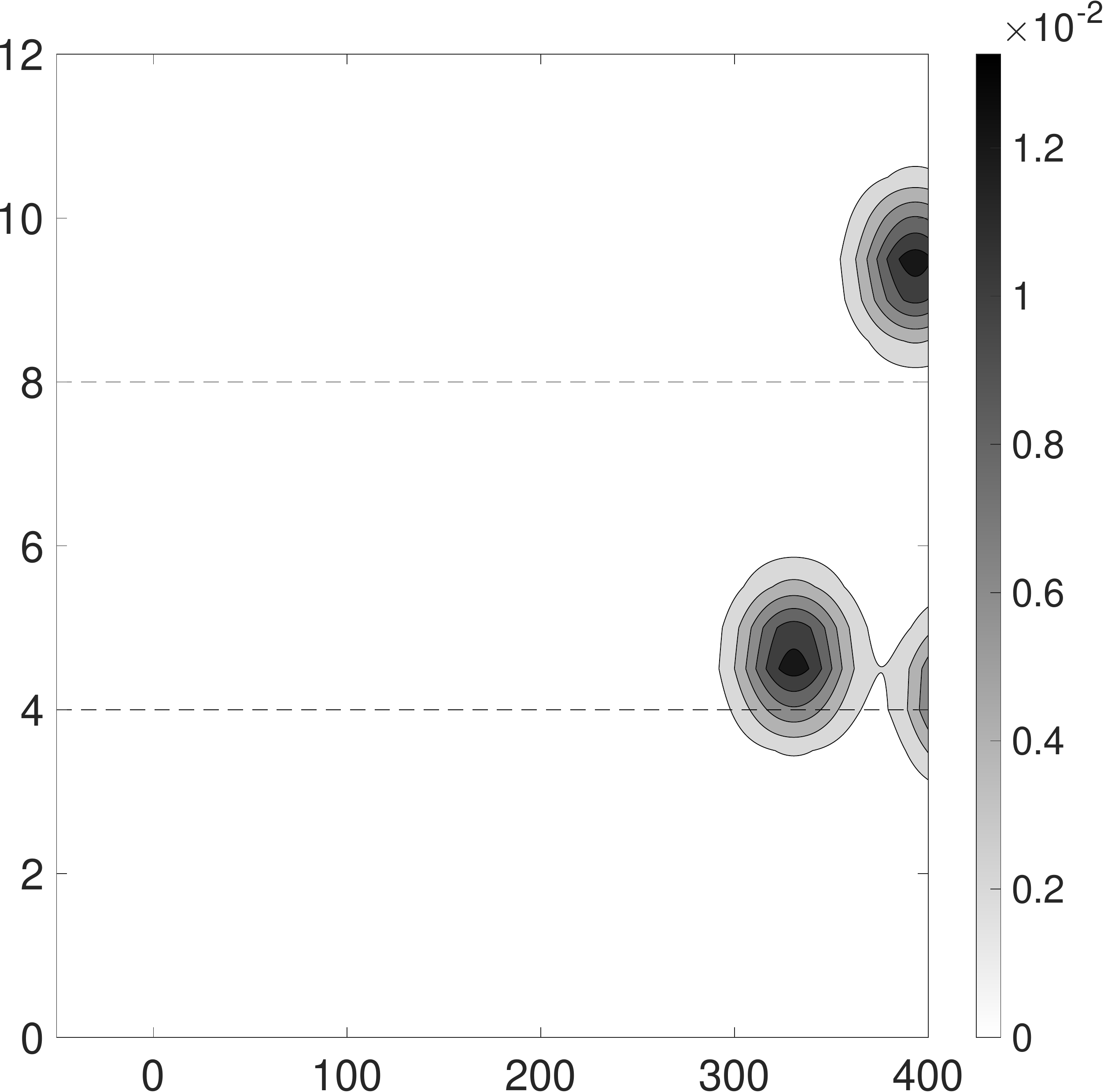}} \\
\subfloat[][Trucks density at time $t=0$.]
{\includegraphics[width=.3\columnwidth]{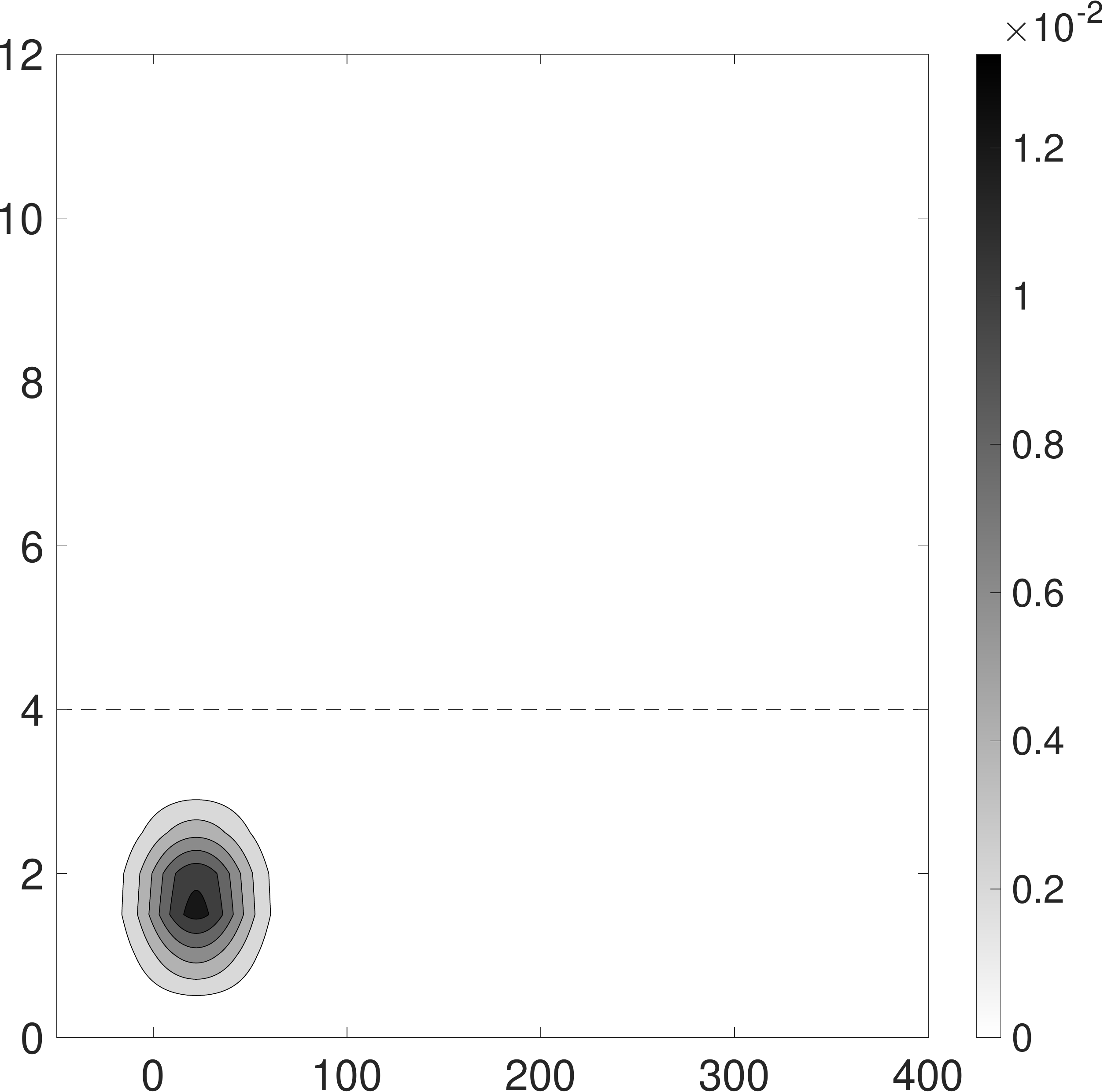}} \qquad
\subfloat[][Simulated trucks density at time $t=T$.]
{\includegraphics[width=.3\columnwidth]{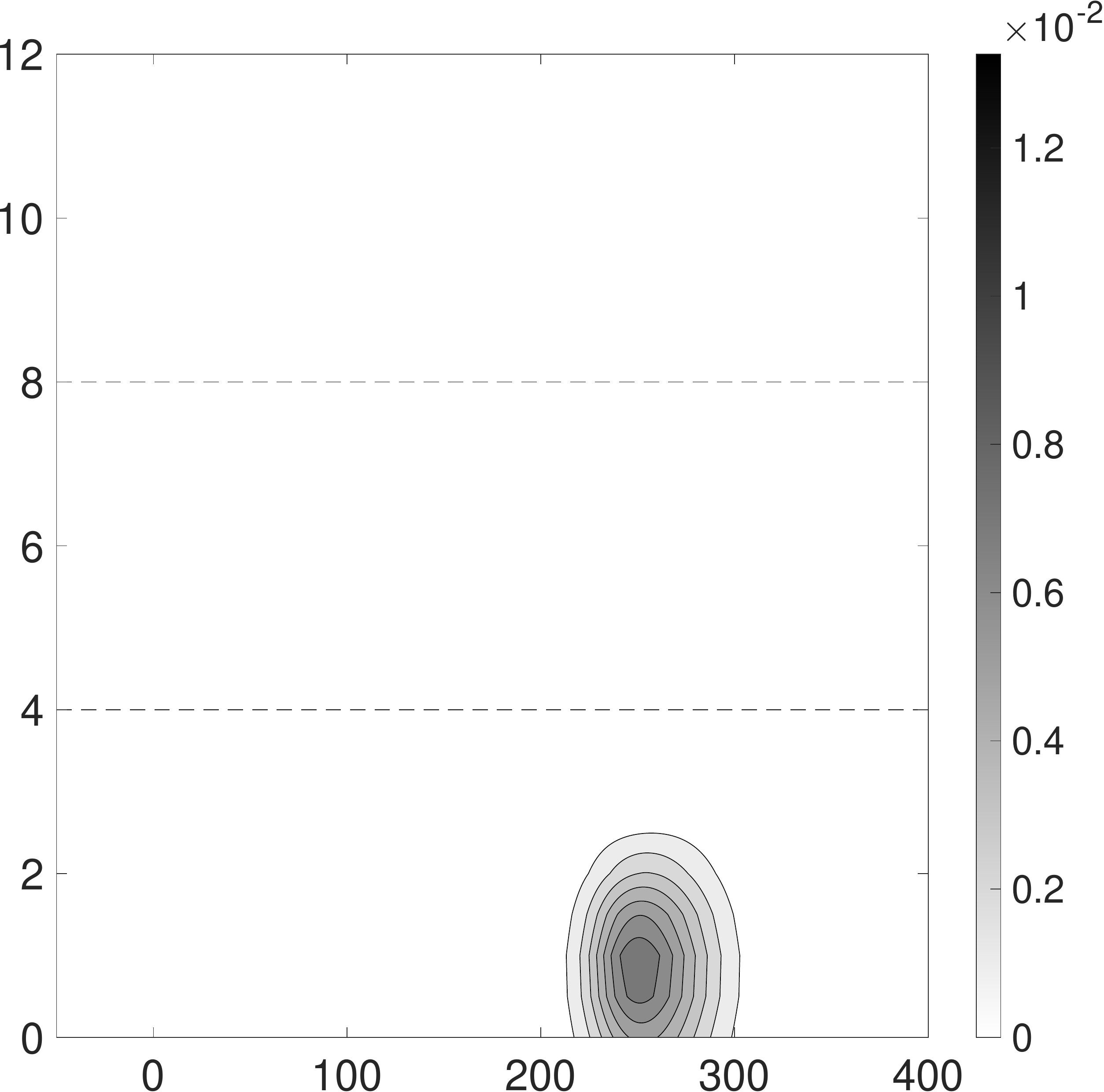}} \qquad
\subfloat[][Real trucks density at time $t=T$.]
{\includegraphics[width=.3\columnwidth]{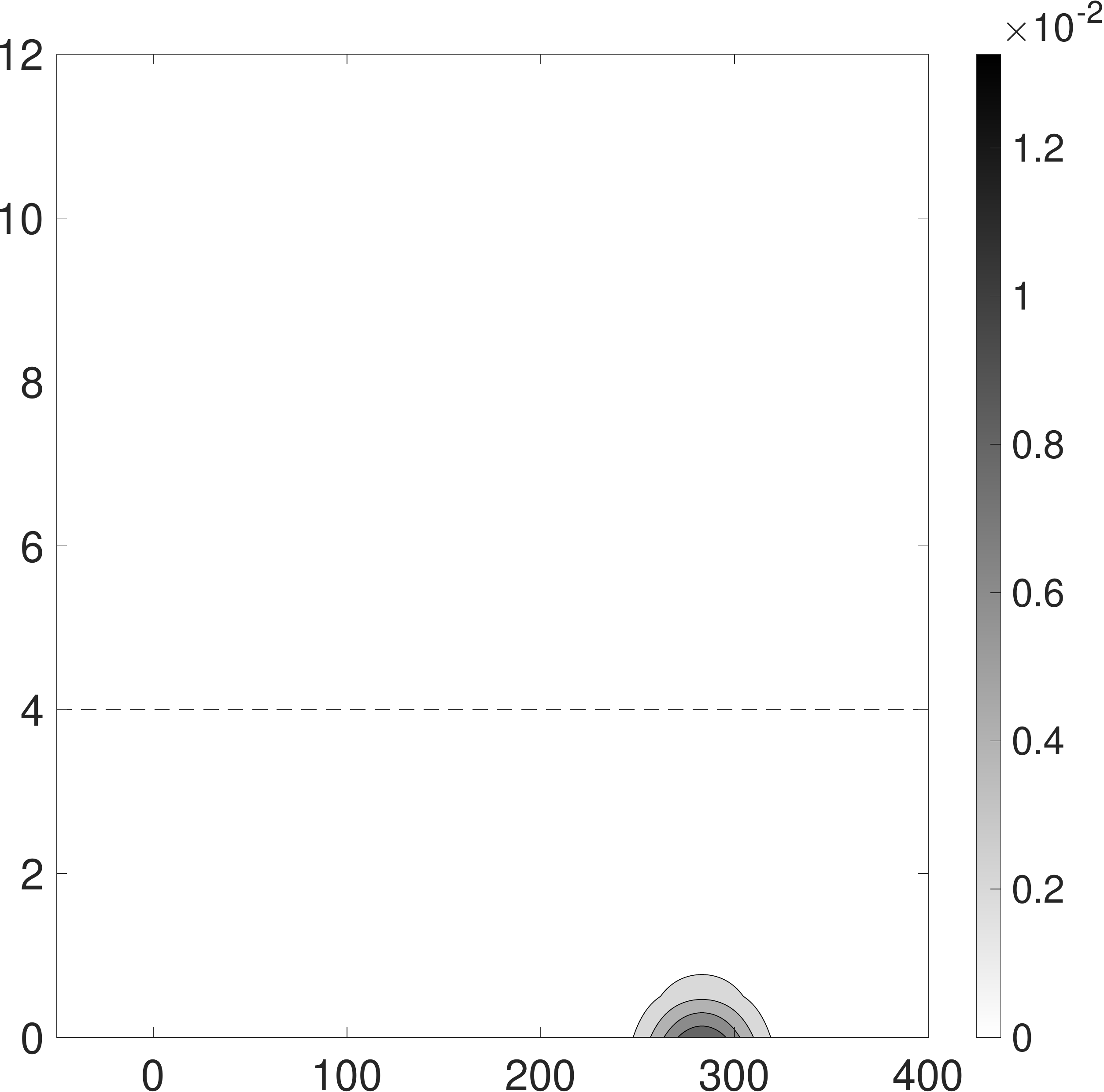}} 
\caption{Contours of the density of cars (top) and trucks (bottom): initial condition at time $t=0$ (left), simulated results at time $t=5\,\mathrm{s}$ (middle) and reconstructed real data at time $t=5\,\mathrm{s}$ (right).}
\label{fig:time1}
\end{figure}

Let us consider the domain $[0,\lx]\times[0,\ly]$ uniformly divided into a numerical grid $\Omega=[0,\nx]\times[0,\ny]$ with $x$-steps of length $\deltax$ and $y$-steps of length $\deltay$ during a time interval $[0,T]$ divided into time steps of length $\deltat$ satisfying \eqref{eq:CFL}. The numerical solutions are computed by means of the numerical scheme introduced in Section \ref{sec:numerica}, and they are denoted by
\[	\rhob^{n}_{ij}=\rhob(x_{i},y_{j},t^{n}), \qquad
	\mub^{n}_{ij}=\mub(x_{i},y_{j},t^{n})
\]
for cars and trucks respectively, with $x_{i}=i\deltax$, $y_{j}=j\deltay$ and $t^{n}=n\deltat$. The ground-truth data are estimated by \eqref{eq:KDE} as explained in Section \ref{sec:KDE}, and they are denoted by 
\[	\rhotr{,n}_{ij}=\rhotr{}(x_{i},y_{j},t^{n}), \qquad
	\mutr{,n}_{ij}=\mutr{}(x_{i},y_{j},t^{n})
\]
for cars and trucks, respectively.
The parameters used in the following test are $\lx=450\,\meter$, $\ly=14\,\meter$, $\deltax=\deltay=0.5\,\meter$, $T=5\,\second$, $h_{x}=22.5$ and $h_{y}=0.7$. The initial configuration of densities is recovered by the ground-truth data \eqref{eq:KDE} starting from the time $\hat t=14\,\second$ of the German dataset.
Therefore, at the beginning, the numerical solution coincides with the ground-truth solution. At time $\hat t$ there are three cars and one truck along the three-lanes highway, so we analyze their dynamics.
In Figure \ref{fig:time1}, we compare the contours of the ground-truth density data with the contours of the reconstructed density by the numerical simulation at the final time $T$. 

The dotted lines divide the road into three lanes. After $5\,\second$ of simulation the real and the numerical configurations of density are quite similar for both populations of vehicles.

\begin{figure}[h!]
\centering
\includegraphics[width=.5\columnwidth]{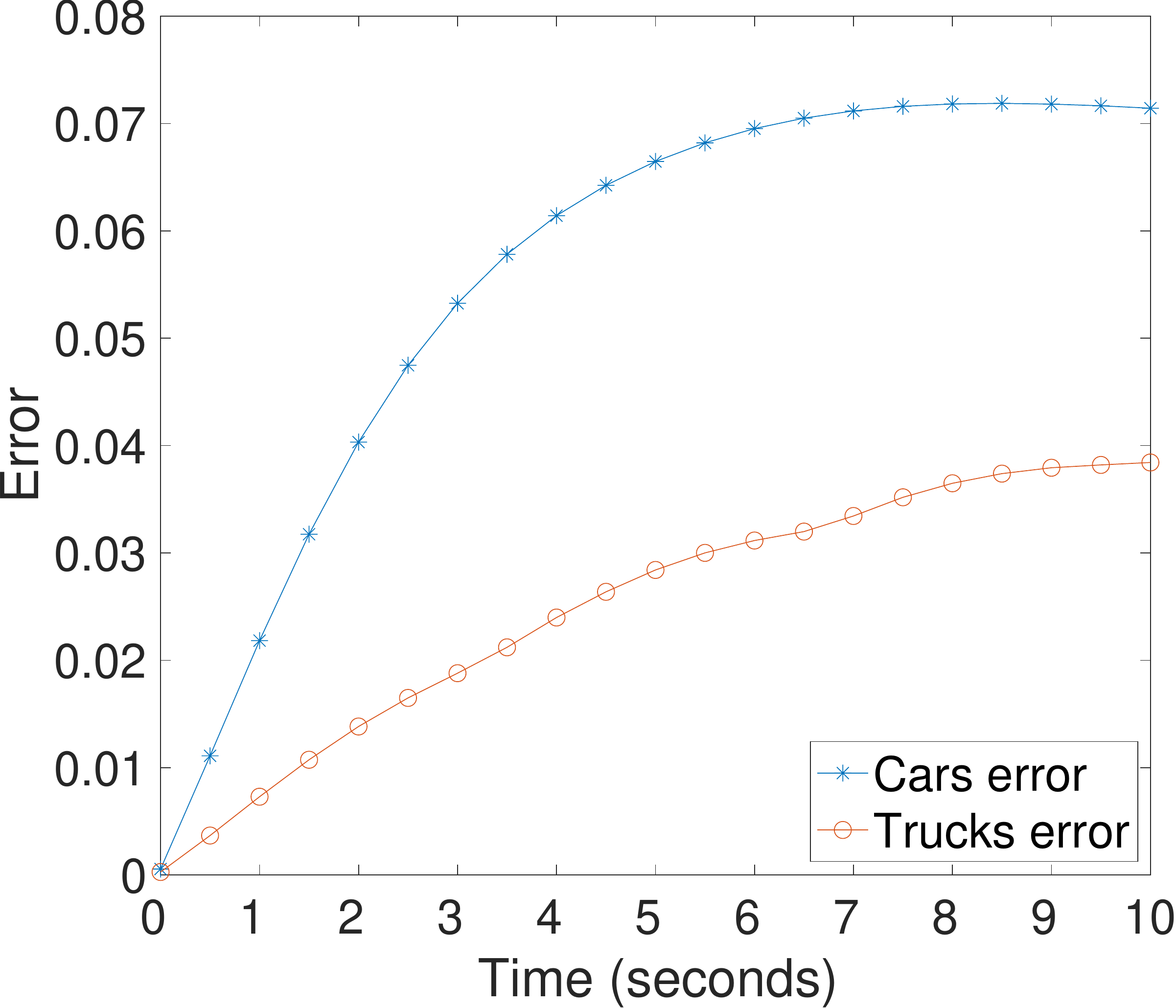} \qquad\qquad
\caption{Error between real and numerical density of cars and trucks during 10 seconds of simulation, computed every 0.5 seconds.}
\label{fig:errore}
\end{figure}

In order to better compare the numerical results with the ground-truth data, we introduce the following errors
\begin{align}
E_{\rho}(t^{n}) &= \norm{\rhotr{}(\cdot,\cdot,t^{n})-\rhob(\cdot,\cdot,t^{n})}_{L^{1}},\label{eq:err1}\\
E_{\mu}(t^{n}) &= \norm{\mutr{}(\cdot,\cdot,t^{n})-\mub(\cdot,\cdot,t^{n})}_{L^{1}}.
\label{eq:err2}
\end{align}
The errors at time $T$ of the previous simulation are $E_{\rho}(T) = 0.06$ and $E_{\mu}(T) = 0.02$ computed with \eqref{eq:err1} and \eqref{eq:err2}, respectively. 
In Figure \ref{fig:errore} we plot the numerical errors between the numerical density and real data during $10$ seconds of simulation computed every 0.5 seconds. We observe that the error related to the truck is lower than the error related to cars and that both errors increase in time. However, they remain of order $10^{-2}$.

\section{Further data analyses}\label{sec:sorpasso}
In the previous section we have seen that the two dimensional multi-class LWR \eqref{eq:lwr2D} with the flux functions defined in \eqref{eq:flussi} and calibrated with real data is able to simulate the dynamics of vehicles. However, our main assumption on the flux functions is that the coefficients $\rm$, $\cx$ and $\cy$ are equal for both the classes of vehicles. This is a strong assumption, since it implies that cars and trucks have the same length and velocity. Therefore, we modify now the definition of the flux functions to differentiate more clearly the dynamics of the two classes.

\subsection{Test with real data}
We consider again the German dataset \cite{germanDataset}. Our aim is to consider different maximum density of cars and trucks, due to the different length of vehicles, and different parameters $\cx_{\rho}$, $\cx_{\mu}$, $\cy_{\rho}$ and $\cy_{\mu}$ in order to take into account velocity functions which depend on the class of vehicles. The maximum density of cars $\rm=\rm_{\rho}=400 \,\mathrm{veh/km}$ coincides with equation \eqref{eq:rmax}, and we assume that the length of trucks is twice that of cars, hence we have $\rm_{\mu}=200 \,\mathrm{veh/km}$. 
We slightly modify the flux functions of \eqref{eq:lwr2D} as
\begin{equation}
\begin{split}
\qxr(\rho,\mu)&=\rho \cx_{\rho}\left(1-\left(\frac{\rho+2\mu}{\rm}\right)\right), \qquad 
q^{y}_{\rho}(\rho,\mu)=\rho \cy_{\rho}\left(1-\left(\frac{\rho+2\mu}{\rm}\right)\right)\\
\qxm(\rho,\mu)&=\mu \cx_{\mu}\left(1-\left(\frac{\rho+2\mu}{\rm}\right)\right),\qquad 
\qym(\rho,\mu)=\mu \cy_{\mu}\left(1-\left(\frac{\rho+2\mu}{\rm}\right)\right).
\label{eq:flussi2}
\end{split} 
\end{equation}
Note that in~\eqref{eq:flussi2} the different maximum densities between cars and trucks is expressed by the term $(\rho+2\mu)/\rm$. 

With the introduction of different coefficients $\cx_{\rho}$, $\cx_{\mu}$, $\cy_{\rho}$ and $\cy_{\mu}$ we are able to better distinguish the behavior of the two classes of vehicles, by means of different maximum velocities for the two classes in both the directions.

Next, we repeat a procedure analogous to the one proposed in Section \ref{sec:FD} to estimate the velocity functions and the fundamental diagrams. In particular, we define different velocity functions for $\rho$ and $\mu$ as
\[	\begin{split}
		\uxt_{\rho}(t_{k})&=\frac{1}{N_{\rho}(t_{k})}\sum_{i=1}^{N_{\rho}(t_{k})}v^{x}_{i},\qquad\, 
		\uyt_{\rho}(t_{k}) =\frac{1}{N_{\rho}(t_{k})}\sum_{i=1}^{N_{\rho}(t_{k})}v^{y}_{i}\\
		\uxt_{\mu}(t_{k})&=\frac{1}{N_{\mu}(t_{k})}\sum_{i=1}^{N_{\mu}(t_{k})}w^{x}_{i},\qquad
		\uyt_{\mu}(t_{k}) = \frac{1}{N_{\mu}(t_{k})}\sum_{i=1}^{N_{\mu}(t_{k})}w^{y}_{i},	
		\end{split}
\]
from which we recover the flux functions as
\begin{equation}\label{eq:q2}
	\begin{split}
		\qxrt(t_{k}) &= \rhot(t_{k})\uxt_{\rho}(t_{k}), \qquad\, 
		\qyrt(t_{k}) = \rhot(t_{k})\uyt_{\rho}(t_{k})\\  
		\qxmt(t_{k}) &= \mut(t_{k})\uxt_{\mu}(t_{k}), \qquad
		\qymt(t_{k}) = \mut(t_{k})\uyt_{\mu}(t_{k}).
	\end{split}
\end{equation}
We estimate the parameters $\cx_{\rho}$, $\cx_{\mu}$, $\cy_{\rho}$ and $\cy_{\mu}$ in order to minimize the $L^{2}$-norm between the flux functions defined in \eqref{eq:flussi2} and the fluxes derived from data in \eqref{eq:q2} and compute
\[\begin{split}
&\min_{\cx_{\rho}} \left(\norm{\qxrt-\qxr(\rho,\mu)}^{2}_{2}\right) \qquad \min_{\cy_{\rho}}\left(\norm{\qyrt-\qyr(\rho,\mu)}^{2}_{2}\right)\\
&\min_{\cx_{\mu}} \left(\norm{\qxmt-\qxm(\rho,\mu)}^{2}_{2}\right) \qquad
\min_{\cy_{\rho}}\left(\norm{\qymt-\qym(\rho,\mu)}^{2}_{2}\right),
\end{split}
\]
using again the $\mathtt{fminbnd}$ MATLAB tool. We obtain $\cx_{\rho}=99.61$, $\cy_{\rho}=-0.40$, $\cx_{\mu}=74.86$ and $\cy_{\mu}=-0.49$. Hence, the cars have a faster velocity than the trucks along the main direction of travel, while the velocity of lane-changing is quite similar between the two classes.

\begin{figure}[h!]
\centering
\subfloat[][Ground-truth speed $\uxt_{\rho}$ and family of speed functions $\uxr$ as $\mu$ changes.]{\label{fig:speedrx}
\includegraphics[width=0.22\linewidth]{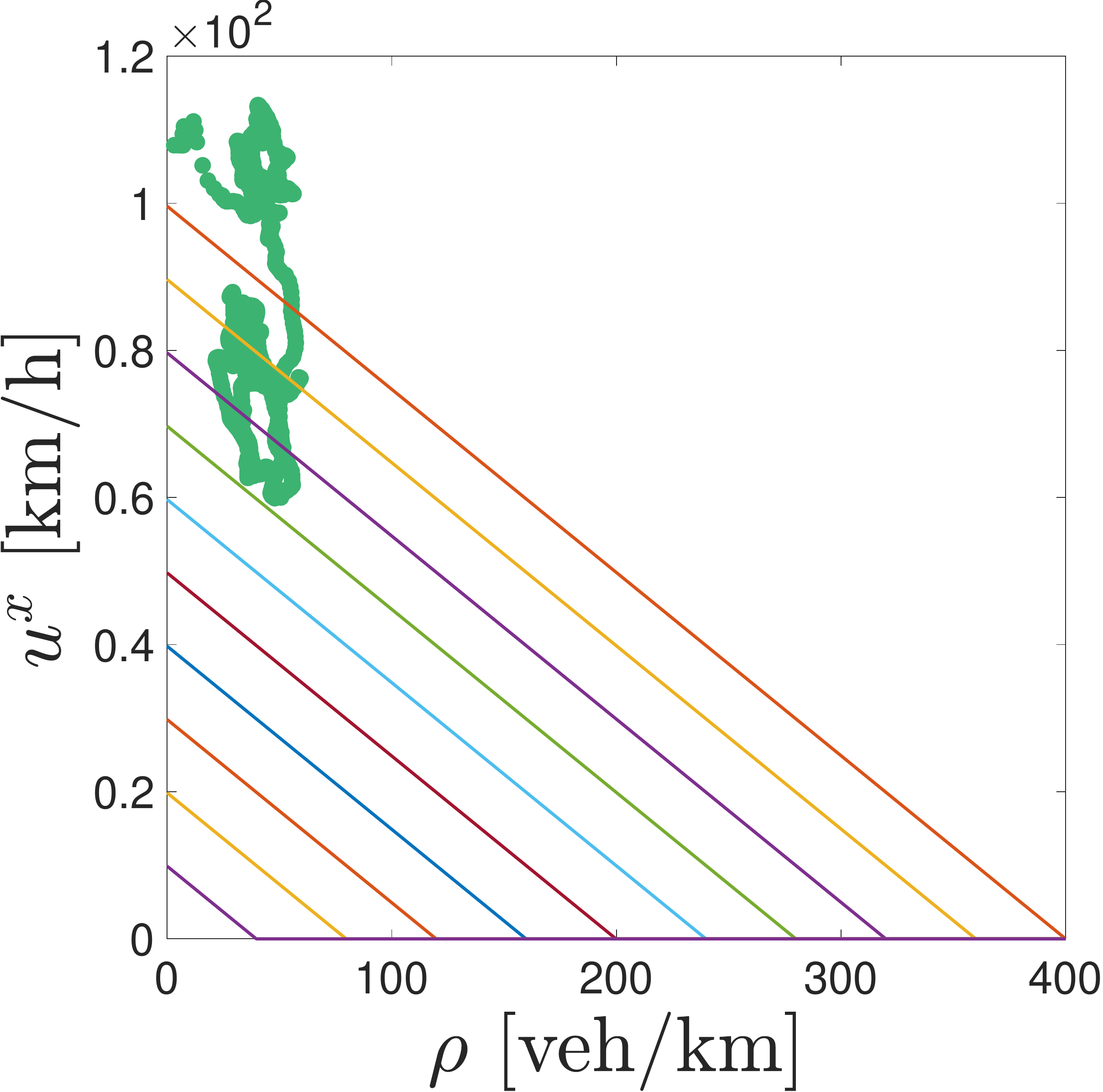}
}\,
\subfloat[][Ground-truth speed $\uxt_{\mu}$ and family of speed functions $\uxm$ as $\rho$ changes.]{\label{fig:speedmx}
\includegraphics[width=0.22\linewidth]{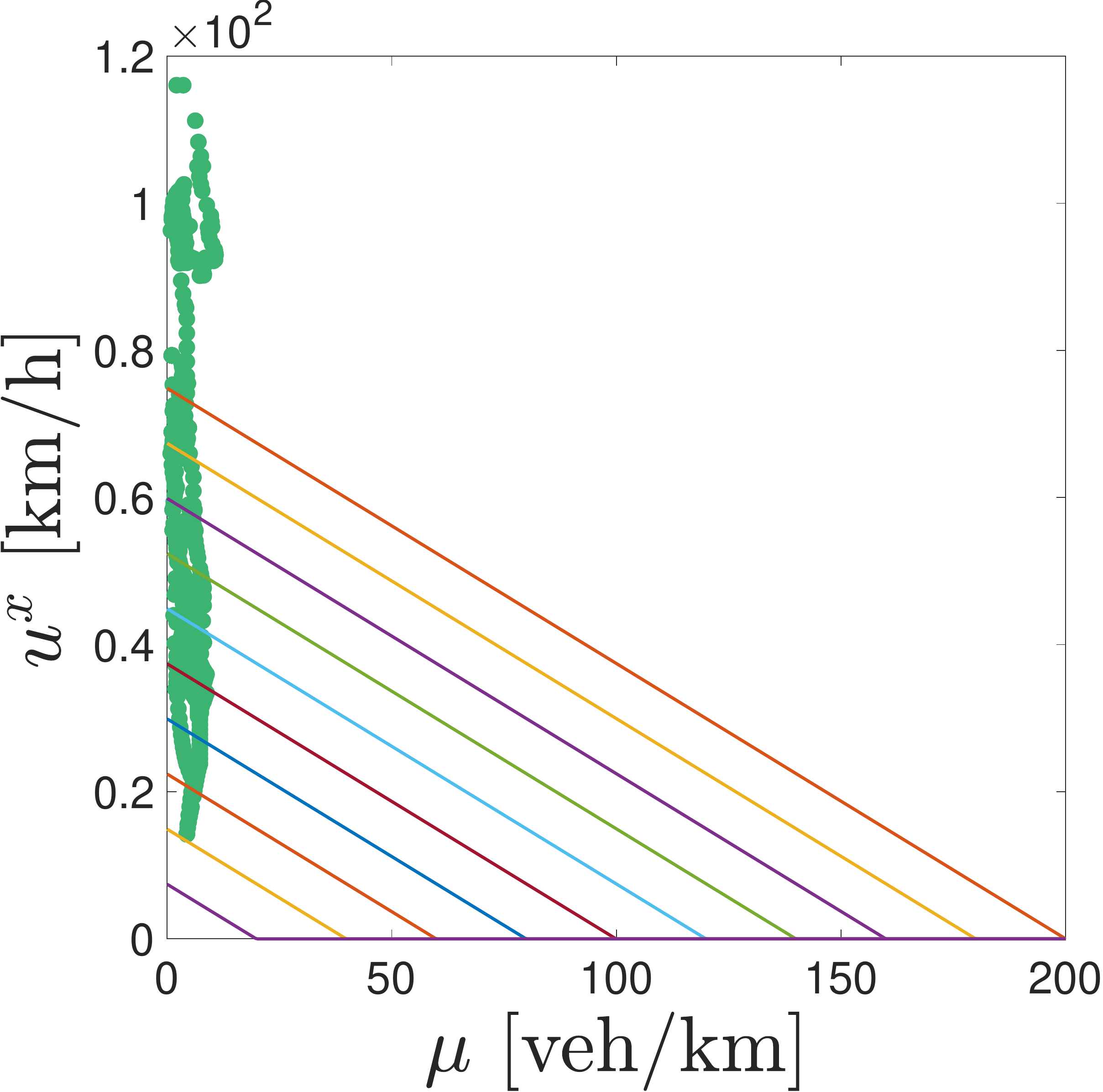}
}\,
\subfloat[][Ground-truth speed $\uyt_{\rho}$ and family of speed functions $\uyr$ as $\mu$ changes.]{\label{fig:speedry}
\includegraphics[width=0.22\linewidth]{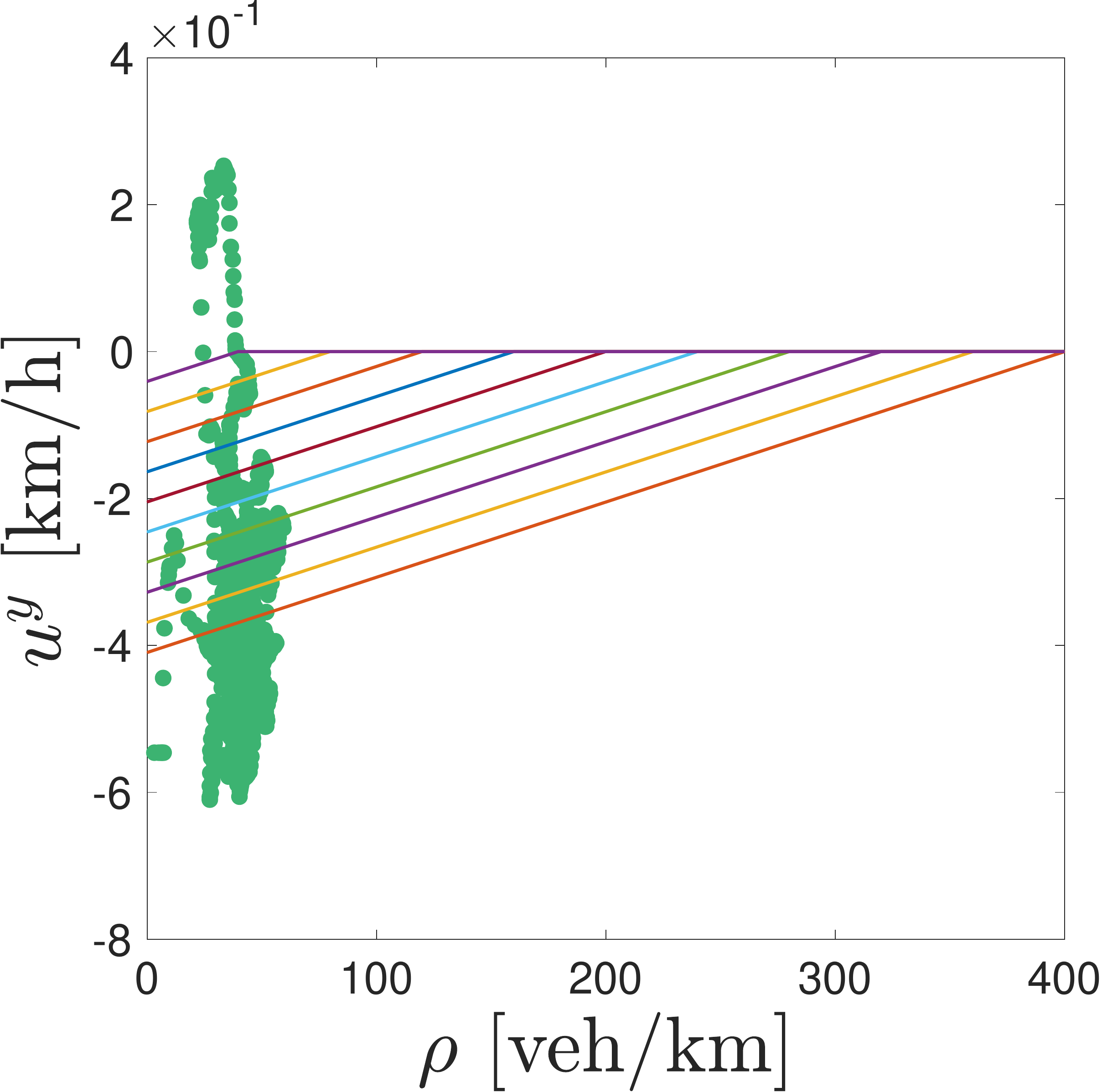}
}\,
\subfloat[][Ground-truth speed $\uyt_{\mu}$ and family of speed functions $\uym$ as $\rho$ changes.]{\label{fig:speedmy}
\includegraphics[width=0.22\linewidth]{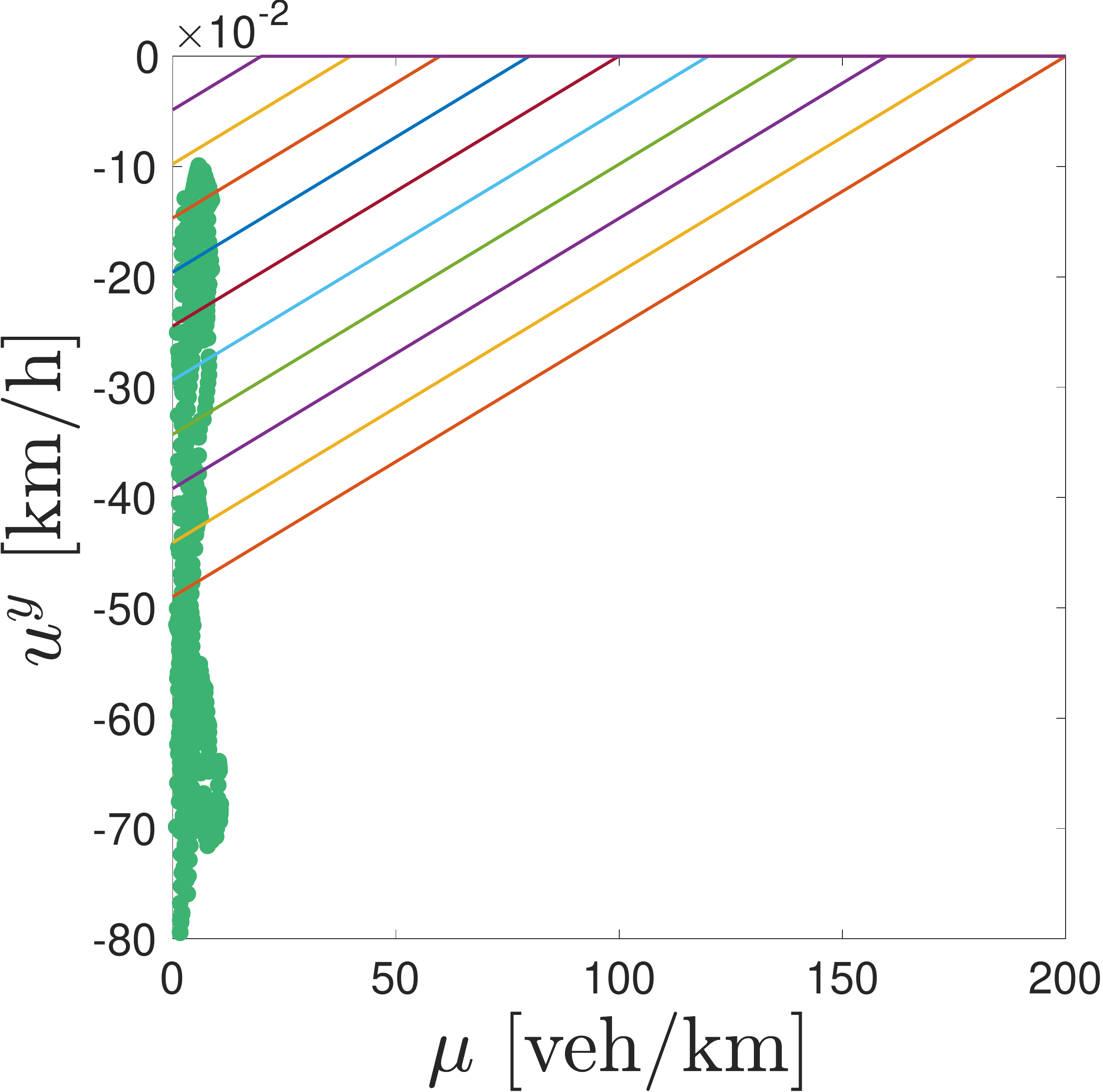}
}\\
\subfloat[][Ground-truth flux $\qxrt$ and family of flux functions $\qxr$ as $\mu$ changes.]{\label{fig:fluxrx2}
\includegraphics[width=0.22\linewidth]{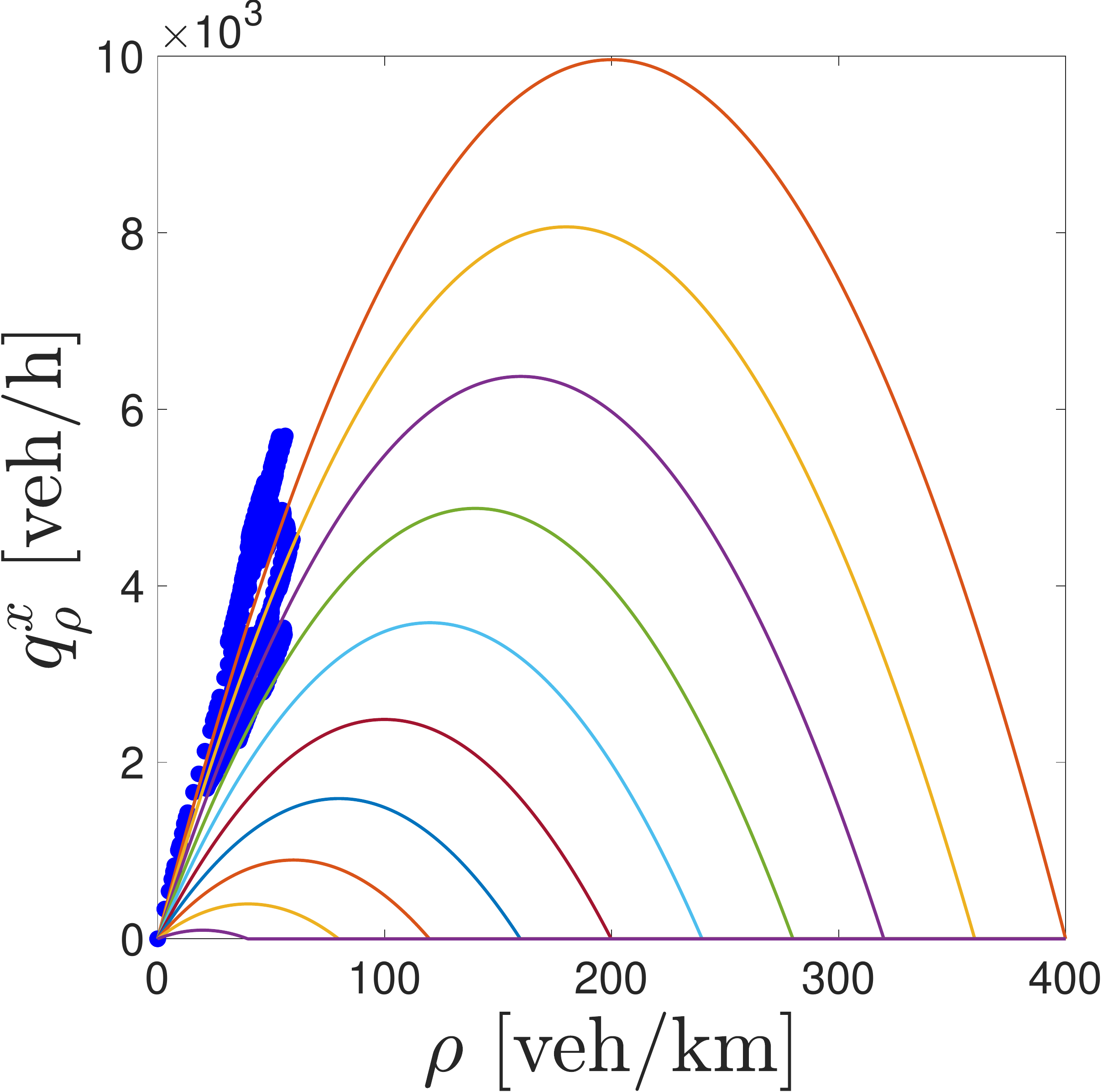}
}\,
\subfloat[][Ground-truth flux $\qxmt$ and family of flux functions $\qxm$ as $\rho$ changes.]{\label{fig:fluxmx2}
\includegraphics[width=0.22\linewidth]{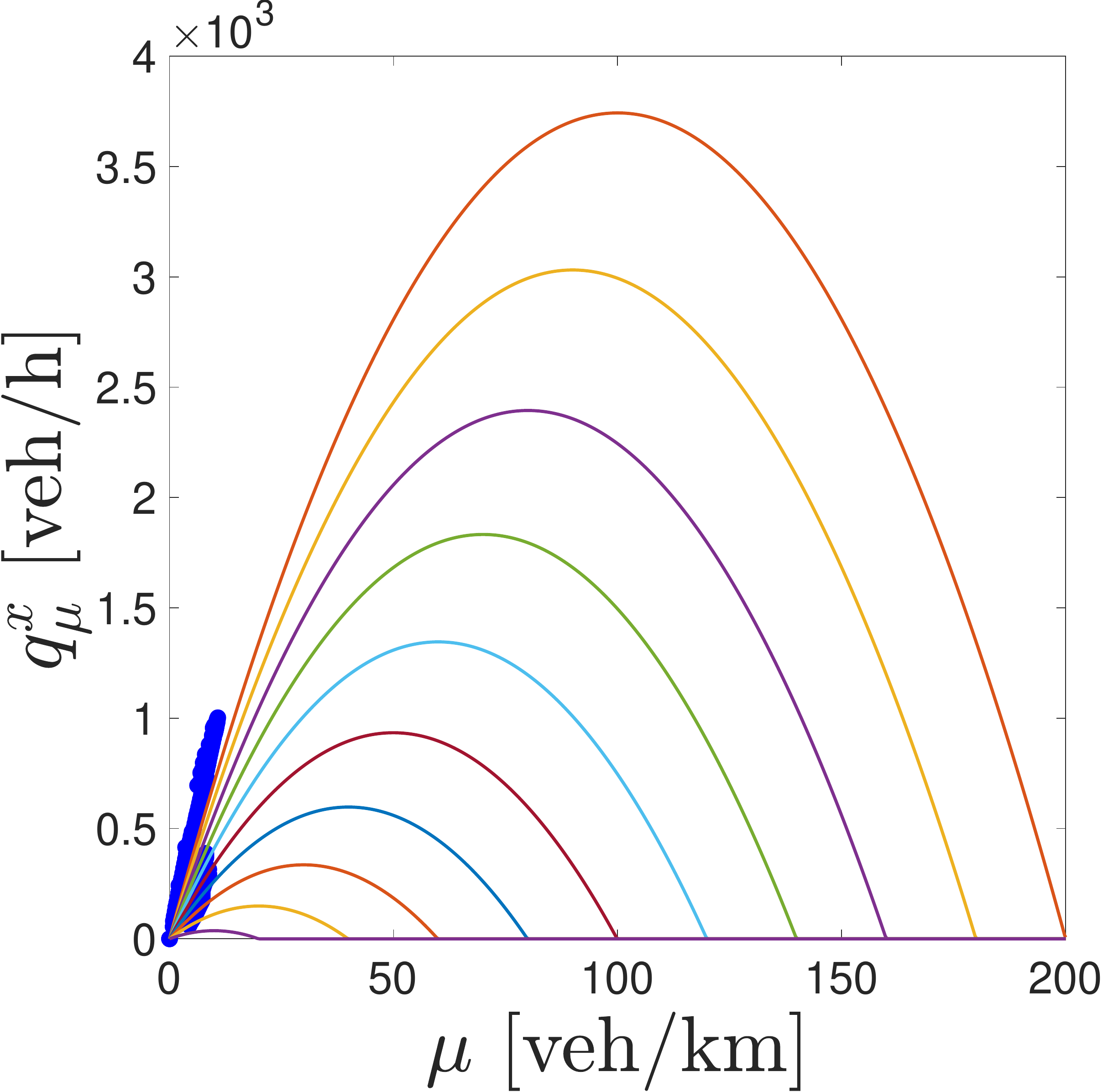}
}\,
\subfloat[][Ground-truth flux $\qyrt$ and family of flux functions $\qyr$ as $\mu$ changes.]{\label{fig:fluxry2}
\includegraphics[width=0.22\linewidth]{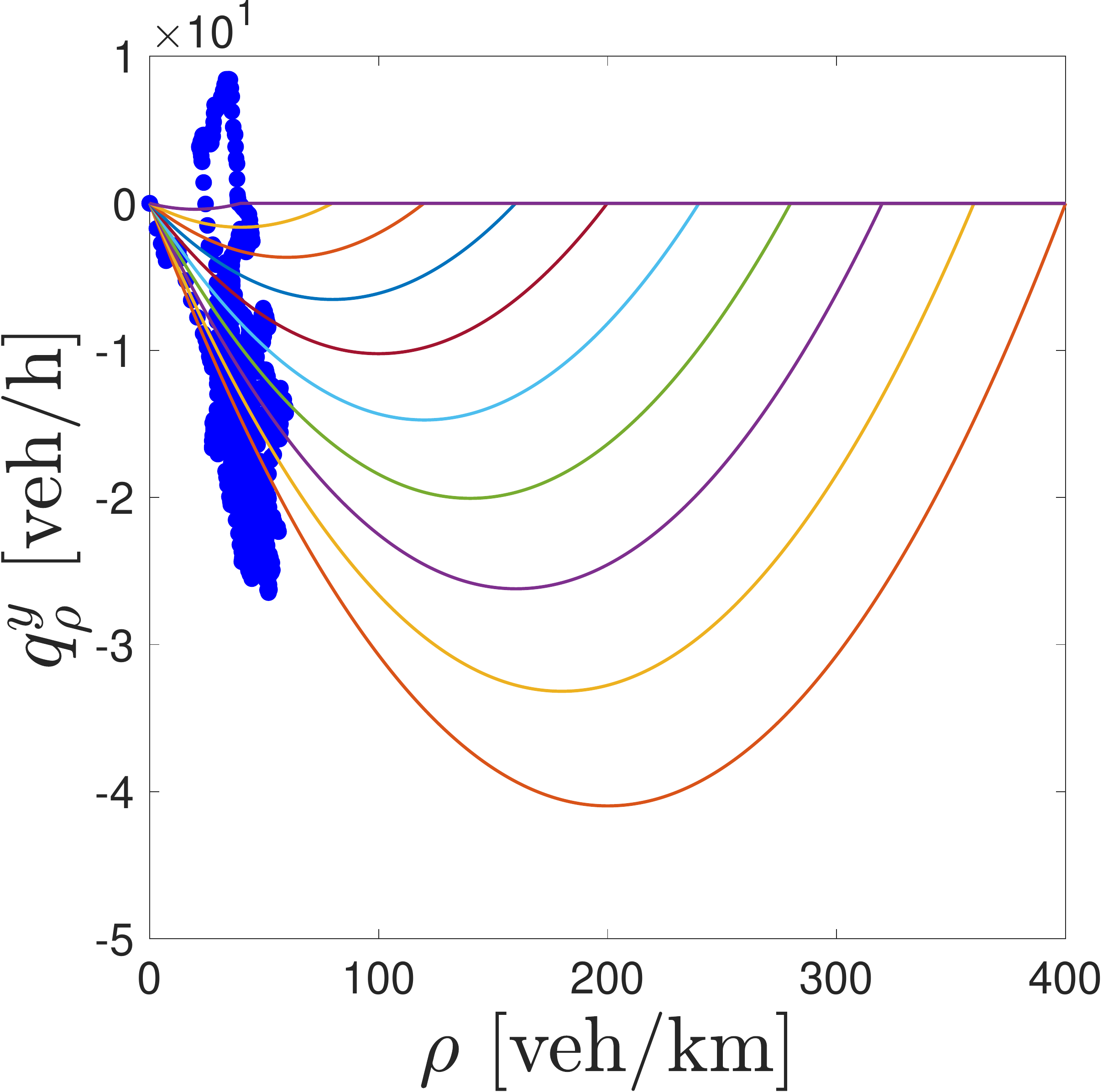}
}\,
\subfloat[][Ground-truth flux $\qymt$ and family of flux functions $\qym$ as $\rho$ changes.]{\label{fig:fluxmy2}
\includegraphics[width=0.22\linewidth]{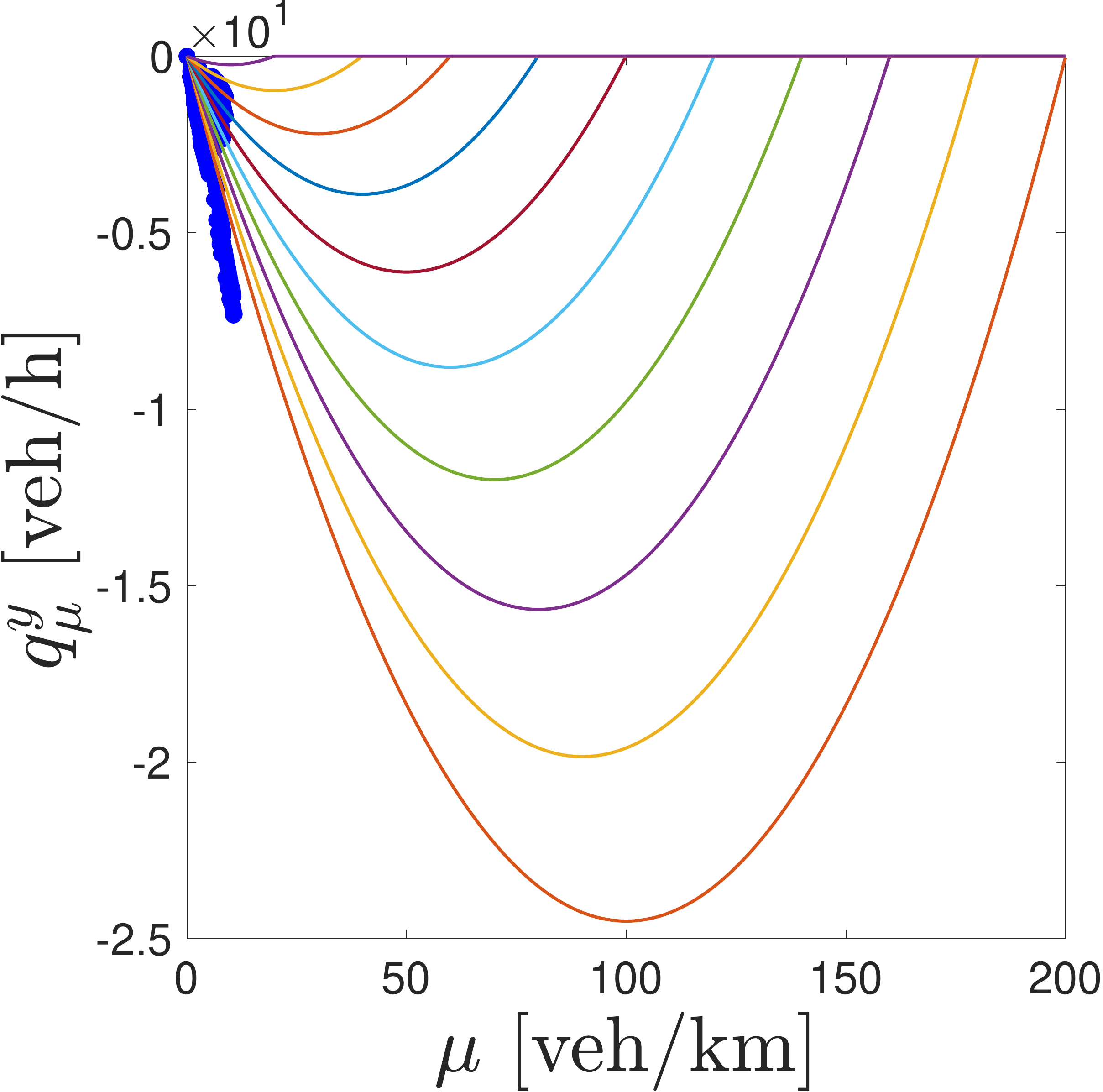}
}
\caption{Speed-density and flow-density diagrams for the two classes defined from real data (green and blue circles) and family of speed and flux functions defined by \eqref{eq:flussi2}.}
\label{fig:fluxspeed2C_2}
\end{figure}

In Figure \ref{fig:fluxspeed2C_2}, we show the family of speed and flux functions obtained with the above described procedure. Note that the speed and flux functions related to trucks, shown in Figures \subref*{fig:speedmx}, \subref*{fig:speedmy}, \subref*{fig:fluxmx2} and \subref*{fig:fluxmy2}, are defined for $\mu\in[0,\rm_{\mu}]$, with $\rm_{\mu}$ being the half of $\rm$. Again, we are able to cover the clouds of real data, but in this case the plots of $\qxm$ and $\qym$ (Figures \subref*{fig:fluxmx2} and \subref*{fig:fluxmy2}) reach lower flux values with respect to Figures \subref*{fig:flussoxm} and \subref*{fig:flussoym}, according to the lower density and velocity of trucks recorded by the dataset. Hence, the overestimation of flux values for the class of trucks is highly reduced with the introduction of $\cx_{\mu}$ and $\cy_{\mu}$ compared to the results obtained in Section \ref{sec:FD}.

\begin{figure}[h!]
\centering
\includegraphics[width=.5\columnwidth]{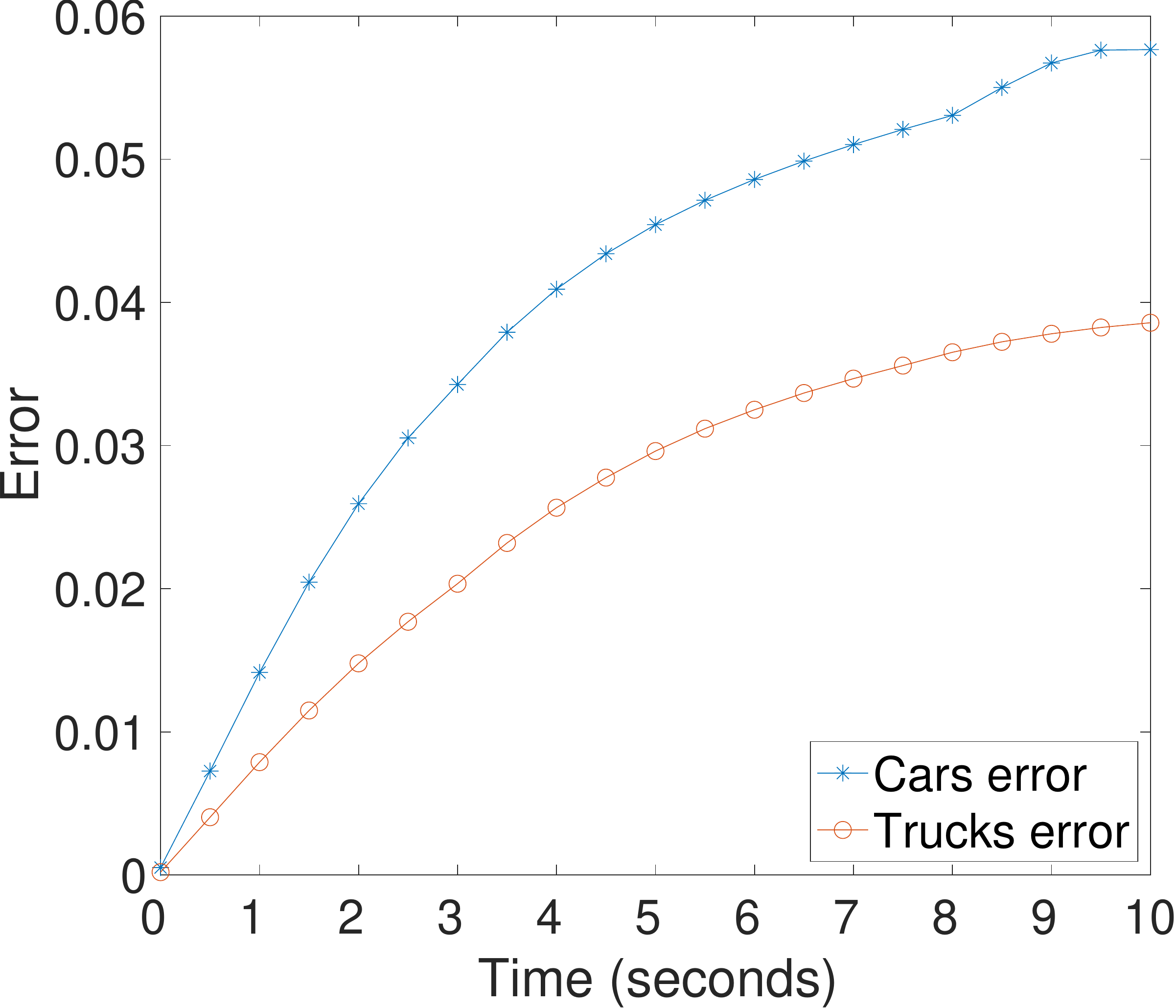} \qquad\qquad
\caption{Error between real density and numerical density of cars and trucks during 10 seconds of simulation, computed every 0.5 seconds.}
\label{fig:errore2}
\end{figure}

We repeat the same numerical test proposed in Section \ref{sec:test1} with the new flux functions \eqref{eq:flussi2}, estimating again the resulting errors with \eqref{eq:err1} and \eqref{eq:err2}. The density plots we obtain are similar to the plots shown in Figure \ref{fig:time1}, thus we omit the picture. However, as shown in Figure \ref{fig:errore2}, we obtain a better estimate of the errors compared to the test done in Section \ref{sec:test1}.   

\subsection{Vehicles overtaking}\label{sec:overtaking}
A further investigation of our model is the testing of the ability of capturing vehicles overtaking. We consider the following flux functions
%
\[\begin{split}
\qxr(\rho,\mu)&=\rho \cx_{\rho}\left(1-\left(\frac{\rho+\mu}{\rm}\right)\right), \qquad 
q^{y}_{\rho}(\rho,\mu)=\rho \cy_{\rho}\left(1-\left(\frac{\rho+\mu}{\rm}\right)\right)\\
\qxm(\rho,\mu)&=\mu \cx_{\mu}\left(1-\left(\frac{\rho+\mu}{\rm}\right)\right),\qquad 
\qym(\rho,\mu)=\mu \cy_{\mu}\left(1-\left(\frac{\rho+\mu}{\rm}\right)\right),
\end{split} 
\]
where we have different parameters $\cx_{\rho}$, $\cx_{\mu}$, $\cy_{\rho}$ and $\cy_{\mu}$ but the same maximum density $\rm$. The idea is to simulate traffic dynamics with different maximum velocities for the two classes and verify if the faster vehicles are able to overtake the slower ones. Indeed, the presence of the component transverse to the main direction of motion naturally lends itself to the modeling of vehicles overtaking. 

We consider a numerical grid $\Omega=[0,\nx]\times[0,\ny]$ with $x$-steps $\deltax$ and $y$-steps $\deltay$ during a time interval $[0,T]$ divided into time steps $\deltat$ satisfying \eqref{eq:CFL}. In particular we work on a road with two lanes, with two cars and one truck. We fix the following parameters: $\lx=100\,\meter$, $\ly=6\,\meter$, $\deltax=\deltay=0.2\,\meter$ and $T=4\,\second$.
Moreover, we assume that $\cx_{\rho}=80$ and $\cy_{\rho}=-0.4$, while $\cx_{\mu}=\cy_{\mu}=0$, thus the truck does not move.

\begin{figure}[h!]
\centering
\subfloat[][Cars and truck density at time $t=0$.]
{ \begin{overpic}[abs,unit=0.5mm,width=.3\columnwidth]{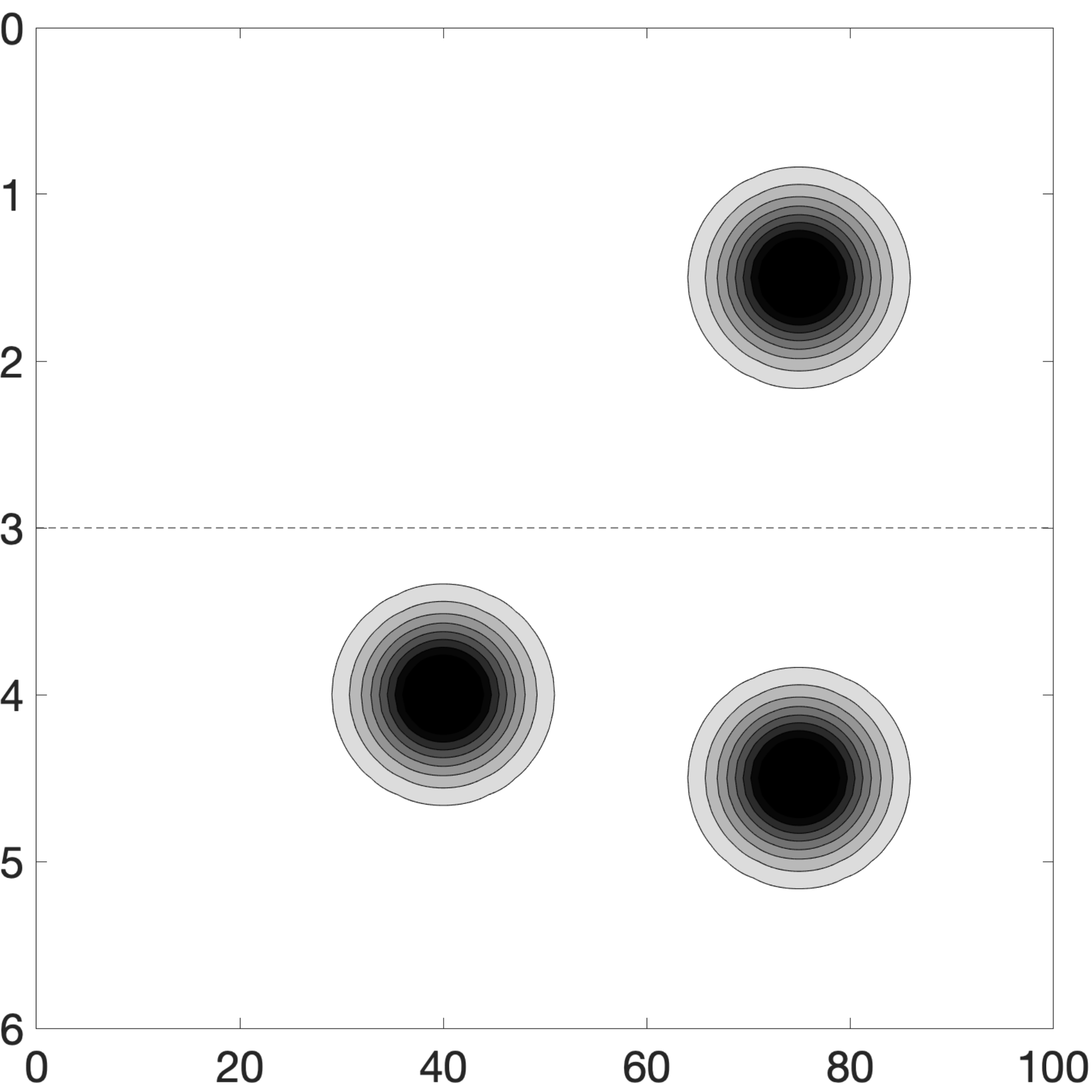}
\put(30,16){\small Car 2}
\put(60,8){\small Truck}
\put(62,53){\small Car 1}
\end{overpic}
\label{fig:test2a}} \quad
\subfloat[][Cars and truck density at time $t=T/4$.]
{\includegraphics[width=.3\columnwidth]{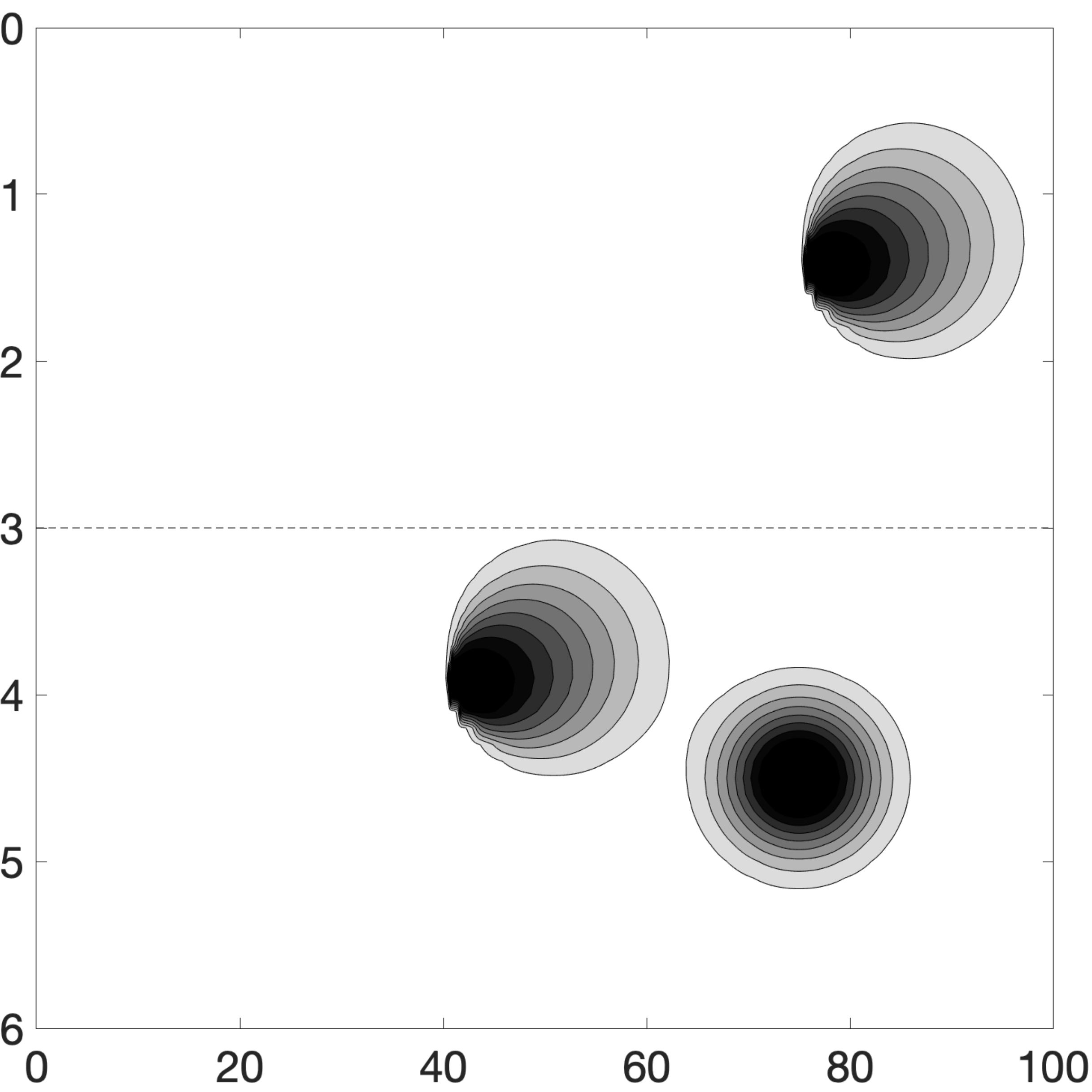}\label{fig:test2b}} \quad
\subfloat[][Cars and truck density at time $t=T$.]
{\includegraphics[width=.3\columnwidth]{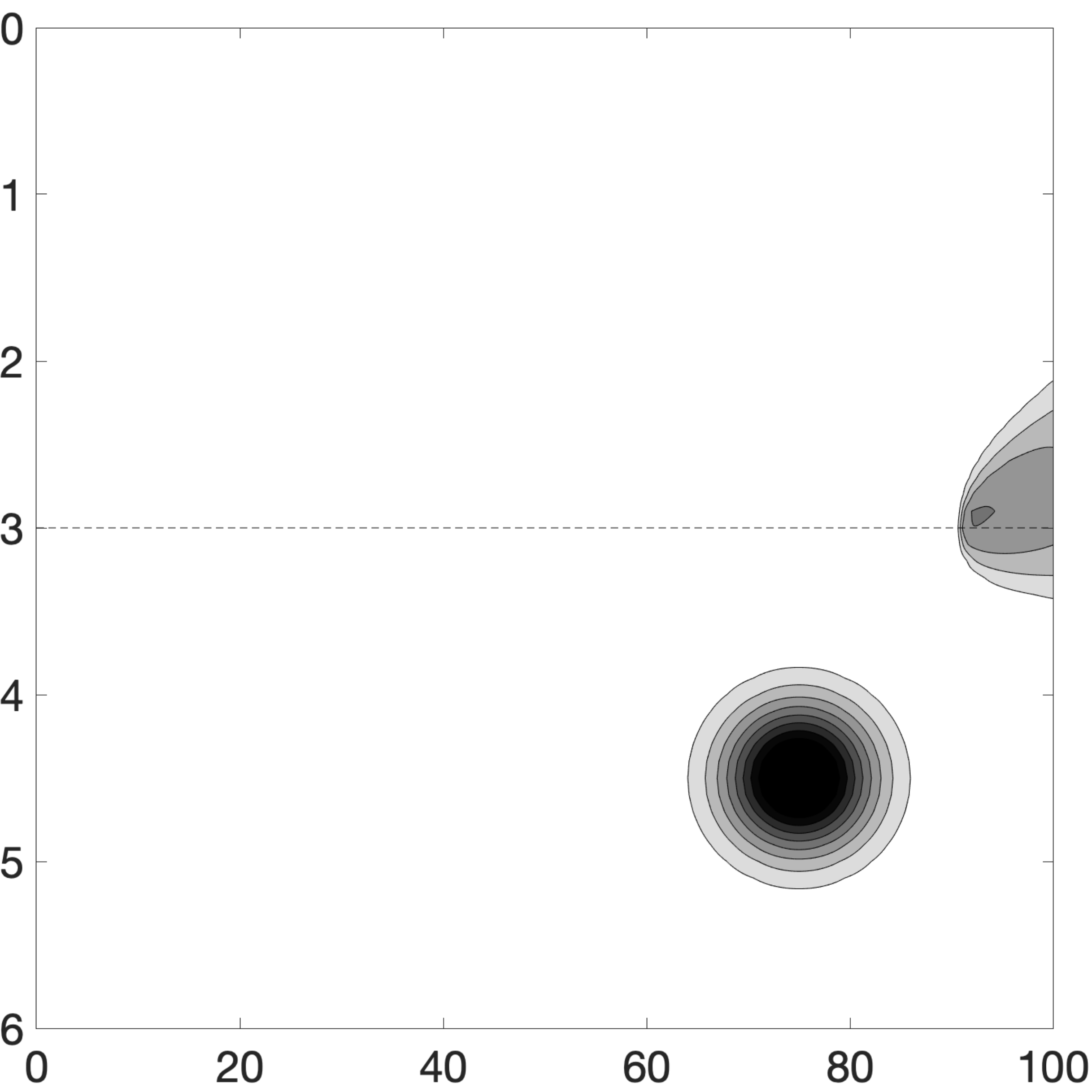}\label{fig:test2c}} \caption{Contours of the density of cars and truck at time $t=0$ (left), $t=T/4$ (middle) and $t=T$ (right). The truck does not move, Car 1 leaves the road during the simulation and Car 2 overtakes the truck and is exiting the road at time $T$ .}
\label{fig:test2}
\end{figure}

As we can see in Figure \ref{fig:test2}, we consider a road with two lanes with two cars and a truck. At the beginning of the simulation there is a car in the top lane and a truck in front of the other car in the bottom lane. Since the truck does not move and the cars are free to move along the $y$-axis, in Figure \subref*{fig:test2b} we see that both cars move towards the north-east direction. In particular, the car in the top lane is leaving the road and the other one starts to overtake the truck, which acts as an obstacle along the main travel direction. Finally, Figure \subref*{fig:test2c} shows that the car has been able to overtake the truck since it is exiting the road while the truck is still inside the domain.

\subsubsection{Comparison with a multi-lane model}\label{sec:multilane}

Finally, we intend to compare the proposed 2D multi-class approach to a first order multi-class multi-lane model. Indeed, the inclusion of lane change dynamics well fits with the simulation of vehicle overtaking and allows us to compare the results obtained with the 2D multi-class model.
Specifically, we extend the multi-lane LWR model proposed in \cite{holden2019SIAM} to a multi-class model. 
Hence, let us consider a road with two lanes and two classes of vehicles $\rho$ and $\mu$. The dynamics on the two lanes is described by
\begin{equation}\label{eq:multilane}
\begin{split}
	\text{Lane 1:}\ &\ \begin{cases}
		\rho^{1}_{t}+(q_{\rho}^{1}(\rho^{1},\mu^{1}))_{x}=-S_{\rho}(\rho^{1},\mu^{1},\rho^{2},\mu^{2})\\
		\mu^{1}_{t}+(q_{\mu}^{1}(\rho^{1},\mu^{1}))_{x}=-S_{\mu}(\rho^{1},\mu^{1},\rho^{2},\mu^{2})\\
	\end{cases}\\
	\text{Lane 2:}\ &\ \begin{cases}
		\rho^{2}_{t}+(q_{\rho}^{2}(\rho^{2},\mu^{2}))_{x}=S_{\rho}(\rho^{1},\mu^{1},\rho^{2},\mu^{2})\\
		\mu^{2}_{t}+(q_{\mu}^{2}(\rho^{2},\mu^{2}))_{x}=S_{\mu}(\rho^{1},\mu^{1},\rho^{2},\mu^{2}),\\
	\end{cases}
\end{split}
\end{equation}
where $\rho^{1}$, $\mu^{1}$ and $q^{1}_{\rho,\mu}$ are the densities and the flux function of the two classes along lane 1, and $\rho^{2}$, $\mu^{2}$ and $q^{2}_{\rho,\mu}$ along lane 2. The functions $S_{\rho,\mu}$ regulate the lane changing and are defined as
\begin{equation}\label{eq:Smultilane}
\begin{split}
	S_{\rho} &= C(\max\{u_{\rho}^{2}(\rho^{2},\mu^{2})-u_{\rho}^{1}(\rho^{1},\mu^{1}),0\}\rho^{1}+\min\{u_{\rho}^{2}(\rho^{2},\mu^{2})-u_{\rho}^{1}(\rho^{1},\mu^{1}),0\}\rho^{2})\\
	S_{\mu} &= C(\max\{u_{\mu}^{2}(\rho^{2},\mu^{2})-u_{\mu}^{1}(\rho^{1},\mu^{1}),0\}\mu^{1}+\min\{u_{\mu}^{2}(\rho^{2},\mu^{2})-u_{\mu}^{1}(\rho^{1},\mu^{1}),0\}\mu^{2}),
\end{split}
\end{equation}
where $u^{1,2}_{\rho,\mu}$ are the velocity functions related to $\rho$ and $\mu$ respectively along the two lanes and $C$ is a constant.

The flux functions are chosen similar to the ones used for the two-dimensional multi-class model. Therefore, 
we define them as
\[\begin{split}
q^{1,2}_{\rho} &= \rho^{1,2} c_{\rho}\left(1-\left(\frac{\rho^{1,2}+\mu^{1,2}}{\rm}\right)\right), \qquad q^{1,2}_{\mu} = \mu^{1,2} c_{\mu}\left(1-\left(\frac{\rho^{1,2}+\mu^{1,2}}{\rm}\right)\right), 
\end{split}\]
where $\rm$ is the maximum density of the two classes. In order to compare such a model with the results obtained with our multi-class two-dimensional model from Section \ref{sec:overtaking}, we replicate an analogous test. Indeed, we consider a road $[0,L]$ with two lanes along which there are two cars and a truck during a time interval $[0,T]$. The parameters of the test are $L=100\,\meter$, $\deltax=0.2\,\meter$, $T=4\,\second$, $C=1$ in \eqref{eq:Smultilane}, $c_{\rho}=80$ and $c_{\mu}=0$. We use a Godunov scheme \cite{godunov1959MS} to approximate problem \eqref{eq:multilane}. 

\begin{figure}[h!]
\centering
\subfloat[][Cars and truck density at time $t=0$ on Lane 1.]
{\includegraphics[width=.3\columnwidth]{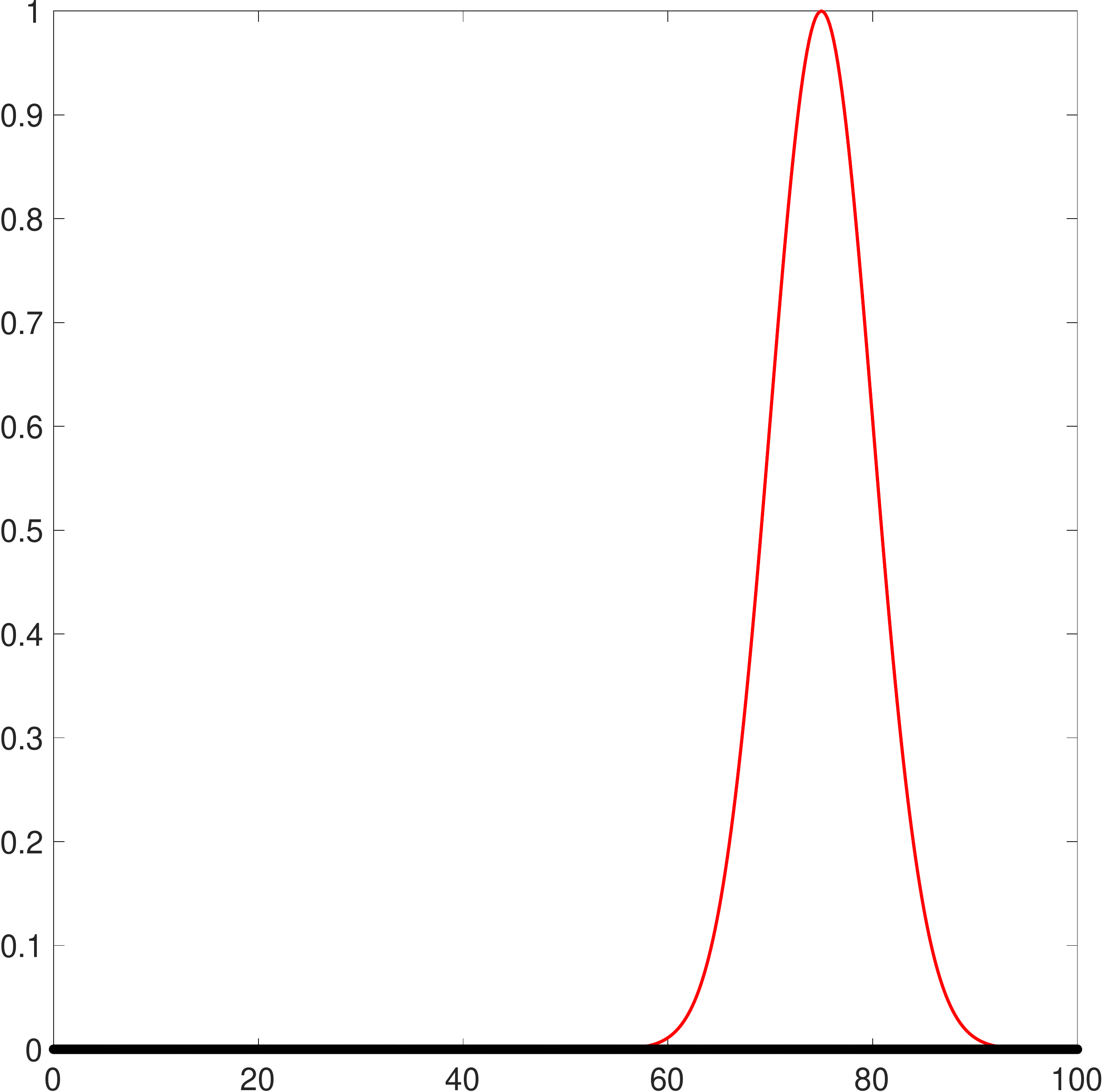}\label{fig:test3a}} \qquad
\subfloat[][Cars and truck density at time $t=T/4$  on Lane 1.]
{\includegraphics[width=.3\columnwidth]{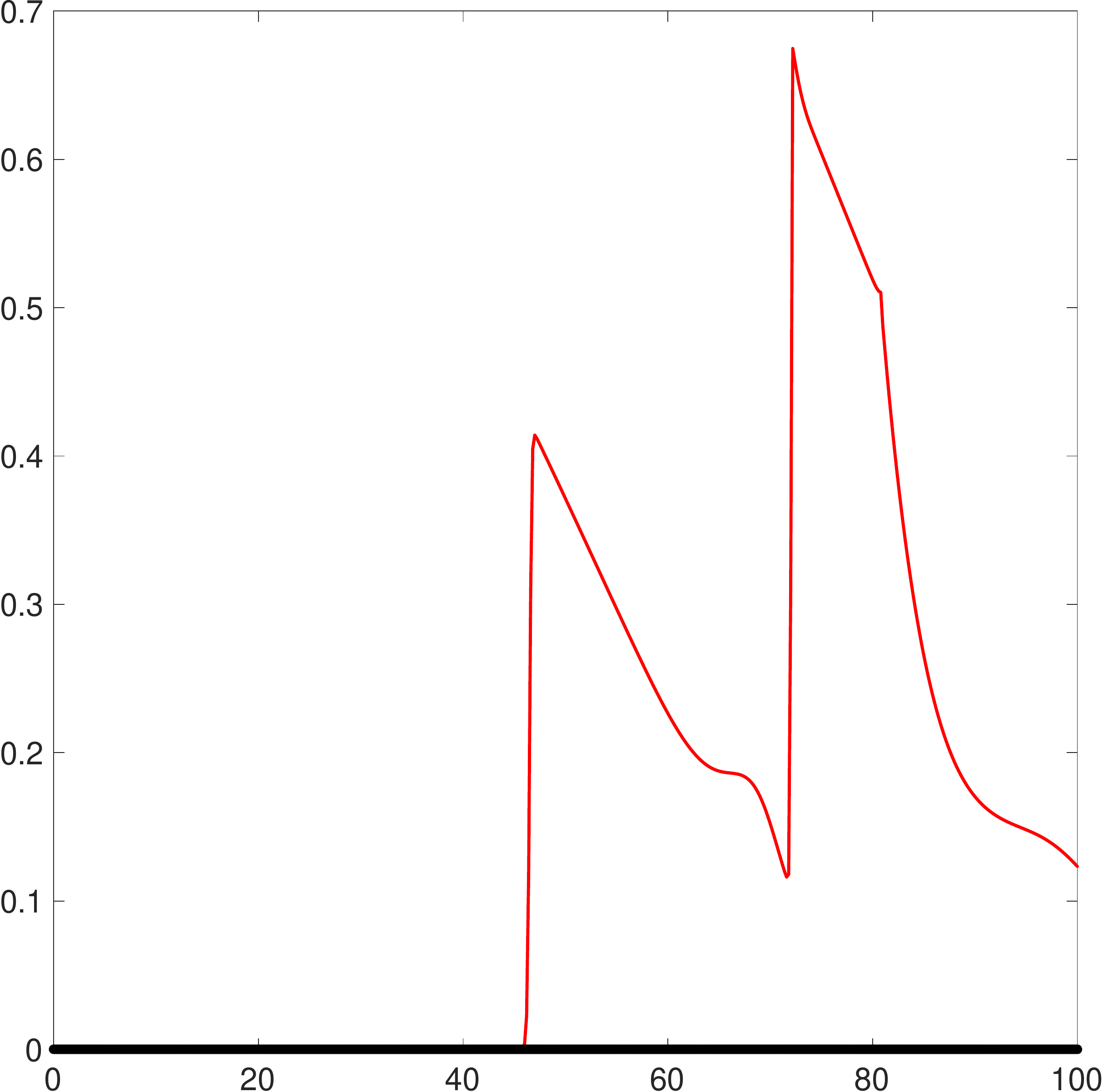}\label{fig:test3b}} \qquad
\subfloat[][Cars and truck density at time $t=T$  on Lane 1.]
{\includegraphics[width=.3\columnwidth]{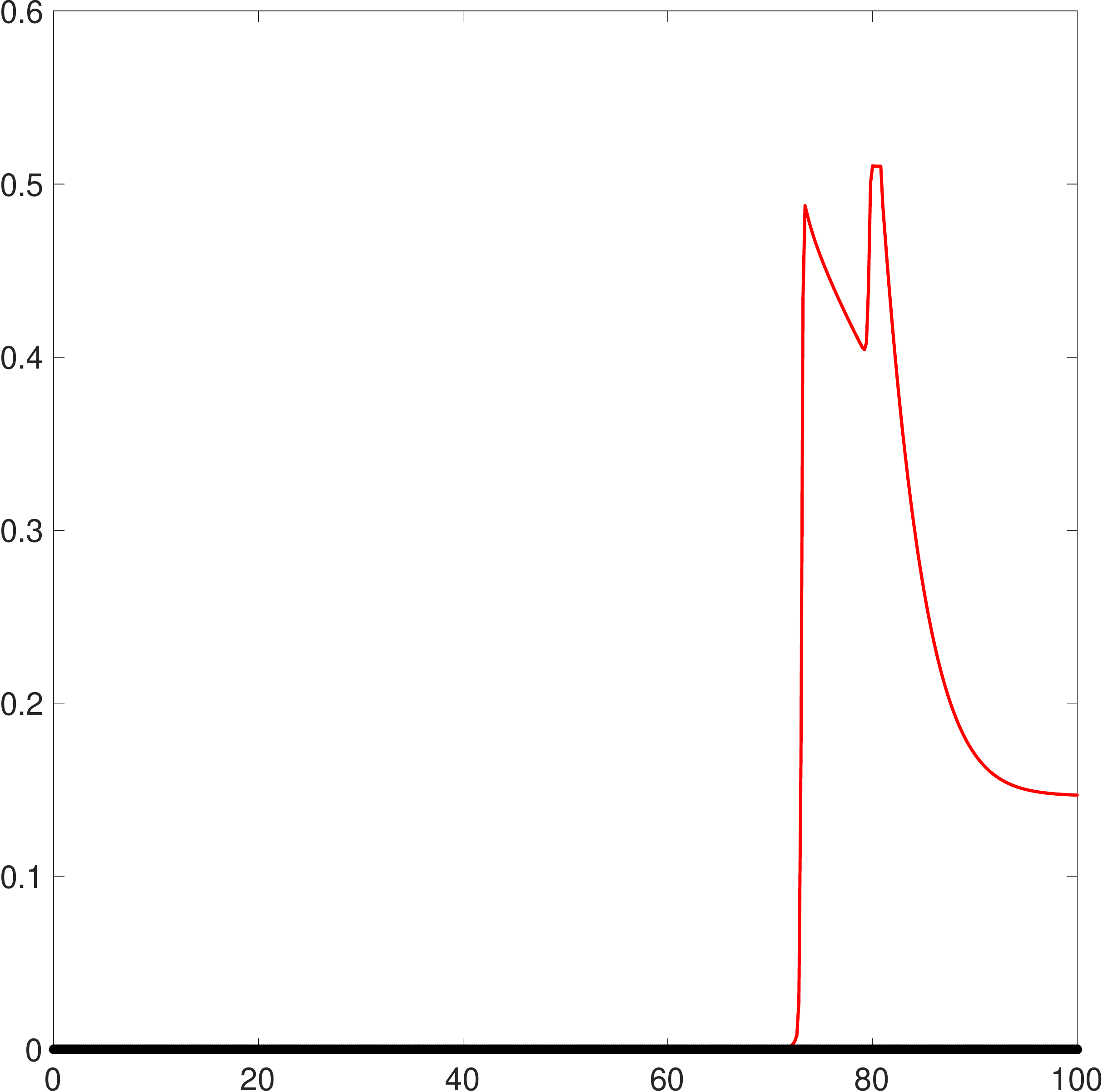}\label{fig:test3c}}\\
\subfloat[][Cars and truck density at time $t=0$  on Lane 2.]
{\includegraphics[width=.3\columnwidth]{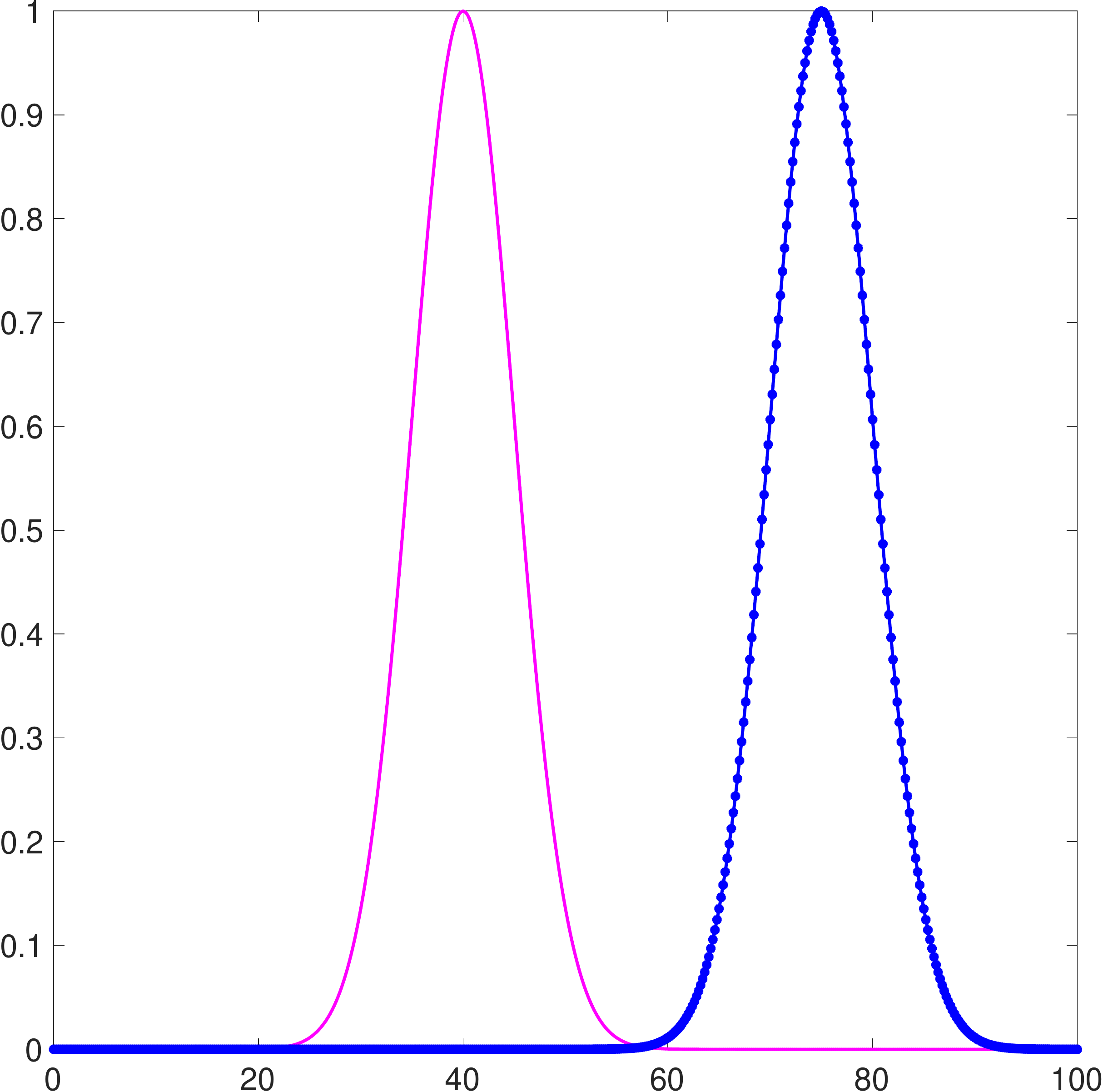}\label{fig:test3d}} \qquad
\subfloat[][Cars and truck density at time $t=T/4$  on Lane 2.]
{\includegraphics[width=.3\columnwidth]{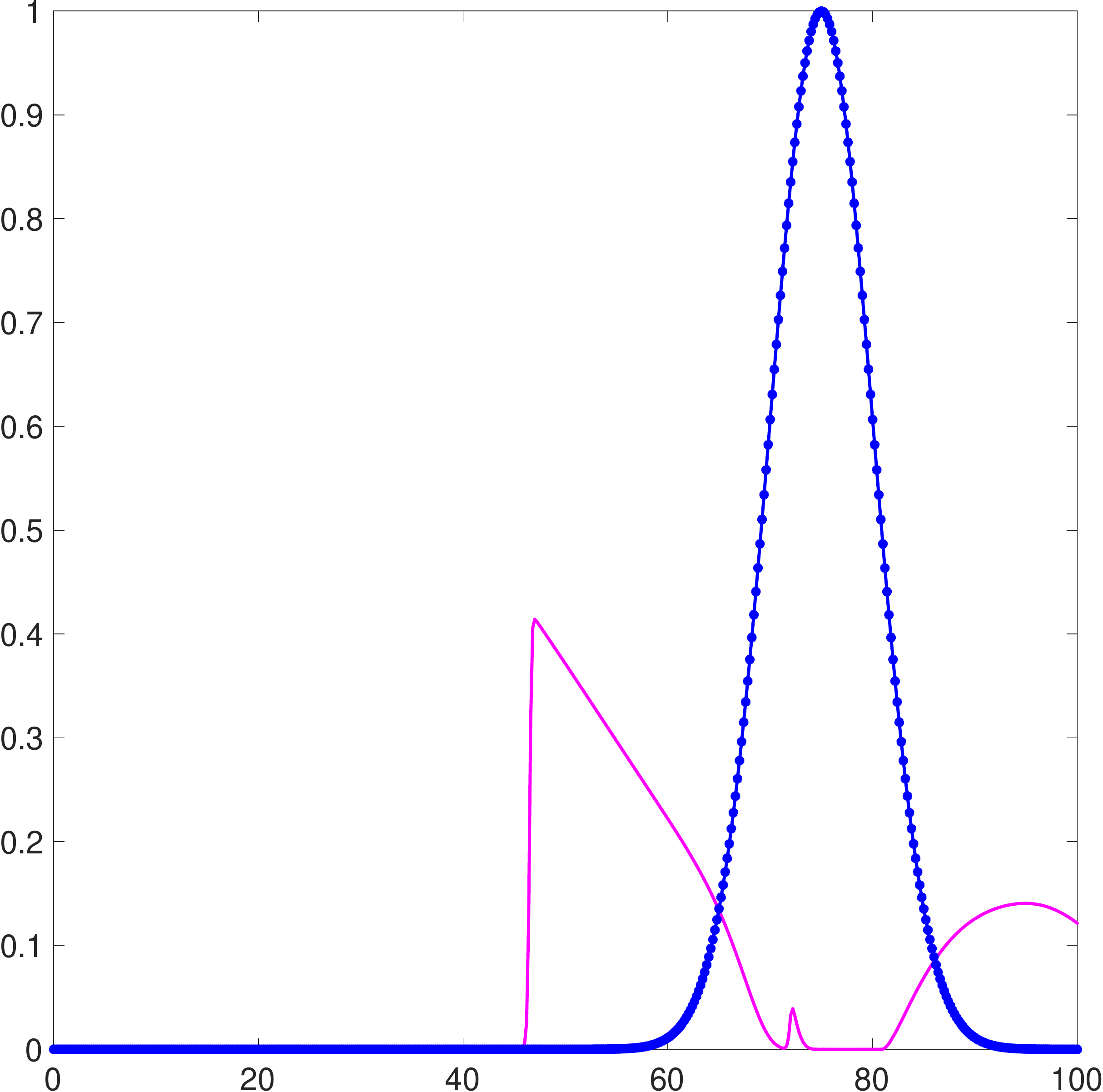}\label{fig:test3e}} \qquad
\subfloat[][Cars and truck density at time $t=T$  on Lane 2.]
{\includegraphics[width=.3\columnwidth]{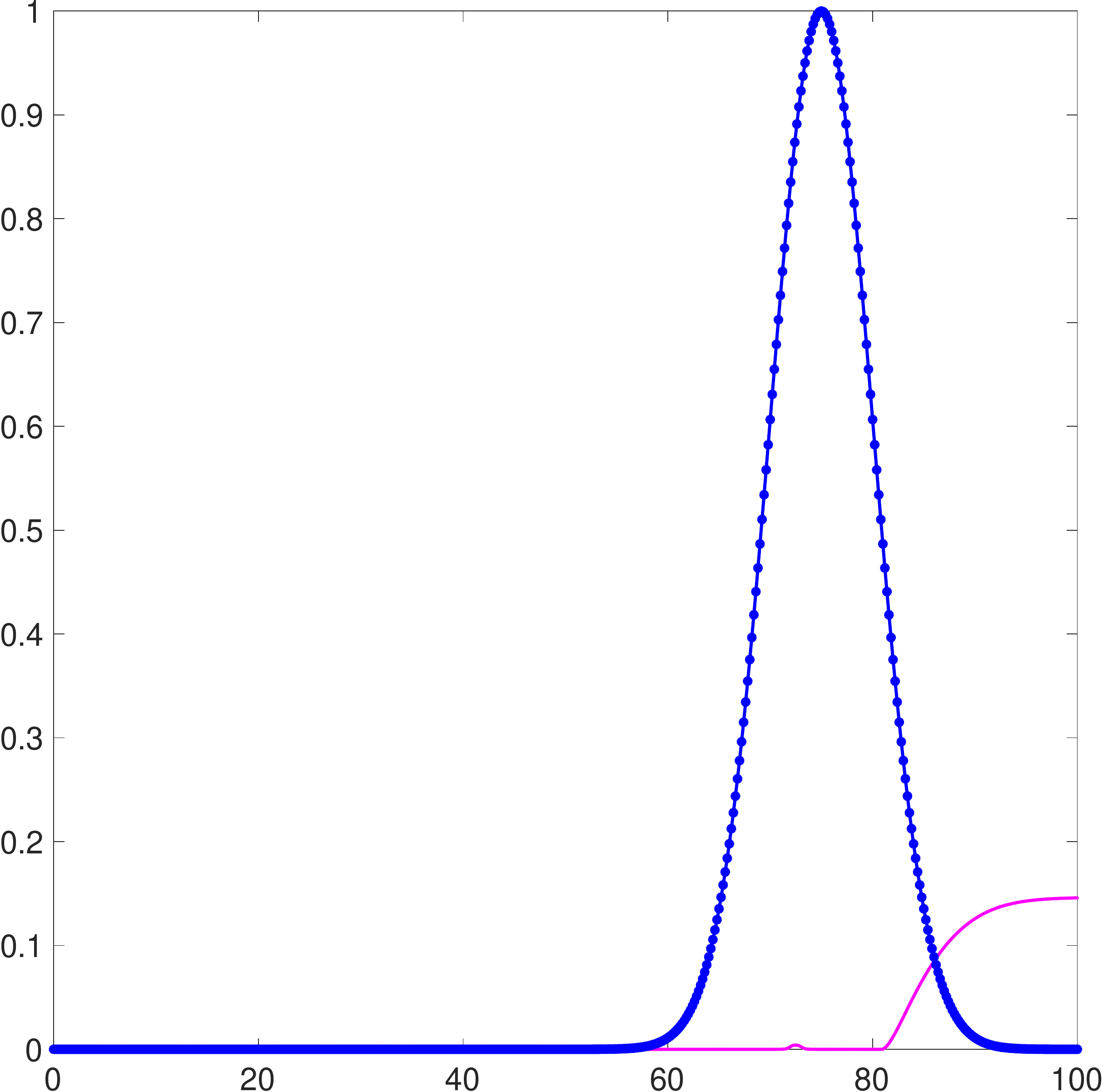}\label{fig:test3f}} 
\caption{Plot of the density of cars and truck on Lane 1 (first row) and Lane 2 (second row) at time $t=0$ (left), $t=T/4$ (middle) and $t=T$ (right).}
\label{fig:test3}
\end{figure}

In Figure \ref{fig:test3}, we show the density of cars and trucks on the two lanes at different times. The plots in the first row show Lane 1, with a car and no trucks, the plots in the second row show Lane 2, with a car and a truck which does not move at all during the simulation since its velocity is 0. The source terms in \eqref{eq:multilane} allow lane changing even if the traffic is not congested and vehicles are free to move along their lane. Hence, in Figure \subref*{fig:test3e} we see that the density of cars increases at the end of Lane 2, due to the source term, while in Figures \subref*{fig:test3b} and \subref*{fig:test3c} the density of cars related to Lane 1 increases since the truck does not move and the cars change lane. At the end of the simulation the density of cars is higher on Lane 1 while it is close to 0 on Lane 2.

From the two tests proposed in Sections \ref{sec:overtaking} and \ref{sec:multilane} we observe that the 2D multi-class model seems to be more suitable for capturing the overtaking of vehicles. Indeed, the two-dimensional description seem to fit better to such a dynamics, as shown in Figures \ref{fig:test2} and \ref{fig:test3}. 

\section{Conclusions}
In this work we have introduced a two-dimensional multi-class traffic model. We have analyzed the two-dimensional Riemann problems related to our model and provided numerical validations with a numerical scheme based on dimensional splitting. Then, we have analyzed the model with a dataset of real trajectories data, focusing on the dynamics of cars and trucks. The dataset has been used to calibrate the flux and velocity functions and to compare the numerical results with ground-truth data. The numerical tests have shown the good approximation of the trajectories with our model, obtaining a numerical error of $10^{-2}$. We have improved the results modifying the flux functions in order to consider different maximum velocity values for the two classes of vehicles. Finally, we have tested the ability of the model to simulate vehicles overtaking, also in comparison with a first order multi-class multi-lane model.  

Future investigations will aim at improving the proposed approach with second order models and deriving the corresponding microscopic model as a system of ordinary differential equations.

\section*{Acknowledgments}
Part of this work was carried out while C. Balzotti was visiting the University of Mannheim within the program IPID4all funded by the German Academic Exchange Service (DAAD).

\bibliography{biblio/biblio}  
\bibliographystyle{siam}

\end{document}